\PassOptionsToPackage{usenames,dvipsnames}{xcolor}
\documentclass[review,onefignum,onetabnum,usenames,dvipsnames]{article}
\makeatletter 
\providecommand*{\input@path}{}
\g@addto@macro\input@path{{../}}
\makeatother

\usepackage{array} % must be before arydshln
\usepackage{arydshln}
\usepackage{epsfig,graphicx}
\usepackage{makeidx}
\usepackage{multicol,multirow}
\usepackage{amsmath,amsfonts,amssymb,mathtools,dsfont,empheq}
\usepackage{amsthm}
\usepackage[normalem]{ulem}
\usepackage{crop,bbm}
\usepackage{xcolor,color}
\usepackage{todonotes}
\usepackage{tikz}
\usepackage{tikz-qtree}
\usetikzlibrary{shapes.geometric}
\usepackage{subfiles}
\usepackage{stmaryrd}
\usepackage{pdflscape}
\usepackage{breqn,dsfont,xspace}
\usepackage{rotating}
\usepackage{hyperref}
\usepackage{cleveref}
\usepackage{xr-hyper}
\usepackage{scalerel}
\usepackage{cancel,upgreek}
\usepackage{enumitem,cite}
\usepackage{subcaption}
\allowdisplaybreaks
\usepackage{rotating}
\usepackage{blkarray}

\usepackage{xparse}
\usepackage{physics}
\usepackage{forest}

\usepackage{array}
\newcolumntype{C}[1]{>{\centering\arraybackslash}m{#1}}  % Vertically-centered & horizontally-centered
\newcolumntype{L}[1]{>{\raggedleft\arraybackslash}m{#1}} % Vertically-centered & horizontally-to-left centered
\newcolumntype{R}[1]{>{\raggedleft\arraybackslash}m{#1}} % Vertically-centered & horizontally-to-right

\usepackage{arydshln}

\usepackage{stackengine}

\usepackage{fancyhdr}

\usetikzlibrary{calc,intersections}

\newif\ifdraft
\drafttrue

\newcommand{\nparts}{{\rm N}}      % Number of partitions for muti-physics systems

\newcommand{\Dt}{\Delta t}

\newcommand{\tree}{{\mathfrak{t}}}
\newcommand{\utree}{{\mathfrak{u}}}
\newcommand{\stree}{{\mathfrak{s}}}
\newcommand{\ktree}{{\mathfrak{k}}}
\newcommand{\rtree}{{\mathfrak{r}}}

\newcommand{\qtree}{{\mathfrak{q}}}
\newcommand{\vtree}{{\mathfrak{v}}}

\newcommand{\bushtree}{\mathfrak{b}}
\newcommand{\laddertree}{\mathfrak{l}}

\newcommand{\treecap}{{\mathrm{T}}}
\newcommand{\treevecinf}{{\vec{\mathrm{T}}}}
\newcommand{\treevecinfs}{{\vec{\mathrm{S}}}}

\newcommand{\xnode}{\mathpgoth{x}}

\newcommand{\unode}{u}
\newcommand{\vnode}{v}

\newcommand{\fforest}{\mathfrak{f}}
\newcommand{\gforest}{\mathfrak{g}}
\newcommand{\hforest}{\mathfrak{h}}
\newcommand{\dforest}{\mathfrak{d}}

\newcommand{\Elemtimes}{% element-wise multiplication
  \resizebox{!}{1.25\fontcharht\font`0}{$\mkern+1mu{\vartimes}\mkern+0mu$}
}
\newcommand{\elemtimes}{% element-wise multiplication
  \resizebox{!}{0.75\fontcharht\font`0}{$\mkern+2mu\vartimes\mkern+2mu$}
}

\forestset{
	btreeold/.style={
		for tree={
			grow=north,
			minimum width=3.5pt,
			minimum height=3.5pt,
			inner sep=0pt,
			outer sep=0pt,
			s sep=3pt,
			circle,
			fill,
			if level=0{
				baseline
			}{}
		},
		before computing xy={
			for tree={
				l=6pt,
			}
		},
		delay={
			for tree={
				if={strequal("r",content)}{rectangle,fill=lightred}{},
				if={strequal("x",content)}{fill=black,draw}{},
				if={strequal("z",content)}{fill=white,draw,thick}{},
				if={strequal("subtree",content)}{triangle,fill=black!20}{},
				content={}
			}
		}
	}
}
\forestset{
	btree/.style={
		for tree={
			grow=north,
			minimum width=4pt,
			minimum height=4pt,
			inner sep=0pt,
			outer sep=0pt,
			s sep=3pt,
			l=4pt,
			l sep=4pt,
			if n children=0{}{
				if level=0{
					baseline
				}{}
			},
			if n children=1{
			    if={mod(level(),2)==0}{
			        grow=60
			    }{
			        grow=120
			    }
			}{}
		},
		delay={
			for tree={
				if={strequal("b",content)}{rectangle,fill=cyan,content}{},
				if={strequal("r",content)}{rectangle,fill=lightred,content}{},
				if={strequal("x",content)}{circle,fill=black,draw,content}{},
				if={strequal("y",content)}{circle,fill=black,draw,content}{},
				if={strequal("z",content)}{circle,fill=white,draw,thick,content}{},
				if={strequal("s",content)}{star,fill=white,draw,thick,content}{},
				if={strequal("m",content)}{circle,fill=white,draw,content=$m$}{},
				if={strequal("n",content)}{circle,fill=white,draw,content=$n$}{},
				if={strequal("p",content)}{circle,fill=white,draw,content=$p$}{},
				if={strequal("k",content)}{circle,fill=white,draw,content=$k$}{},
				if={strequal("c0",content)}{circle,fill=white,draw,content=$0$}{},
				if={strequal("c1",content)}{circle,fill=white,draw,content=$1$}{},
				if={strequal("c2",content)}{circle,fill=white,draw,content=$2$}{},
				if={strequal("c3",content)}{circle,fill=white,draw,content=$3$}{},
				if={strequal("c4",content)}{circle,fill=white,draw,content=$4$}{},
				if={strequal("cy",content)}{circle,fill=white,draw,content=$\lambda$}{},
				if={strequal("cN",content)}{circle,fill=white,draw,content=$\nparts$}{},
				if={strequal("bm",content)}{regular polygon, regular polygon sides=4,fill=white,draw,content=$m$}{},
				if={strequal("bn",content)}{regular polygon, regular polygon sides=4,fill=white,draw,content=$n$}{},
				if={strequal("bp",content)}{regular polygon, regular polygon sides=4,fill=white,draw,content=$p$}{},
				if={strequal("bk",content)}{regular polygon, regular polygon sides=4,fill=white,draw,content=$k$}{},
				if={strequal("tm",content)}{regular polygon, regular polygon sides=3,fill=white,draw,inner sep=-0.1em,content=$m$}{},
				if={strequal("tn",content)}{regular polygon, regular polygon sides=3,fill=white,draw,inner sep=-0.1em,content=$n$}{},
				if={strequal("tp",content)}{regular polygon, regular polygon sides=3,fill=white,draw,inner sep=-0.1em,content=$p$}{},
				if={strequal("tk",content)}{regular polygon, regular polygon sides=3,fill=white,draw,inner sep=-0.1em,content=$k$}{},
				if={strequal("tk1",content)}{regular polygon, regular polygon sides=3,fill=white,draw,inner sep=-0.1em,content=$k_1$}{},
				if={strequal("tk2",content)}{regular polygon, regular polygon sides=3,fill=white,draw,inner sep=-0.1em,content=$k_2$}{},
				if={strequal("ttk",content)}{regular polygon, regular polygon sides=3,fill=white,draw,inner sep=-0.1em,content=$\ktree$}{},
				if={strequal("ttk1",content)}{regular polygon, regular polygon sides=3,fill=white,draw,inner sep=-0.1em,content=$\ktree_1$}{},
				if={strequal("ttk2",content)}{regular polygon, regular polygon sides=3,fill=white,draw,inner sep=-0.1em,content=$\ktree_2$}{},
				if={strequal("dm",content)}{diamond,fill=white,draw,content=$m$}{},
				if={strequal("dn",content)}{diamond,fill=white,draw,content=$n$}{},
				if={strequal("dp",content)}{diamond,fill=white,draw,content=$p$}{},
				if={strequal("dk",content)}{diamond,fill=white,draw,content=$k$}{},
				if={strequal("subtree",content)}{isosceles triangle,fill=black!20,draw,scale=1.5,yshift=0pt}{},
				if={strequal("t",content)}{circle split,draw,content}{},
				if={strequal("gth",content)}{regular polygon, regular polygon sides=3,fill=cyan,draw,content}{},
				if={strequal("gth1",content)}{regular polygon, regular polygon sides=3,fill=cyan,draw,content,label={[xshift=0.65ex,yshift=0.4ex]center:{\tiny $1$}}}{},
				if={strequal("gth2",content)}{regular polygon, regular polygon sides=3,fill=cyan,draw,content,label={[xshift=0.65ex,yshift=0.4ex]center:{\tiny $2$}}}{},
				if={strequal("gthr",content)}{regular polygon, regular polygon sides=3,fill=magenta,draw,content}{},
				if={strequal("gthr1",content)}{regular polygon, regular polygon sides=3,fill=magenta,draw,content,label={[xshift=0.65ex,yshift=0.4ex]center:{\tiny $1$}}}{},
				if={strequal("gthr2",content)}{regular polygon, regular polygon sides=3,fill=magenta,draw,content,label={[xshift=0.65ex,yshift=0.4ex]center:{\tiny $2$}}}{},
				if={strequal("gthg",content)}{regular polygon, regular polygon sides=3,fill=green,draw,content}{},
				if={strequal("gthg1",content)}{regular polygon, regular polygon sides=3,fill=green,draw,content,label={[xshift=0.65ex,yshift=0.4ex]center:{\tiny $1$}}}{},
				if={strequal("gthg2",content)}{regular polygon, regular polygon sides=3,fill=green,draw,content,label={[xshift=0.65ex,yshift=0.4ex]center:{\tiny $2$}}}{},
                if={strequal("w",content)}{circle,fill=red,draw,content}{},
                if={strequal("bc",content)}{circle,fill=blue,draw,content}{},
                if={strequal("mc",content)}{circle,fill=magenta,draw,content}{},
				if={strequal("e",content)}{cross out, line width = 0.25ex,draw = black,content}{},
				if={strequal("e1",content)}{cross out, line width = 0.25ex,draw = red,content}{},
				if={strequal("e2",content)}{cross out, line width = 0.25ex,draw = blue,content}{},
				if={strequal("ep",content)}{cross out, line width = 0.25ex,draw = magenta,content}{},
				if={strequal("d",content)}{diamond,fill=blue,draw,content}{}
			}
		}
	}
}
\newcommand{\btree}[1]{\Forest{btree #1}}
\forestset{
	btreelarge/.style={
		for tree={
			grow=north,
			minimum width=8pt,
			minimum height=8pt,
			inner sep=0pt,
			outer sep=0pt,
			s sep=5pt,
			l=7pt,
			l sep=7pt,
			if n children=0{}{
				if level=0{
					baseline
				}{}
			},
			if n children=1{
			    if={mod(level(),2)==0}{
			        grow=60
			    }{
			        grow=120
			    }
			}{}
		},
		delay={
			for tree={
				if={strequal("b",content)}{rectangle,fill=cyan,content}{},
				if={strequal("r",content)}{rectangle,fill=lightred,content}{},
				if={strequal("x",content)}{circle,fill=black,draw,content}{},
				if={strequal("y",content)}{circle,fill=black,draw,content}{},
				if={strequal("z",content)}{circle,fill=white,draw,thick,content}{},
				if={strequal("s",content)}{star,fill=white,draw,thick,content}{},
				if={strequal("m",content)}{circle,fill=white,draw,content=$m$}{},
				if={strequal("n",content)}{circle,fill=white,draw,content=$n$}{},
				if={strequal("p",content)}{circle,fill=white,draw,content=$p$}{},
				if={strequal("subtree",content)}{isosceles triangle,fill=black!20,draw,scale=1.5,yshift=0pt}{},
				if={strequal("t",content)}{circle split,draw,content}{},
				if={strequal("gth",content)}{regular polygon, regular polygon sides=3,fill=blue,draw,content}{},
                if={strequal("w",content)}{circle,fill=red,draw,content}{},
                if={strequal("bc",content)}{circle,fill=blue,draw,content}{},
                if={strequal("mc",content)}{circle,fill=magenta,draw,content}{}
                if={strequal("dots",content)}{text,\(\dots\)}{}
			}
		}
	}
}

\newcommand{\btreega}[1]{\!\begin{gathered}\btree{#1}\end{gathered}}

\newcommand{\Id}[1][]{% Identity matrix
   \ifthenelse{ \equal{#1}{} }
      {{\mathbf{I}}}
      {{\mathbf{I}_{#1 \times #1}}}
}
\newcommand{\ee}[1][]{% One column of the identity matrix
   \ifthenelse{ \equal{#1}{} }
      {{\mathbf{e}}}
      {{\mathbf{e}_{#1}}}
}
\newcommand{\one}[1][]{% Vector/matrix of ones
   \ifthenelse{ \equal{#1}{} }
      {{\mathbf{1}}}
      {{\mathbf{1}_{#1}}}
}
\newcommand{\zero}[1][]{% Vector/matrix of zeros
   \ifthenelse{ \equal{#1}{} }
      {{\mathbf{0}}}
      {{\mathbf{0}_{#1}}}
}

\newcommand{\sfrac}[2]{\mbox{\footnotesize$\displaystyle\frac{#1}{#2}$}} % small fraction for tables
\newcommand{\bigOh}{\mathcal{O}}

\newcommand{\impl}{{\scaleto{\mathbf{I}}{4pt}}}
\newcommand{\expl}{{\scaleto{\mathbf{E}}{4pt}}}
\NewDocumentCommand{\component}{m m m}{^{ \left\{ #1 \right\} #2 \IfBooleanT{#3}{\mathrm{T}} }}
\NewDocumentCommand{\comp}{m O{} s}{\component{#1}{#2}{#3}}
\NewDocumentCommand{\I}{O{} s}{\component{\impl}{#1}{#2}}
\NewDocumentCommand{\E}{O{} s}{\component{\expl}{#1}{#2}}
\NewDocumentCommand{\II}{O{} s}{\component{\impl, \impl}{#1}{#2}}
\NewDocumentCommand{\IE}{O{} s}{\component{\impl, \expl}{#1}{#2}}
\NewDocumentCommand{\EI}{O{} s}{\component{\expl, \impl}{#1}{#2}}
\NewDocumentCommand{\EE}{O{} s}{\component{\expl, \expl}{#1}{#2}}

\newcommand{\fun}{\mathbf{f}}

\newcommand{\gun}{\mathbf{g}}

\newcommand{\y}{\mathbf{y}}
\newcommand{\YY}   {\mathbf{Y}}

\newcommand{\Vertiii}{\vert\kern-0.22ex\vert\kern-0.22ex\vert}
\newcommand{\Vertiv}{\vert\kern-0.22ex\vert\kern-0.22ex\vert\kern-0.22ex\vert}



\definecolor{lightblue}{rgb}{0,0,1}
\definecolor{lightred}{rgb}{1,0,0}

\newcolumntype{T}{>{$}c<{$}} % math-mode version of "l" column type

\newcolumntype{V}{>{$}p{0.2cm}<{$}} % math-mode version of "l" column type
\newcolumntype{L}{>{$}p{2.25cm}<{$}} % math-mode version of "l" column type
\newcolumntype{Q}{>{$}p{11.8cm}<{$}} % math-mode version of "l" column type

\newcolumntype{H}{>{$}p{2.45cm}<{$}} % math-mode version of "l" column type
\newcolumntype{G}{>{$}p{11.60cm}<{$}} % math-mode version of "l" column type

\newcolumntype{U}{>{$}p{1.9cm}<{$}} % math-mode version of "l" column type
\newcolumntype{K}{>{$}p{12cm}<{$}} % math-mode version of "l" column type

\usepackage{macros_BSeries}
\usepackage[draft]{csl}

\usepackage{cases}
\usepackage{algpseudocode}
\usepackage{authblk}
\usepackage[numbers]{natbib}

\DeclareMathAlphabet{\mathpgoth}{OT1}{pgoth}{m}{n}

\numberwithin{equation}{subsection}

\usepackage{bm}
\usetikzlibrary{arrows.meta}
\usepackage {pgfplots}
\pgfplotsset{compat=1.17}
\newcommand\scalemath[2]{\scalebox{#1}{\mbox{\ensuremath{\displaystyle #2}}}}

\makeatletter 
\newtheorem{definition}{Definition}[section]
\crefname{definition}{Definition}{Definitions}
\newtheorem{remark}{Remark}[section]
\crefname{remark}{Remark}{Remarks}
\newtheorem{example}{Example}[section]
\newtheorem{theorem}{Theorem}[section]

\newtheorem{lemma}{Lemma}[section]
\newtheorem{property}{Property}[section]
\crefname{property}{Property}{Properties}
\newtheorem{corollary}{Corollary}[section]
\crefname{corollary}{Corollary}{Corollaries}

\crefname{conjecture}{Conjecture}{Conjectures}

\crefname{hypothesis}{Hypothesis}{Hypotheses}

\crefname{section}{Section}{Sections}
\crefname{lemma}{Lemma}{Lemmas}
\crefname{theorem}{Theorem}{Theorems}
\makeatother

\providecommand{\keywords}[1]{\small\textbf{\textit{Keywords---}} #1}

\newcommand{\ignore}[1]{}

\newcommand{\z}{\mathbf{z}}

\newcommand{\papertitle}{\BSeriesText Theory for Continuous Time Integration Methods}
\title{%
\papertitle%
\thanks{%
\funding{This work was funded by the National Science Foundation (NSF) award DMS-2411069 and by the Computational Science Laboratory.}
}}

\title{\papertitle}
\author[1]{Reid J. Gomillion}
\author[1]{Adrian Sandu}

\affil[1]{Department of Computer Science, Virginia Tech, 620 Drillfield Dr. Blacksburg, VA 24061. (rjg18@vt.edu, asandu7@vt.edu)}

\begin{document}

%CSL Report stuff

\cslyear{26}
\cslreportnumber{5}

\cslemail{rjg18@vt.edu, asandu7@vt.edu}

% Once the paper has been accepted, you can add the journal citation
% \cslciteas{citation}

\csltitlepage

\maketitle

\begin{abstract}
    \BSeriesText are a common tool for constructing numerical integration methods and partitioned variations. In this paper, we expand the operations of \BSeriesText to describe the formal expansion of modified differential equations found in multirate ODE solvers. Typically, such modified differential equations have time dependent approximations using summed fixed constant states. To create the \BSeriesText we introduce a new operation called replacement to induce the constant states and add a time parameterization to the labels of the \BSeriesText. By doing so we introduce new analysis via the permanent and block matrices which provide alternative formulations to composition and substitution. To support the time dependent labels properties are provided, the exact solution is described in Tree Iterated Integrals, and the underlying approximate representations are analyzed.
\end{abstract}

\tableofcontents

\section{Introduction}
\label{section:Introduction}

We consider the problem of a ODE with \(\nparts\) partitions given as:
\begin{equation} 
    \label{eq:PartitionedODE}
    \ddifffrac{t}\,\begin{bmatrix}
        \y\comp{1}\\
        \y\comp{2}\\
        \vdots \\
        \y\comp{\nparts}
    \end{bmatrix} =
    \begin{bmatrix}
        \fun\comp{1}\bigl(t,\y\comp{1},\y\comp{2}, \dots, \y\comp{\nparts}\bigr)  \\
        \fun\comp{2}\big(t,\y\comp{1},\y\comp{2}, \dots, \y\comp{\nparts}\big)  \\
        \vdots \\
        \fun\comp{\nparts}\big(t,\y\comp{1},\y\comp{2}, \dots, \y\comp{\nparts}\big)  \\
    \end{bmatrix}, \hspace{1cm} \begin{bmatrix}
        \y\comp{1}(t_\initstate)\\
        \y\comp{2}(t_\initstate) \\
        \vdots\\
        \y\comp{\nparts}(t_\initstate)
    \end{bmatrix} = \begin{bmatrix}
        \y\comp{1}_\initstate\\
        \y\comp{2}_\initstate \\
        \vdots\\
        \y\comp{\nparts}_\initstate
    \end{bmatrix},
\end{equation}
over the domain of \(t \in [t_\initstate, t_\finalstate]\) where \(\y\comp{i} \in \outspace^{n_{i}}\) and \(\fun\comp{i}: \outspace \times \prod^{\nparts}_{j=1} \outspace^{n_j} \rightarrow \outspace^{n_i}\). All RHS \(\fun\comp{i}\) are considered smooth w.r.t all variables. For this paper we aim to explain the approximation to \eqref{eq:PartitionedODE} via the approximate ODE 
\begin{equation}
    \label{eq:ApproximationProblem}
    \ddifffrac{\tau}\,\begin{bmatrix}
        \vg\comp{1}\\
        \vg\comp{2}\\
        \vdots \\
        \vg\comp{\nparts}
    \end{bmatrix} =
    \begin{bmatrix}
        \gun\comp{1}\big(\tau,\vg\comp{1},\vg\comp{2}, \dots, \vg\comp{\nparts}\big)  \\
        \gun\comp{2}\big(\tau,\vg\comp{1},\vg\comp{2}, \dots, \vg\comp{\nparts}\big)  \\
        \vdots \\
        \gun\comp{\nparts}\big(\tau,\vg\comp{1},\vg\comp{2}, \dots, \vg\comp{\nparts}\big)  \\
    \end{bmatrix}, \hspace{1cm} \begin{bmatrix}
        \vg\comp{1}(\tau_\initstate)\\
        \vg\comp{2}(\tau_\initstate) \\
        \vdots\\
        \vg\comp{\nparts}(\tau_\initstate)
    \end{bmatrix},
\end{equation}
over the time domain \(\tau \in [\tau_\initstate, \tau_\finalstate]\) where \(\gun\comp{i}: \outspace \times \times^{\nparts}_{j=1} \outspace^{n_j} \Rightarrow \outspace^{n_i}\). The functions \(\gun\comp{i}\) and states \(\vg\comp{i}\) are restricted to both having asymptotic expansions only in terms of the derivatives of \(\fun\comp{i}\). We describe such asymptotic expansions by \BSeriesText \cite{Butcher_2021_Bseries-book, Chartier_2010_Bseries, hairerOrderConditionsNumerical1981, butcherAlgebraicTheoryIntegration}. Such series are used to explain error and additional properties of numerical integration methods for ODEs. A more formal description of the asymptotic expansions for \cref{eq:ApproximationProblem} are provided in \cref{eq:modifiedAsymptotic} and this manuscript will build the necessary definitions for it.

The aim of this work is to construct operations for \BSeriesText to explain how \cref{eq:ApproximationProblem} approximates \cref{eq:PartitionedODE}. We allow for \textit{some} \(\gun\comp{i} = \fun\comp{i}\). Additionally, we deal with the cases that the system of equations \eqref{eq:ApproximationProblem} is either solved \textit{exactly} or \textit{approximately}. In most methods this approximate system will only be a single \textit{stage} of the scheme. Such stages are used to improve the approximations of the other stages to produce an improved result. Such stages will be denoted with subscripts \(\vg\comp{i}_{j}\) where \(j\) is some stage value.
An example of a method is the order \(2\) Multirate Infinitesimal method \cite{Sandu_2019_MRI-GARK} for \(\nparts = 2\) in component form
\begin{equation}
    \begin{split}
        \ddifffrac{\tau_{1}}\vg\comp{1}_1(\tau_{1}) &= \fun\comp{1}(\vg\comp{1}_1, \vg\comp{2}_1), \qquad \ddifffrac{\tau_{1}}\vg\comp{2}_1(\tau_{1}) = \fun\comp{2}(\y\comp{1}_{\initstate}, \y\comp{2}_{\initstate}),\\
        \vg\comp{1}_1(\tau_{1,\initstate}) &= \y\comp{1}_{\initstate}, \qquad \vg\comp{2}_1(\tau_{1,\initstate}) = \y\comp{2}_{\initstate}, \qquad \tau_{1} \in [t_{\initstate}, t_{\initstate} + \frac{1}{2}\,\Dt ],\\
        \ddifffrac{\tau_{2}}\vg\comp{1}_2(\tau_{2}) &= \fun\comp{1}(\vg\comp{1}_2, \vg\comp{2}_2), \qquad \ddifffrac{\tau_{2}}\vg\comp{2}_2(\tau_{2}) = -\fun\comp{2}(\y\comp{1}_{\initstate}, \y\comp{2}_{\initstate}) + 2\,\fun\comp{2}(\vg\comp{1}_1(\tau_{1,\finalstate}), \vg\comp{2}_1(\tau_{1,\finalstate})),\\
        \vg\comp{1}_2(\tau_{2,\initstate}) &= \vg\comp{1}_1(\tau_{1,\finalstate}), \qquad \vg\comp{2}_2(\tau_{1,\initstate}) = \vg\comp{2}_1(\tau_{1,\finalstate}), \qquad \tau_{2} \in [t_{\initstate} + \frac{1}{2}\,\Dt, t_{\finalstate}],\\
        \y\comp{1}_{\finalstate} &= \vg\comp{1}_{2}(\tau_{2,\finalstate}), \qquad \y\comp{2}_{\finalstate} = \vg\comp{2}_{2}(\tau_{2,\finalstate})
    \end{split}
\end{equation}
where the step size \(\Dt = t_\finalstate - t_\initstate\). Each stage has an approximate form however, each stage improves the approximation cumulatively for the final result. Another example is the Average Vector Field (AVF) \cite{quispelNewClassEnergypreserving2008,hairer2010energy} method given as 
\begin{equation}
    \begin{split}
        \ddifffrac{\tau}\vg\comp{1}(\tau) &= \fun(\vg\comp{2}) \qquad \ddifffrac{\tau}\vg\comp{2}(\tau) = \frac{\y_{\finalstate} - \y_{\initstate}}{\Dt}\\
        \vg\comp{1}(\tau) &= \y_{\initstate}, \qquad \vg\comp{2}(\tau) = \y_{\initstate}, \qquad \tau \in [t_{\initstate}, t_{\finalstate}],\\
        \y_{\finalstate} &= \vg\comp{1}(t_{\finalstate}).
    \end{split}
\end{equation}
Some numerical methods beyond these examples which are encapsulated by \eqref{eq:ApproximationProblem} are 
\begin{itemize}
    \item Multirate Infinitesimal methods (MRI/MIS) and variations \cite{Sandu_2019_MRI-GARK,guntherMultirateMethodsOrdinary2025,wenschMultirateInfinitesimalStep2009,robertsCoupledMultirateInfinitesimal2020, fish2024implicit,rice1960split}.
    \item Continuous Stage Runge-Kutta methods (CSRK)\cite{miyatakeCharacterizationEnergyPreservingMethods2016,hairer2010energy,quispelNewClassEnergypreserving2008,cai2018partitionedavf,eidnesOrderTheoryDiscrete2022} and by extension Continuous Runge-Kutta methods (CRK)\cite{owrenDerivationEfficientContinuous1992a}.
    \item Waveform relaxation \cite{Burrage_1995_book,whiteWAVEFORMRELAXATIONTHEORYa,bellenUseRungeKuttaFormulae1993} also sometimes referred to as Dynamic Iteration\footnote{From literature the term ``Dynamic Iteration'' was used originally to describe the linear differential equation splittings w.r.t. iterative linear solvers (\cite{miekkalaConvergenceDynamicIteration1987a}). However, nomenclature has evolved to use the same term in the nonlinear context just as Waveform Relaxation. In more recent works this can be seen. For simplicity for this work we only utilize the term ``Waveform Relaxation''.}\cite{miekkalaConvergenceDynamicIteration1987a,arnoldPreconditionedDynamicIteration2001a}.
    \item Exponential methods \cite{hochbruckExponentialIntegrators2010,Sandu_2020_Partitioned-Exp_Stiff,Sandu_2021_Partitioned-Exp_Classic}.
    \item Operator Splitting methods \cite{blanes2024splitting,Mclachlan_2002_splitting,Blanes_2010_splitting,Sandu_2021_GLM-ADI,Sandu_2022_GARK_splittingSchemes}.
    \item Continuous Co-Simulation schemes of extrapolating/interpolating coupling inputs to differential equations \cite{gomes2018co, buschContinuousApproximationTechniques2016,friedrich2011parallel}.
    \item Generalized Additive Runge-Kutta (GARK)\cite{Sandu_2015_GARK,Sandu_2016_MR-GARK} methods.
    \item Modified Differential Equations\footnote{In multirate literature (i.e. \cite{guntherMultirateMethodsOrdinary2025}) the same term of ``Modified Differential Equations'' is used to describe the sub-problems solved with approximate couplings. However, such modified differential equations in multirate literature differ then the ones from geometric integration. In the context of this work the asymptotic series encapsulates both variations thus, we do not differentiate between such.} \cite{Hairer_2006_geometric-book,chartierNumericalIntegratorsBased2007}.
\end{itemize}
This is not an exhaustive list of all possible methods encapsulated by \eqref{eq:ApproximationProblem}. We take extra care to focus on the case of MRI methods.

In \cref{section:Background} we start by introducing trees, colors, operations on them, sub-sets, and their usage in \BSeriesText. We introduce a new node type called ``stumps'' and their representations. Next in \cref{section:Combinatorics} we define generalized outer products and the matrix permanent \cite{minc1984permanents} which are not typically seen in \BSeriesText analysis. With the foundations established, \cref{section:Composition} provides the definition of \BSeriesText composition \cite{Chartier_2010_Bseries,Butcher_2021_Bseries-book,Hairer_2006_geometric-book}, additional properties, and the repeated matrix formulation. Next in \cref{section:Replacement} a new operation called replacement is provided with discussions on its similarity to \BSeriesText substitution \cite{Chartier_2010_Bseries,Hairer_2006_geometric-book}. \cref{section:TFDEQ} defines the time dependent labels and their underlying representations. Lastly, in \cref{section:TII} we introduce the exact solution of time dependent series via Tree Iterated Integrals \cite{kreimer1999chen}. Throughout the sections recursive functions for these operations will be provided to allow for easy of implementation. Colored trees and some of these operations are provided in \cite{roberts2025integreat}.

With each operation defined for transforming \eqref{eq:ApproximationProblem} into a single \BSeriesText we provide some examples of methods explained via this analysis. The authors note that TII \cite{kreimer1999chen}, composition, substitution and problems of the form \eqref{eq:ApproximationProblem} can be found in other context such as Hopf algebra \cite{calaque2011two,Chartier_2010_Bseries,chanRelationsHopfAlgebra2014}, Lie series \cite{muruaHopfAlgebraRooted2006,muruaWordSeriesDynamical2017}, quantum field theory \cite{brouder2000runge,brouderTreesRenormalizationDifferential2004,connesHopfAlgebrasRenormalization1998}, and rough differential equation analysis \cite{gubinelli2010ramification}. This is not to be an extensive list of connecting fields of research. For this work we take inspiration from them but focus on the numerical methods side and a more combinatorics route to their derivations.

Some proofs are reserved for the appendix. Additional supporting triangular block matrix analysis is provided in the appendix as well. This is to support \BSeriesText composition, substitution, coefficient properties, and various block formulations.

Contributions of this text are as follows:
\begin{itemize}
    \item Continuous parameterization of \BSeriesText with their underlying approximate expansions.
    \item Generalizations of binary operations on trees which preserve symmetry.
    \item Formulations of composition and substitution which do \textbf{not} rely on ordered trees. All new operations will be defined on unordered trees.
    \item New replacement operation which allows for constant vector field transformations.
    \item Matrix form of analysis for composition and backward error analysis.
    \item Exact solution and partial constructed from Tree Iterated Integrals for \BSeriesText and how they relate to the underlying series forms.
    \item Example applications of the theory applied to continuous methods.
\end{itemize}
\begin{remark}
    Note that though the focus has been on component partitioning \cref{eq:PartitionedODE} all results are applicable to the additive partitioning \cite{guntherMultirateMethodsOrdinary2025}. 
\end{remark}

\section{Trees and \BSeriesText}
\label{section:Background}

For this section we introduce trees and their usage in \BSeriesText analysis\cite{Butcher_2021_Bseries-book, Chartier_2010_Bseries, hairerOrderConditionsNumerical1981, butcherAlgebraicTheoryIntegration, Hairer_2006_geometric-book}. Such analysis is common when designing numerical integrators for ODEs. There will be a further emphasis on multicolored trees \cite{araujo1997symplectic, muruaFormalSeriesNumerical1999, murua1995metodos, hairerOrderConditionsNumerical1981} for their usage in \NBSeriesText and \PSeriesText for partitioned numerical methods. As such each sub-system in \cref{eq:PartitionedODE} is represented as a different color thus, such series are of the main focus.

Definitions to trees, functions for trees, and series are provided. New material relating to ``stump'' nodes are provided for explaining the approximate sub-system functions \cref{eq:ApproximationProblem}. 

\subsection{Tree Spaces}
\begin{definition}[Tree]\label{def:tree}
    A tree is an object with nodes and edges to form a graph such that there are no cycles. All trees have a single node specified as the root. Visually this is shown as the bottom of the tree. For the rest of the paper only \textbf{unordered} trees are considered, unless specified.\par
    An unordered tree is defined as given some node with children then all possible permutations of children are equal. As an example, consider 
    \[\btree{[x[x[x]][x]]} = \btree{[x[x][x[x]]]}.\]
    The set of all possible unordered rooted trees is \(\tree \in \TT\).
    A tree may have different \textbf{colors} at each node. As an example
    \[\{\btree{[x[x]]} ,\btree{[x[z]]} ,\btree{[z[x]]} ,\btree{[z[z]]}\}.\]
    The set of all possible unordered rooted multicolored trees is \(\tree \in \NTTC\) with colors \(\Tcolors\). 
    If it doesn't matter whether the tree is colored or not then as shorthand \(\tree \in \AbstractT\) is used.
\end{definition}
More formal definitions on graphs are defined in \cite{Butcher_2021_Bseries-book,butcherAlgebraicTheoryIntegration}.
Next the notion of a forest (discrete set of trees) will be defined.
\begin{definition}[Forest of Trees \cite{Butcher_2021_Bseries-book}]\label{def:forest}
    A "Forest" is described as a set of trees where order of the trees does not matter. Typically, this is done through multiplication where product of the trees is commutative. The forests are represented as 
    \[\underbrace{\tree_1, \dots \tree_1}_{q_1},\underbrace{\tree_2, \dots \tree_2}_{q_2}, \dots, \underbrace{\tree_n, \dots \tree_n}_{q_n} = \forestTQt \in \ForestSpace{\AbstractT}\]
    where $\tree_i \in \AbstractT\, \forall i$, $\tree_i \ne \tree_j  \text{ for } i \ne j$, $q_i \in \mathbb{Z}^+$ and \(\tree_i \in \AbstractT, \forall i\). \(\ForestSpace{\AbstractT}\) is the forest over the space of trees \(\AbstractT\). A more restricted form is given with \(\ForestSpace{\AbstractT}_{p}\) where all forests must satisfy the constraint \(\sum_{i=1}^{n} q_i = p\).
\end{definition}
This type of notation allows for a dense representation of trees and their children. Next the tree sets are recursively defined.
\begin{definition}[Multi-Colored Trees]\label{def:coloredTreesSet}
    Let $c \in \Tcolors$ be a set of colors which are unique and \(\tau_c \in \NTTC\) are the single nodes with some color. Then the colored trees are formally given as:
    \[\tau_c[\tree_1^{q_1} \tree_2^{q_2} \dots \tree_n^{q_n}] \in \NTTC,\]
    where $\tree_1^{q_1} \tree_2^{q_2} \dots \tree_n^{q_n} \in \ForestSpace{\NTTC}$.
\end{definition}
Such a definition can be used recursively for the children of each node.
\begin{remark}[Root Notation]
    Throughout the paper some common denotations for single node trees are \(\tau\), \(\tau_c\), and \(\psi\). All trees are specified with \(\tau_{c}[\fforest]\) notation where \(\fforest\) is a forest. This is different from standard notations of colored trees where bracket wise subscripts for the color \([\fforest]_{c}\) or \(\prescript{}{c}{[\fforest]}\) are used instead.
\end{remark}
In literature there is a special tree that is used to denote the empty tree.
\begin{definition}[Empty Tree (Ground)]\label{def:emptytree}
    The empty tree is denoted as \(\emptyset\) or also referred to as the ground. The node has the following properties:
    \begin{equation}
        \begin{split}
            \tree\, \emptyset &= \tree,\\
            \tau[\emptyset] &= \tau,
        \end{split}
    \end{equation}
    where \(\tree,\tau \in \AbstractT\) and \(\tau\) is a single node tree. Thus, trees can have children of \(\emptyset\), but they resolve to the tree without them. 

    Whether \(\emptyset \in \NTTC\) is inconsistent in literature; however, for the purposes here it is assumed to be. Note that at times there can be empty trees with a color associated with them. Given the colors \(c \in \Tcolors\) then a colored ground is \(\emptyset_c\) which has the same properties. Such distinctions help in PSeries \cite{Hairer_2006_geometric-book}; however, there are transformations to remove these terms. We focus on the non-colored one but will mention the colored variant on occasion.
\end{definition}
A common representation that will be utilized throughout is given as follows.
\begin{definition}[Forest Space \cite{butcherTreeForestSpaces2010,butcherNewApproachAlgebraic,chanRelationsHopfAlgebra2014,Butcher_2021_Bseries-book}]
    % [Section 2.4]
    A ``Forest Space'' is described as a linear span of forest. Formally as 
    \begin{equation}
        \sum_{i=1}^{m} a_i\, \fforest_i = \sum_{i=1}^{m} a_i \left[\prod_{j=1}^{n} \tree_{i,j}^{q_{i,j}}\right] \in \ForestSpace{\AbstractT}
    \end{equation}
    where \(a_i \in \outspace\), \(\fforest_i \in \ForestSpace{\AbstractT}\), and \(\tree_{i,j} \in \AbstractT\).
\end{definition}
\begin{definition}[Zero or Null Tree]
    With the ``Forest Space'' there will be an introduction of the ``null tree''. The tree has no colors and is represented by \(0\) with the properties:
    \begin{equation}
        \begin{split}
            a \, 0 &= 0,\\
            0 \, \tree &= 0,\\
            0 + \tree &= \tree,\\
            \tau[0] &= 0,
        \end{split}
    \end{equation}
    where \(a \in \outspace\), \(\tree,\tau \in \AbstractT\) and \(\tau\) is a single node tree. The zero tree is not in any set of trees.
\end{definition}
\begin{remark}[Trees/Forest on Functions]\label{remark:treeFunc}
    Let \(\fun\) be some function that takes a tree and/or forest as input and outputs something. Unless specified otherwise the function must always respect the unordered properties of every tree. As an example 
    \[\fun\left(\btreega{[x[x[x]][x]]}\right) = \fun\left(\btreega{[x[x][x[x]]]}\right).\]
    If the function does not have this property then it is not valid over the set of unordered trees but instead ordered trees. For this paper we actively avoid such operations such that only unordered trees are used.
\end{remark}
\begin{definition}[Tree Leaf]\label{def:leaf}
    Tree leaves are defined as nodes that have no children. Furthermore, there is always a path starting at the root of the tree \(\tree \in \AbstractT\) that must end at a leaf.
\end{definition}
\begin{definition}[Ladder (Tall) Trees]\label{def:talltree}
    A ladder or tall tree is defined as a tree in which every node can only have \(1\) or \(0\) children. This is independent of color. Some examples are as follows:
    \begin{equation}
        \left\{\btreega{[x]}, \btreega{[x[x]]}, \btreega{[x[x[x]]]}, \btreega{[x[x[x[x]]]]}, \dots\right\} = \LadderTrees{\TT}.
    \end{equation}
    Furthermore, there is always a path starting at the root of some arbitrary tree \(\tree \in \AbstractT\) that must end at a subtree which has the structure of a ladder tree.
\end{definition}
\begin{definition}[Tall Bush Trees]\label{def:bushtree}
    A tall bush tree is defined as only the root node can have more than \(1\) child, but all other nodes can only have \(1\) or \(0\) children. Thus, only ladder trees can be children of the root. Some examples are as follows:
    \begin{equation}
        \left\{\btreega{[x]}, \btreega{[x[x]]}, \btreega{[x[x[x]]]}, \btreega{[x[x[x[x]]]]}, \dots, \btreega{[x[x][x]]}, \btreega{[x[x][x[x]]]}, \btreega{[x[x][x[x[x]]]]}, \dots, \btreega{[x[x][x][x]]}, \btreega{[x[x][x][x[x]]]}, \dots \right\} = \TallBushTrees{\TT}.
    \end{equation}
    A tree from the set has the representation of \(\bushtree = \tau[\forestLadder]\) where \(\laddertree_i \in \LadderTrees{\TT}\).
    Furthermore, there is always a path starting at the root of some arbitrary tree \(\tree \in \AbstractT\) that must end at a subtree which has the structure of a tall bush tree.
\end{definition}
 Note that typically tall bush trees are given as tree with only single nodes as children of the root. This typically refers to only bush trees instead of tall bush trees. We push this notion to ladder trees as this introduces more complicated symmetry factors versus the simple case.
 \begin{definition}[\(n\)th Tall Bush Trees]\label{def:nthbushtree}
    An \(n\)th tall bush tree is defined as only nodes at a distance less than or equal to \(n\) of the root node can have more than \(1\) child. Distance being the number of edges needed to cross along the path from the current node to the root node. Otherwise, nodes which are at a distance greater than \(n\) can only have \(1\) or \(0\) children. The case when \(n = 0\) are tall bush trees and \(n = -1\) are ladder trees. An example is as follows:
    \begin{equation}
        \left\{\btreega{[x[x[x][x]]]}, \btreega{[x[x[x][x]][x]]}, \dots \right\} = \TallBushTreesN{\TT}{1} \setminus \TallBushTreesN{\TT}{0}.
    \end{equation}
    A tree from the set provides \(\tree = \tau[\forestTQt] \in \TallBushTreesN{\TT}{n}\) where \(\tree_i \in \TallBushTreesN{\TT}{n-1}\).
    Furthermore, there is always a path starting at the root of some arbitrary tree \(\tree \in \AbstractT\) that must end at a subtree that is an \(n\)th tall bush tree.
 \end{definition}
Each of these sub-sets demonstrate an increasing interaction of symmetry onto the trees. These sub-sets of trees will be utilized throughout and for proofs.

\subsection{Operations on Trees}
We provide some common functions for trees that are utilized throughout.
\begin{definition}[Tree Cardinality]\label{def:Cardinality}
    The cardinality of a tree \(|\tree|\) is defined as the sum of all nodes in the tree \(\tree \in \AbstractT\). The empty tree has the cardinality \(|\emptyset| = 0\).
\end{definition}
For the case of trees in this article the cardinality is the order.
\begin{definition}[Tree Order]\label{def:Order}
    The order of a tree \(\tree \in \NTTC\) is defined as its cardinality \(\treeorder{\tree}\).
\end{definition}
Unordered trees have a notion of symmetry associated with them.
\begin{definition}[Tree Symmetry]\label{def:Sym}
    Tree symmetry is defined as 
    \begin{equation}
        \begin{split}
            \sigma(\tree) = \sigma(\tau[\tree_1^{q_1} \tree_2^{q_2} \dots \tree_n^{q_n}]) &= \prod_{i=1}^n q_i!\, \sigma(\tree_i)^{q_i},\\
            \sigma(\tau) &= 1,\\
            \sigma(\emptyset) &= 1,
        \end{split}
    \end{equation}
    where $t = \tau[\tree_1^{q_1} \tree_2^{q_2} \dots \tree_n^{q_n}] \in \NTTC $.
\end{definition}
The gamma (or sometimes referred to as density or tree factorial  \cite{Butcher_2016_book}) of the tree is given.
\begin{definition}[Tree Density]\label{def:Gamma}
    Tree density is defined as
    \begin{equation}
        \begin{split}
            \gamma(\tree) = \tree! = \gamma(\tau[\tree_1^{q_1} \tree_2^{q_2} \dots \tree_n^{q_n}]) =  \treeorder{\tree}\, \prod_{i=1}^n (\gamma(\tree_i))^{q_i},
        \end{split}
    \end{equation}
    where $\tree = \tau[\tree_1^{q_1} \tree_2^{q_2} \dots \tree_n^{q_n}] \in \NTTC $. Additionally, given a set of colors \(\Tcolors\) let the restricted form be defined as
    \begin{equation}
        \gammaFun{\Tcolors}(\tau[\tree_1^{q_1} \tree_2^{q_2} \dots \tree_n^{q_n}]) = \begin{cases}
                \treeorder{\tree} \, \prod_{i=1}^n (\gamma(\tree_i))^{q_i} & \textbf{ if } \tau \in \Tcolors,\\
                \infty & \textbf{ otherwise}.
        \end{cases}
    \end{equation}
    When taking the inverse of this restricted form \(1/\gammaFun{\Tcolors}\) then nodes not in the set of colors resolve the inverse to \(0\).
\end{definition}
\begin{definition}[Monotonic Labels 
    \cite{Butcher_2021_Bseries-book,hairerOrderConditionsNumerical1981}]\label{def:monotonic}
    The number of monotonic labelings of a tree is defined as follows:
    \begin{equation}
        \alpha(\tree) = \frac{\treeorder{\tree}!}{\gamma(\tree)\,\sigma(\tree)}.
    \end{equation}
    Furthermore, when observing the forest case the following holds
    \begin{equation}
        \alpha(\tau[\forestTQt]) = (\treeorder{\tree}-1)! \prod_{i=1}^{n}\frac{1}{q_i ! \gamma(\tree_i)^{q_i} \sigma(\tree_i)^{q_i}}.
    \end{equation}
\end{definition}
\begin{remark}[Alternative Proof]
    Note if one has a tree \(\bushtree \in \TallBushTreesN{\NTTC}{0}\) then the expression for \(\alpha\) matches \cite{charalambidesEnumerativeCombinatorics2002}[Theorem 2.10]. More explicitly given \\\(\bushtree = \tau[\forestLadder] \in \TallBushTreesN{\NTTC}{0}\) where \(\laddertree_i \in \LadderTrees{\NTTC}\) then only linear increases in monotonic labels are possible thus giving
    \[\frac{\treeorder{\bushtree}}{\treeorder{\bushtree}} \left((\treeorder{\bushtree}-1)! \prod_{i=1}^{n}\frac{1}{q_i ! (\treeorder{\laddertree_i}!)^{q_i}} \right) = \alpha(\bushtree)\]
    since the \(\alpha(\laddertree_i) = \sigma(\laddertree_i) = 1\). Continuing the iteration by inductively applying over sets \(\TallBushTreesN{\NTTC}{n}\) gives the final result.
\end{remark}

A common tool used in the analysis of \BSeriesText is the ``Butcher Product'' or in other terms the ``beta-product''.
\begin{definition}[Butcher Product (Beta-Product)]\label{def:butcherprod}
    The Butcher Product is defined as a (non-associative) binary operation \cite{Butcher_2021_Bseries-book,Hairer_2006_geometric-book} on trees as \(\ButcherProduct: \AbstractT \times \AbstractT \rightarrow \AbstractT\). The product is given as
    \begin{equation}
        \begin{split}
            \tree \ButcherProduct \stree &= \tau[\forestTQt] \ButcherProduct \stree = \tau[\forestTQt \stree],\\
            \tau[\forestTQt] \ButcherProduct \forestTQs &= \tau[\forestTQt \forestTQs]
        \end{split}
    \end{equation}
    where \(\tree, \stree, \in \AbstractT\). 
\end{definition}

Additionally, a function is provided to grab the color of the root node of a tree.
\begin{definition}[Root Color of Tree] \label{def:treeroot}
    For a tree \(\tree\) the operation to get its root's color is
    \begin{equation}
        \treeRoot(\tree) = \treeRoot(\tau_{c}[\forestTQt]) = c.
    \end{equation}
\end{definition}

The use of tree based analysis in differential equations is closely tied to the differential operator. The idea that an ordinary differential equation expansion can be represented by the expansion of trees.
\begin{definition}[Differential Operator on Trees\cite{hairerOrderConditionsNumerical1981}]\label{def:diffop}
    For a tree \\\(\tree = \tau[\forestTQt] \in \NTTC\) the differential operator it represents is defined as 
    \begin{equation}
        \scalemath{1.0}{
            \begin{split}
                \BDiffOp(\tree)(\y) &= \fun\comp{c}_{(\y\comp{c_{1}})^{q_1} \dots (\y\comp{c_{n}})^{q_n}}(\y) (\underbrace{\BDiffOp(t_1)(\y), \dots, \BDiffOp(\tree_1)(\y)}_{q_1}\\
                &, \underbrace{\BDiffOp(\tree_2)(\y), \dots, \BDiffOp(\tree_2)(\y)}_{q_2}, \dots, \underbrace{\BDiffOp(\tree_n)(\y), \dots, \BDiffOp(\tree_n)(\y)}_{q_n}),
            \end{split}
        }
    \end{equation}
    where the operator is evaluated at some state \(\y \in \outspace^{p}\) and \(c_{i}\) are the colors of the subtrees. The special case of \(\emptyset_c\), where \(c \in \Tcolors\), returns
    \begin{equation}
        \BDiffOp(\emptyset_c)(\y) = \y_c
    \end{equation}
    where \(\y_c\) is the partitioning of \(\y\) that gives the color \(c\). More specifically \(\y_c \in \outspace^{p_c}\) and \(\y \in \outspace^{p}\) where \(p = \sum_{c \in \Tcolors} p_c\).
\end{definition}

\subsection{Stumps}
Stumps are a new addition (for vector fields) to the tree operations. We take inspiration from the works \cite{butcherTreesStumpsApplications2018,Butcher_2021_Bseries-book} by extending to multiple variables and partitions. Going forward we use the node \(\xnode\) to denote a stump node. Such a stump is considered colorless and when colors are required they are given as \(\xnode_{c}\). Introducing colors allows restrictions on operators due to partitioning and the colorless stump acts as catch-alls when colors do not matter.

Some set notations of stumps are provided as follows.
\begin{definition}[Stump Grafting]\label{def:stumpTreeGraft}
    Given a tree \(\tree \in \NTTC\) and a set of stumps colors defined as \(\TStumpcolors\) then for \(\dforest \in \ForestSpace{\TStumpcolors}\) the set \(\tree\lceil \dforest \rceil\) is the set of all possible trees \(\tree\) with the stumps \(\dforest\) grafted on to arbitrary nodes. The set of trees \(\tree\lceil \ForestSpace{\TStumpcolors} \rceil\) is the set of all possible stump forest grafted sets (including the empty set). The set \(\NTTC\lceil \dforest \rceil\) is the set union of all possible trees with the stump forest \(\dforest\) grafted on. The empty set can only be grafted with one stump node. Lastly \(\NTTC\lceil \ForestSpace{\TStumpcolors} \rceil\) is the combination of both and \(\ForestSpace{\NTTC}\lceil \ForestSpace{\TStumpcolors} \rceil\) extends this over forest.
\end{definition}
To count the amount of stumps in the trees the following functions are provided.
\begin{definition}[Amount of Stumps in a Tree]\label{def:stumpAmount}
    Given a tree \(\tree \in  \NTTC\lceil \ForestSpace{\TStumpcolors} \rceil\) then the total amount of stumps in the tree is given as \(\stumpCount(\tree)\). For the amount for a specific stump \(\xnode\) this is given as \(\stumpCountRes{\xnode}(\tree)\). As shorthand for factorial products \(\stumpCountFactorial(\tree) = \prod_{\xnode \in \TStumpcolors} \stumpCountRes{\xnode}(\tree)!\). For shorthand for converting forest to colored stump of the same amount of base nodes \(\stumpForestConvert(\forestTQt) = \prod_{j=1}^{n} \xnode_{\treeRoot(\tree_{j})}^{q_j} = \prod_{c \in \Tcolors} \xnode_{c}^{x_c}\) where \(x_{c} = \sum_{j=1}^{n} \delta_{\treeRoot(\tree_{j}),c} q_j\).
\end{definition}
Before going into the effects on the elementary differentials we provide a symmetric product over tensors for our intentions over trees.
\begin{definition}[Symmetric Tree Product]\label{def:symProduct}
    Given a tensor \(\mathbf{A} \in (\outspace^{n})^{\otimes (p+1)}\) where \(p \ge 1\) then for some vector \(\mathbf{v} \in \outspace^{n}\) the symmetric product is
    \begin{equation}
        (\mathbf{A} \symprod \mathbf{v})_{i} = (\mathbf{A} \mathbf{v})_{i} = \frac{1}{p}\sum_{k=1}^{p} \sum_{j_{k}=1}^{n} \mathbf{A}_{i,j_1, \dots, j_p}\mathbf{v}_{j_{k}} = \frac{1}{p}\left(\sum_{k=1}^{p} (\mathbf{A}(\underbrace{\dots}_{k-1}, \mathbf{v}, \underbrace{\dots}_{p - k}))_{i}\right).
    \end{equation}
    If another vector \(\mathbf{w} \in \outspace^{n}\) is applied at the same time then \((\mathbf{A} \symprod \mathbf{v}) \symprod \mathbf{w} = \mathbf{A} \symprod (\mathbf{v} \mathbf{w})\) where the product of vectors \(\mathbf{v} \mathbf{w}\) are treated as intermediates until they are applied to a tensor. As shorthand for vector inputs \(\mathbf{A} (\mathbf{v} \mathbf{w})\) can be used instead.
    If instead \(\mathbf{A} \in \prod^{m}_{j=1} (\outspace^{n_j})^{\otimes p_i}\) then given a vector \(\mathbf{v} \in \outspace^{n_q}\) the symmetric tree product is 
    \begin{equation}
        (\mathbf{A} \symprod \mathbf{v})_{i} = (\mathbf{A} \mathbf{v})_{i} = \frac{1}{p_{q}}\sum_{k=1}^{p_q} \sum_{j_{k,p_{q}} = 1}^{n_{q}} \mathbf{A}_{i, j_{(1,1)}, \dots, j_{(m,p_{m})}}\mathbf{v}_{j_{k}} = \frac{1}{p_q}\left(\sum_{j=1}^{p_q} \mathbf{A}(\dots,\underbrace{\dots}_{j-1}, \mathbf{v}, \underbrace{\dots}_{p - j},\dots)\right)_{i}.
    \end{equation}
    In the case the dimensions are the same, but the contractions are treated differentially then various dimensions are still treated differentially and an index \(a\) is provided to select them as \(\symprod_{a}\) where \(1 \le a \le m\).

    For when tensor on tensor contractions are needed then if instead a tensor \(\mathbf{B} : (\outspace^{n})^{\otimes q}\) and \(q \ge 1\) the symmetric product is 
    \begin{equation}
        (\mathbf{A} \symprod \mathbf{B})_{i} = \frac{1}{p}\sum_{k=1}^{p} \sum_{j_{k}=1}^{n} \mathbf{A}_{i,j_1, \dots, j_p}\mathbf{B}_{j_{k}, \hcancel{l}_1, \dots, \hcancel{l}_{q}} = \frac{1}{p}\left(\sum_{i=1}^{p} \mathbf{A}(\underbrace{\dots}_{i-1}, \mathbf{B}, \underbrace{\dots}_{p - i})\right).
    \end{equation}
    In the case of contractions over tensor inputs then the simplified representation of vector inputs \(\mathbf{A} \mathbf{B}\) no longer holds. For two tensors to prevent contractions on top of another tensor \(\mathbf{C}\) is represented as
    \begin{equation}
        \begin{split}
            (\mathbf{A} \symprod (\mathbf{B} \circledcirc \mathbf{C}))_{i} &= \frac{1}{p(p-1)}\sum_{k=1}^{p} \sum_{a=1}^{p-1}\sum_{j_{a}=1}^{n}\sum_{j_{k}=1}^{n} \mathbf{A}_{i,j_1, \dots, j_p}\mathbf{B}_{j_{k}, \hcancel{l}_1, \dots, \hcancel{l}_{q_1}}\mathbf{C}_{j_{a}, \hcancel{i}_1, \dots, \hcancel{i}_{q_2}} \\
            &=\frac{1}{p(p-1)}\left(\sum_{i=1}^{p} \sum_{j=1, i\ne j}^{p-1} \mathbf{A}(\dots, \mathbf{B}, \dots \mathbf{C})\right).
        \end{split}
    \end{equation}
    The \(\circledcirc\) is used for the terms \(\mathbf{B}\) and \(\mathbf{C}\) to maintain them as indeterminates. This prevents matrix multiplications and vector symmetric products while allowing commutative inputs on symmetric product. All the same properties of multiplication apply to \(\circledcirc\). Same results hold for when different size dimensions are needed.
\end{definition}
With the symmetric tree product given for tensors then the stumps are defined via the differential operator. 
\begin{definition}[Differential Operator for Stumps]\label{def:stumpDiffOp}
    The elementary differentials for stumps are given as
    \begin{equation}
        \fun\comp{c}_{(\y\comp{c_{1}})^{b_1} \dots (\y\comp{c_{p}})^{b_p}}\BDiffOp(\xnode_1^{b_1} \xnode_2^{b_2} \dots \xnode_p^{b_p})(\y) = \fun\comp{c}_{(\y\comp{c_{1}})^{b_1} \dots (\y\comp{c_{p}})^{b_p}}.
    \end{equation}
    These nodes are considered as only taking the derivative and attaching no new terms. If the stump has no parent, then
    \begin{equation}
        \BDiffOp(\xnode)(\y_0) = \mathbf{I}.
    \end{equation}
    Stump nodes cannot have children. Given stumps are treated as tensors then for a tree with stumps \(\dforest = \xnode_1^{b_1} \xnode_2^{b_2} \dots \xnode_p^{b_p}\) attached \(\tree \in \NTTC\lceil \dforest \rceil\) and for trees \(\forestTQs \in \ForestSpace{\NTTC}\) then the symmetric tree product gives
    \begin{equation}
        \BDiffOp(\tree)(\y_0) \symprod \left(\prod_{i=1}^{m} \BDiffOp(\stree_{i})(\y_0)^{q_i}\right).
    \end{equation}
    This applies each tree \(\stree_i\) to a corresponding appropriate stump to graft on to. All properties of \cref{def:symProduct} hold over stumps.
\end{definition}
\begin{example}[Examples of Stumps]
    Some examples of stumps in action are as follows:
    \begin{align}
        &\BDiffOp(\tau[\xnode^2])(\y_0) \BDiffOp(\tau)(\y_0) \BDiffOp(\tau)(\y_0) = \BDiffOp(\tau[\tau^2])(\y_0),\\
        &\BDiffOp(\tau[ \tau[\xnode] \xnode])(\y_0)\BDiffOp(\tau)(\y_0) \BDiffOp(\tau[\tau])(\y_0) = \frac{1}{2} \left(\BDiffOp(\tau[ \tau[\tau[\tau]] \tau])(\y_0) + \BDiffOp(\tau[ \tau[\tau]^2])(\y_0)\right),\\
        &\BDiffOp(\tau_1[\tau_1[\xnode_1] \xnode_2])(\y_0)\BDiffOp(\tau_1)(\y_0) \BDiffOp(\tau_2)(\y_0) = \BDiffOp(\tau_1[\tau_1[\tau_1] \tau_2])(\y_0),\\
        &\BDiffOp(\tau[\tau \xnode])(\y_0)^2 = \BDiffOp(\tau[\tau \tau[\tau \xnode]])(\y_0),\\
        &\BDiffOp(\tau[\xnode^2])(\y_0) (a_1\BDiffOp(\tau)(\y_0) + a_2\BDiffOp(\tau[\tau])(\y_0))^2 = \\
        & \qquad a_1^2 \BDiffOp(\tau[\tau^2])(\y_0) + 2 a_1\, a_2\, \BDiffOp(\tau[\tau \tau[\tau]])(\y_0) + a_2^2 \BDiffOp(\tau[\tau[\tau]^2])(\y_0)\\
        & \BDiffOp(\tau[\tau[\xnode] \xnode])(\y_0) \symprod (\BDiffOp(\tau[\xnode^2])(\y_0) \circledcirc \BDiffOp(\tau[\xnode])(\y_0)) \\
        &= \frac{1}{2}\left(\BDiffOp(\tau[\tau[\tau[\xnode^2]] \tau[\xnode]])(\y_0) + \BDiffOp(\tau[\tau[\tau[\xnode]] \tau[\xnode^2]])(\y_0)\right)
    \end{align}
\end{example}
We note similarities to the stumps in \cite{butcherTreesStumpsApplications2018,Butcher_2021_Bseries-book} and grafting operations in \cite{merson1957operational} however, a key difference is we do not consider the product of trees but only the products over elementary differentials. This allows for more expressive labels to be provided and preserve commutativity between trees by extending the set of trees.

\subsection{\BSeriesText}
\begin{definition}[Label/Coefficient]\label{def:coefficient}
    A label/coefficient is defined as some function of trees \cite{Butcher_2021_Bseries-book}.
    The label function is to be considered element-wise operation when applied to trees. The label has the following properties
    \begin{equation}
        \begin{split}
            a(c_1 \fforest_1 + c_2 \fforest_2) &= c_1\, a(\fforest_1) + c_2\, a(\fforest_2),\\
            a(\vec{\tree}) &= \begin{bmatrix}
                a(\tree_1)\\
                a(\tree_2)\\
                \vdots\\
                a(\tree_n)
            \end{bmatrix},\\
            a(\treecap) &= \begin{bmatrix}
                a(\tree_{1,1}) & a(\tree_{1,2}) & \dots & a(\tree_{1,m})\\
                a(\tree_{2,1}) & a(\tree_{2,2}) & \dots & a(\tree_{2,m})\\
                \vdots & \vdots & \ddots & \vdots\\
                a(\tree_{n,1}) & a(\tree_{n,2}) & \dots & a(\tree_{n,m})
            \end{bmatrix},
        \end{split}
    \end{equation}
    where \(c_1,c_2 \in \outspace\), \(\tree_i, \tree_{i,j} \in \AbstractT\), \(\fforest_1, \fforest_2 \in \ForestSpace{\AbstractT}\), \(\vec{\tree} \in \AbstractT^n\), \(\treecap \in \AbstractT^{n \times m}\).
\end{definition}

A \BSeriesText, as defined in \cite{Butcher_2021_Bseries-book, Hairer_2006_geometric-book, Hairer_book_I, Hairer_book_II}, is given as follows.
\begin{definition}[\BSeriesText]\label{def:ClassicBSeries}
    For a label \(a\) the classical \BSeriesText is defined as
    \begin{equation}
        \BSeries{\Dt}(a, \y) = \sum_{\tree \in \TT}  a(\tree)\,\frac{\Dt^{\treeorder{\tree}}}{\sigma(\tree)}\BDiffOp(\tree)(\y).
    \end{equation}
    The series matches the exact Taylor series of the ODE \(\frac{\ddiff \y}{\ddiff t} = \fun(\y)\) when \(a(\tree) = \frac{1}{\gamma(\tree)}\).
\end{definition}
Note that in subsequent sources they keep \(\emptyset\) as separate treatment. For this definition it is treated all as one since \(a(\emptyset)\frac{\Dt^{0}}{\sigma(\emptyset)}\BDiffOp(\emptyset)(\y) = a(\emptyset)\y\).\par

For partitioned integration the use of multiple colors are introduced with the following series.
\begin{definition}[\NBSeriesText]\label{def:NBSeries}
    Given a label \(a\) and set of colors \(\Tcolors\) an \NBSeriesText is defined as follows:
    \begin{equation}
         \NBSeries{\Dt}(a, \y) = \sum_{\tree \in \NTTC} a(\tree)\frac{\Dt^{\treeorder{\tree}}}{\sigma(\tree)}\BDiffOp(\tree)(\y).
    \end{equation}
    The series matches the exact Taylor series of the ODE when \(a(\tree) = \frac{1}{\gamma(\tree)}\).
\end{definition}
For multicolored trees there are two common series in literature to choose from. \PSeriesText \cite{Hairer_2006_geometric-book, hairerOrderConditionsNumerical1981, murua1995metodos, muruaFormalSeriesNumerical1999} and \NBSeriesText \cite{araujo1997symplectic, muruaFormalSeriesNumerical1999}. PSeries is typically used for component partitioned methods whereas additive partitions utilize \NBSeriesText. However, \PSeriesText is built with only \(2\) component partitions in mind whereas \NBSeriesText is for arbitrary partitions. Thus, \NBSeriesText will be used for the rest of the paper. Note that for the remainder of the paper for simplicity \NBSeriesText will be referred to as \BSeriesText textual since the only difference is the set of trees they sum over. However, for mathematical consistency \(\NBSeries{\Dt}(a, \y)\) is utilized in equations for indicating there is more than one color.

Next an extension of \NBSeriesText for stumped series.
\begin{definition}[\NBSeriesText for Stumps]\label{def:NBSeriesStumped}
    Given a label \(a\), set of colors \(\Tcolors\), and stumps \(\dforest \in \ForestSpace{\TStumpcolors}\) then the stumped variation of \NBSeriesText is as follows:
    \begin{equation}
         \NBSeries{\Dt}^{\dforest}(a, \y) = \sum_{\tree \in \NTTC\lceil \dforest \rceil} a(\tree)\frac{\Dt^{\treeorder{\tree}}}{\sigma(\tree)}\BDiffOp(\tree)(\y).
    \end{equation}
    If \(\dforest = \emptyset\) then the standard \NBSeriesText is given.
\end{definition}
Such series will be used to describe the underlying approximations of sub-systems \cref{eq:ApproximationProblem}. All operations defined for \BSeriesText in this section hold for stumped variations.
Next some properties of \BSeriesText.
\begin{property}[\NBSeriesText Time Scaling]\label{property:timescale}
    For a label \(a\) and a given time step \(\alpha \Dt\) the time scaling is
    \begin{equation}
        \NBSeries{\alpha \Dt}(a, \y) = \NBSeries{\Dt}(\tilde{a}, \y)
    \end{equation}
    where \(\tilde{a}(\tree) = \alpha^{\treeorder{\tree}}a(\tree)\).
\end{property} 
A separate derivation of \BSeriesText will be utilized throughout this text. This takes the \BSeriesText definition and utilizes an inner product to expand the sum instead \cite{Butcher_2021_Bseries-book}. Let \(\treevecinf\) be an infinite size vector where each entry is a tree. A \BSeriesText can be described in inner product form as follows:
\begin{equation}
    \begin{split}
        \BSeries{\Dt}(a, \y) &=\left( \BDiffOp(\treevecinf, \y) \frac{\Dt^{\treeorder{\treevecinf}}}{\sigma(\treevecinf)} \right)^{T} a(\treevecinf)\\
        &=\scalemath{1.0}{\begin{bmatrix}
            \BDiffOp(\tree_1, \y) \frac{\Dt^{\treeorder{\tree_1}}}{\sigma(\tree_1)} & \BDiffOp(\tree_2, \y) \frac{\Dt^{\treeorder{\tree_2}}}{\sigma(\tree_2)} & \BDiffOp(\tree_3, \y) \frac{\Dt^{\treeorder{\tree_3}}}{\sigma(\tree_3)} & \dots 
        \end{bmatrix} \begin{bmatrix}
            a(\tree_1)\\
            a(\tree_2)\\
            a(\tree_3)\\
            \vdots
        \end{bmatrix}}.
    \end{split}
\end{equation}
The inner product form provides a linear operator insight into some more general operations. A similar formulation is used for \NBSeriesText.
For simplicity throughout the text a color compression is introduced. The compression conforms all colors down to one color by representing a sub-vector of the full space.
\begin{definition}[Color Compression]\label{def:compression}
    For a single color \(c\) let it be represented as a sum or partition of multiple colors \(d \in \Tcolors \setminus \{c\}\). Let the vector \(\treevecinf\) be the vector of all trees of colors \(\Tcolors \setminus \{c\}\). For a tree \(\tree\) with at least one node of color \(c\) it represents a sub-vector of \(\treevecinf\) for all possible unique permutation of the colors notational indexed by \(\blockIndex{\treevecinf}{\tree}\). Matrices of trees can be indexed by them to obtain block matrices as \(\blockIndex{\mathrm{A}}{\tree, \stree}\). All label operations still operate element-wise.
\end{definition}
This compression of colors will be utilized throughout to induce block-wise operations. 
\begin{example}
    As an example let the color \(\btree{[x]} = \{\btreega{[c1]}, \btreega{[c2]}\}\) then the trees represent 
    \begin{equation}
        \begin{split}
            \btree{[x[x]]} &= \begin{bmatrix}
                \btreega{[c1[c1]]} & \btreega{[c1[c2]]} & \btreega{[c2[c1]]} & \btreega{[c2[c2]]}
            \end{bmatrix}^T,\\
            \btree{[x[x][x]]} &= \begin{bmatrix}
                \btreega{[c1[c1][c1]]} & \btreega{[c1[c2][c1]]} & \btreega{[c1[c2][c2]]} & \btreega{[c2[c1][c1]]} & \btreega{[c2[c2][c1]]} & \btreega{[c2[c2][c2]]}
            \end{bmatrix}^T.\\
        \end{split}
    \end{equation}
    Evaluating the compressed colors over \(\BDiffOp(\tree)(\y)\, a(\tree)\, \frac{\Dt^{\treeorder{\tree}}}{\sigma(\tree)}\) gives that
    \begin{equation}
        \scalemath{1.0}{
        \begin{split}
            \BDiffOp\left(\btreega{[x[x]]}\right)(\y)\, a\left(\btreega{[x[x]]}\right) \Dt^{2} &= \Dt^2 a\left(\btreega{[c1[c1]]}\right)\BDiffOp\left(\btreega{[c1[c1]]}\right)(\y) + \Dt^2 a\left(\btreega{[c1[c2]]}\right)\BDiffOp\left(\btreega{[c1[c2]]}\right)(\y) \\
            &+ \Dt^2 a\left(\btreega{[c2[c1]]}\right)\BDiffOp\left(\btreega{[c2[c1]]}\right)(\y) + \Dt^2 a\left(\btreega{[c2[c2]]}\right)\BDiffOp\left(\btreega{[c2[c2]]}\right)(\y),\\
            \BDiffOp\left(\btreega{[x[x][x]]}\right)(\y)\, a\left(\btreega{[x[x][x]]}\right) \frac{\Dt^{3}}{2} &= \dots + \frac{\Dt^3}{2} a\left(\btreega{[c1[c1][c1]]}\right)\BDiffOp\left(\btreega{[c1[c1][c1]]}\right)(\y)\\
            &+ \Dt^3 a\left(\btreega{[c1[c1][c2]]}\right)\BDiffOp\left(\btreega{[c1[c1][c2]]}\right)(\y) + \frac{\Dt^3}{2} a\left(\btreega{[c1[c2][c2]]}\right)\BDiffOp\left(\btreega{[c1[c2][c2]]}\right)(\y) + \dots .
        \end{split}
        }
    \end{equation}
\end{example}
Lastly a common technique to achieve order conditions for arbitrary stage number \(s\) is to ``stack'' \BSeriesText. This can be thought as letting a coefficient be a vector where all operations on it are element-wise.
\begin{definition}[\BSeriesText Stacking]\label{def:stacking}
    For the labels \(a_{i}\) where \(i \in \{1 .. s\}\) the stacked \BSeriesText is defined as
    \begin{equation}
        \NBSeries{\Dt}(a,\y) = \begin{bmatrix}
            \NBSeries{\Dt}(a_{1},\y)\\
            \NBSeries{\Dt}(a_{2},\y)\\
            \vdots\\
            \NBSeries{\Dt}(a_{s},\y)\\
        \end{bmatrix}.
    \end{equation}
    Note that multiplication of labels \(a\) with other labels of the same size acts as elementwise multiplication.
\end{definition}
Stacking allows for a simplification of order conditions down to matrix vector products. This is used commonly to simplify order conditions \cite{Sandu_2015_GARK}.

\section{Combinatorics Tools\label{section:Combinatorics}}
Next some new combinatorics tools will be introduced for the assistance in defining tree operations. To the best of the authors' knowledge the tools are not currently used in \BSeriesText analysis literature. We provide definitions and properties of the \textbf{matrix permanent} \cite{minc1984permanents} for trees.

\subsection{Generalized Outer Product}
To start a generalized outer product is defined. Typically, the outer product is defined as multiplying the entries to form a matrix between two vectors. A generalization is formed when multiplication is replaced by an arbitrary binary operation.
\begin{definition}[Generalized Outer Product]\label{def:outerproduct}
    For $f(a,b) : \mathcal{A} \cross \mathcal{B} \rightarrow \mathcal{C}$ and the vectors $a = [a_1 \, a_2 \, a_3 \, \dots \, a_n] \in \mathcal{A}^{n}$ and $b = [b_1 \, b_2 \, b_3 \, \dots \, b_m] \in \mathcal{B}^{m}$ the generalized outer product is defined as
    \begin{equation}
        a \outergen{f} b = \begin{bmatrix}
            f(a_1,b_1) & f(a_1,b_2) & \cdots & f(a_1,b_m)\\
            f(a_2,b_1) & f(a_2,b_2) & \cdots & f(a_2,b_m)\\
            \vdots & \vdots & \ddots & \vdots\\
            f(a_n,b_1) & f(a_n,b_2) & \cdots & f(a_n,b_m)
        \end{bmatrix}.
    \end{equation}
\end{definition}
Though the vector spaces are kept arbitrary here the preferred space is forest. Each forest has a similar representation where each element of the vector is a tree instead. The unordered nature allows the vector to be the same under arbitrary permutations. Such effects are mentioned further on.
\begin{remark}[Forest Outer Product Padding]
    We note that for the cases we consider with forests are only for square matrices. Thus, the tree \(\emptyset\) is considered as a ``padding'' such that the resulting matrix is square. As an example, for the forest \(\fforest = \forestTQt\) and \(\gforest = \forestTQs\) where \(\sum_{i=1}^{n} q_{i} = \bar{n} > \sum_{i=1}^{m} k_{i} =  \bar{m}\) the generalized outer product gives
    \begin{equation}
        (\fforest \outergen{f} \gforest)_{i,j} = f(\tree_{\pi(i)}^{q_{\pi(i)}}, \emptyset) \quad \forall\, j > \bar{n} - \bar{m}
    \end{equation}
    where \(\pi\) is some integer mapping. Such padding may be denoted by \(\gforest\, \emptyset^{\bar{n} - \bar{m}}\). The case is symmetric for when \(\bar{m} > \bar{n}\). The empty sets used for padding do not consider coloring. All binary operations for this paper will require that any amount of padding, that keeps the matrix square, will give the same result.
\end{remark} 

\subsection{Symmetric Sets}
\begin{definition}[Symmetric Set]\label{def:symmetricgroup}
    The group \(S_n\) is the set of all permutations of \(\{1 .. n\}\). The set has the size \(|S_n| = n!\).
\end{definition}
\begin{example}
    Given the set \(S_3\) contains the sets
    \[\{(1,2,3), (2,1,3), (2,3,1), (3,2,1), (3,1,2), (1,3,2)\} = S_3.\]
    For the element \((2,1,3) = \sigma \in S_3\) the indexings are \(\sigma(1) = 2, \sigma(2) = 1, \sigma(3) = 3\).
\end{example}
This generalizes to the case of sets with duplicates called multi-sets.
\begin{definition}[Multi-Set Permutations Set] \label{def:multisetpermutations}
    Let \(\omega_i\) be some unique character indexed by \(i \in \{1 .. n\}\) and let the multi-set \(f = \omega_1^{k_1} \omega_2^{k_2} \dots \omega_n^{k_n}\) be the set irrespective of permutations where \(\bar{n} = \sum_{i=1}^{n} k_i\). The set \(S_f\) is defined as the set of all unique permutations of the multiset \(f\). Thus, \(S_f\) has a total size of \(|S_f| = \frac{\bar{n}!}{\prod_{i=1}^{n} k_i!}\).

    The same notion of \(\sigma \in S_f\) works as before with symmetric sets. Let \(\sigma(j)\) be the mapping of some index of a multi-set that returns the index to the corresponding characters with input \( 1 \le j \le \bar{n}\) and output one of the \(n\) characters. Furthermore, let \(\mu_f(p) = k_p\) be a function which returns the number of occurrences in the sequence.
\end{definition}
\begin{example}
    To demonstrate multi-sets for \(f = \omega_1^2 \omega_2^2\) the multi-set is 
    \begin{equation}
        \begin{split}
            &\{(\omega_1, \omega_1, \omega_2, \omega_2), (\omega_1, \omega_2, \omega_1, \omega_2), (\omega_1, \omega_2, \omega_2, \omega_1), (\omega_2, \omega_2, \omega_1, \omega_1),\\
            &(\omega_2, \omega_1, \omega_2, \omega_1), (\omega_2, \omega_1, \omega_1, \omega_2)\}
        \end{split} = S_f
    \end{equation}
    where each element is a unique permutation of the characters \(\omega_1,\omega_2\).
    With the indexing notation if \((\omega_1, \omega_2, \omega_1, \omega_2) = \sigma \in S_f\) then \(\sigma(1) = 1, \sigma(2) = 2, \sigma(3) = 1, \sigma(4) = 2\).
    Furthermore, the symmetry factors of the multiset are \(\mu_f(1) = 2, \mu_f(2) = 2\).
\end{example}
Even though this was defined over characters this will be utilized with trees instead.

\subsection{Matrix Permanent}
\begin{definition}[Matrix Permanent\cite{minc1984permanents}]\label{def:permanent}
    Let $\permanent: \mathcal{A}^{n \times n} \rightarrow \mathcal{A}$ where \(\mathcal{A}\) is a commutative ring w.r.t. multiplication.
    Over a symmetric set let \(\sigma \in S_n\) be the symmetric group and \(\sigma\) is the permutation mapping. Given a matrix $A \in \mathcal{A}^{n \times n}$ the permanent is defined as 
    \begin{equation}
        \permanent(A) = \sum_{\sigma \in S_n} \prod_{j=1}^{n} a_{\sigma(j),j}.
    \end{equation}
    The Laplace (co-factor) expansion for the \(\permanent\) is defined as
    \begin{equation}
        \permanent(A) = \sum_{j=1}^n a_{j,i}\,\permanent(A_{[j,i]^{-1}}) \quad \text{ for some column } i,
    \end{equation}
    where $A_{[j,i]^{-1}}$ is the matrix \(A\) with row \(j\) and column \(i\) omitted.
\end{definition}
\begin{remark}[Permanent + Outer Product]\label{remark:permanentOuter}
    Given some function $f : \mathcal{A} \cross \mathcal{A} \rightarrow \mathcal{A}$, $a = [a_1 \, a_2 \, a_3 \, \dots \, a_n] \in \mathcal{A}^{n}$, and $b = [b_1 \, b_2 \, b_3 \, \dots \, b_n] \in \mathcal{A}^{n}$ the combination of the permanent (\cref{def:permanent}) and the generalized outer product (\cref{def:outerproduct}) gives
    \begin{equation}
        \permanent(a \outergen{f} b).
    \end{equation}
    This will give a polynomial result of all possible combination of sets $a$ and $b$ applied to $f$ combined.
\end{remark}
\begin{example}\label{example:permanent}
    As an example for the list $v = [a,b,c]$ and $u = [d,e,f]$ with \(\phi(x,y) = x\, y\) the permanent + outer product gives 
    \begin{equation}
        \begin{split}
            \permanent(v \outergen{\phi} u) &= \permanent\left(\begin{bmatrix}
                ad & ae & af\\
                bd & be & bf\\
                cd & ce & cf\\
            \end{bmatrix}\right) \\
            = &(ad)[(be)(cf) + (bf)(ce)] + (bd)[(ae)(cf) + (af)(ce)] +\\
            & \quad  (c d)[(ae)(bf) + (af)(be)].
        \end{split}
    \end{equation}
\end{example}
Properties of the permanent and with combinations of the outer product are provided.
\begin{property}[Permanent Permutations {\cite[Chapter 43.1]{hogben2006handbook}\cite{minc1984permanents}}]\label{property:permanentpermute}
    The permanent of any matrix \(A\) multiplied by any permutation matrices \(P, Q\) will give equivalent results
    \begin{equation}
        \permanent(A) = \permanent(P A Q) \hspace*{1cm} \forall\, P,\, Q.
    \end{equation}
\end{property}
\begin{remark}[Forest and Permanent]
    Note that since the permanent is invariant to any permutations of the states, then this preserves order invariance in the forest. To demonstrate this for forests \(\fforest = \forestTQt\) and \(\gforest = \forestTQs\) with vector permutations \(a,b\) for each respectively then due to \cref{property:permanentpermute} this gives 
    \[\permanent(a \outergen{\star} b) = \permanent(P (a \outergen{\star} b) Q) = \permanent(\fforest \outergen{\star} \gforest ).\]
\end{remark}
\begin{property}[Permanent Scalar Multiplication  {\cite[Chapter 43.1]{hogben2006handbook}}]\label{property:permanentconstant}
    For a matrix \(A \in \mathcal{A}^{n,n}\) and scalar \(c \in \mathcal{A}\) the permanent gives
    \begin{equation}
        \permanent(c A) = c^{n} \permanent(A).
    \end{equation}
    This comes from principal-minors expansion of the permanent.
    Furthermore, if the matrix has repeated columns as \(B_{i,j} = b, \quad \forall\, i,\,j \in \{1,..n\}\) the permanent gives
    \begin{equation}
        \permanent(B) = n!\, b^{n}.
    \end{equation}
\end{property}
Another property comes from the case of whether a permanent is zero or not. In \cite[Section 3.2]{minc1984permanents} this is proven for non-negative matrix for both directions. It is noted that the proof performed in \cite[Section 3.2;Theorem 2.2]{minc1984permanents} is for an ``if and only if'' condition. When looking at only the converse the proof does not specify the need for non-negativity nor does it require it. For results in this paper only the converse case is utilized.
\begin{property}[Permanent Zero {\cite[Section 3.2;Theorem 2.2]{minc1984permanents}}]\label{property:permanentzero}
     Let \(A\) be some square matrix of size \(n\) and suppose there exists a submatrix of dimensions \(s \times t\) that is a zero matrix. If \(s + t \ge n + 1\) then \(\permanent(A) = 0\).
\end{property}
\begin{lemma}[Permanent Diagonal Scaling]\label{lemma:permanentscale}
    For a matrix \(A_{i,j} = a_{i,j}\) \(\forall\, i,\,j \in \{1,..n\}\) and a diagonal matrix \(\mathrm{Diag}(\vec{c})\), where the diagonal is the vector, the permanent simplifies to
    \begin{equation}
        \permanent(\mathrm{Diag}(\vec{c}) \, A) = \permanent(A \, \mathrm{Diag}(\vec{c}) ) = \permanent(A) \, \permanent(\mathrm{Diag}(\vec{c})) = \permanent(A) \, \left(\prod_{i=1}^{n} c_i \right).
    \end{equation}
\end{lemma}
\begin{proof}
    By the definition via symmetric groups
    \[\permanent(\mathrm{Diag}(\vec{c})\, A) =  \permanent(A\, \mathrm{Diag}(\vec{c}) ) = \sum_{\sigma \in S_n} \prod_{j=1}^{n} a_{\sigma(j),j} c_j.\]
    Since the product of \(c\) is formed over each point on \(j\) it can be extracted out as product giving 
    \[= \sum_{\sigma \in S_n} \prod_{j=1}^{n} a_{\sigma(j),j}\, \left( \prod_{k=1}^n c_k\Big) = \permanent(A)\Big(\prod_{i=1}^{n} c_i \right).\]
\end{proof}

Next for properties combining the permanent and outer product.
\begin{lemma}[Permanent, Outer Product, Repeated Columns]\label{lemma:permanentmultiset}
    Let \(\star\) be some arbitrary binary operation. Given the matrix \(\mathbf{K} = \vec{a} \outergen{\star} \vec{d}\) where \(\vec{a} = [a_1 \, a_2 \,  \dots \, a_{\bar{m}}]\) and \\\(\vec{d} = [b_1 \vec{1}_{k_1},\, \,b_2 \vec{1}_{k_2},\, \, \dots,\, \,b_m \vec{1}_{k_m}]\) where \(\bar{m} = \sum_{i=1}^{m}k_i\). Given the vector \(\vec{d}\) with the forest representation \(g = b_1^{k_1} \, b_2^{b_2} \, \dots \, b_m^{k_m}\) with a mapping \(\pi\) based on \cref{def:multisetpermutations} the combination of permanent, outer product, and repeated columns gives
    \begin{equation}
        \permanent(\mathbf{K}) = \sum_{\sigma \in S_{\bar{m}}} \prod^{\bar{m}}_{j=1}  (a_{j} \star g_{\pi(\sigma(j))}) = \left(\prod_{i=1}^{m}k_i!\right) \sum_{\sigma \in S_{g}} \prod_{j=1}^{\bar{m}} (a_{j} \star g_{\sigma(j)}).
    \end{equation}
    Thus, changing the operations from a symmetric set \(S_{n}\) to a multi-symmetric set of the forest \(S_{g}\). Additionally, for the symmetric case, letting \(\mathbf{K}' =  \vec{d} \outergen{\star} \vec{a}\) gives the result
    \begin{equation}
        \permanent(\mathbf{K}') = \sum_{\sigma \in S_{\bar{m}}} \prod^{\bar{m}}_{j=1}  (g_{\pi(\sigma(j))} \star a_{j} ) = \left(\prod_{i=1}^{m}k_i!\right) \sum_{\sigma \in S_{g}} \prod_{j=1}^{\bar{m}} (g_{\sigma(j)} \star  a_{j}).
    \end{equation}
\end{lemma}
\begin{proof}
    For simplicity let \(\sigma \in S_{\bar{m}}\) and \(\tilde{\sigma} \in S_{g}\). Given some sequence of characters with repeats \(\vec{d}\) with the forest representation \(g\) the map \(\pi\) will take a permutation set operation and map to forest. However, under the permutation if similar characters are swapped under various permutations they will produce the same product. This is formalized by the following: suppose that \(\sigma_1, \sigma_2 \in S_{\bar{m}}\) where \(\sigma_2\) is a single permutation difference of the same character i.e. swapping \(b_i\) with another \(b_i\) then 
    \[\pi(\sigma_1(j)) = \pi(\sigma_2(j)) \quad \forall\, j \in \{1,\bar{m}\}\]
    since there will be \(k_i\) such repeated permutations that produce the same result. Let \(\tilde{\sigma} \in  S_{g}\) be the forest representation and \(\sigma_{i,q} \in S_{\bar{m}}\) be the set of results where only the character \(b_i\) is permuted and matches the sequence of \(\tilde{\sigma}\). This changes the indexing to 
    \[\tilde{\sigma}(j) = \pi(\sigma_{i,q}(j)) \quad \forall\, j \in \{1,\bar{m}\} \quad q \in \{1,k_i\}\]
    where for the characters \(b_i\) there are \(k_i!\) duplications summed. Multiplying over all duplications gives the results of all permutations of duplications. This leads to a leading factor of \(\left(\prod_{i=1}^{m}k_i!\right)\). Converting the rest of the expression down to multi-sets gives the final result. The same argument holds for the case of \(K'\) given the co-factor expansion of the permanent.
\end{proof}
\begin{lemma}[Permanent and Trees]\label{lemma:permanentTrees}
    Let \(\mathrm{T}\) be some set of trees. Given the operator \(\star: \ForestSpace{\mathrm{T}} \times \ForestSpace{\mathrm{T}} \rightarrow \ForestSpace{\mathrm{T}}\) and the forest \(\fforest = \forestTQt \in \ForestSpace{\mathrm{T}}_{\bar{n}}\) where  \(\bar{n} = \sum_{i=1}^{n} q_i\)then the following holds  
    \begin{equation}
            \sum_{\gforest = \forestTQs \in \ForestSpace{\mathrm{T}}_{\bar{n}}} \frac{1}{\prod_{j=1}^{m} k_j!} \permanent\left(\fforest \, \outergen{\star} \, \gforest \right) = \prod_{i=1}^{n} \left( \sum_{\stree \in \mathrm{T}} (\tree_{i} \star \stree)\right)^{q_i}.
    \end{equation}
\end{lemma}
\begin{proof}
    Starting from RHS of the equality. Let \(\pi\) be an index mapping to the trees in \(\fforest\). Since the order does not matter then the choice of \(\pi\) is arbitrary giving
    \begin{equation}
        \prod_{i=1}^{n} \left( \sum_{\stree \in \mathrm{T}} (\tree_{i} \star \stree)\right)^{q_i} = \prod_{i=1}^{\bar{n}} \left( \sum_{\stree \in \mathrm{T}} (\fforest_{\pi(i)} \star \stree)\right).
    \end{equation}
    Next since there are no exponents left, then the product is expanded to a sum over all tree possibilities
    \begin{equation}
        \sum_{\stree_1, \stree_2, \dots \stree_{\bar{n}} \in \mathrm{T}} \left[\prod_{i=1}^{\bar{n}} (\fforest_{\pi(i)} \star \stree_{i}) \right].
    \end{equation}
    If the restricted forest set \(\gforest \in \ForestSpace{\mathrm{T}}_{\bar{n}}\) is taken then all unordered permutations are given. This is equivalent of summing over the symmetric multi-set \(S_{\gforest}\) sum giving
    \begin{equation}
        \sum_{\gforest \in \ForestSpace{\mathrm{T}}_{\bar{n}}} \sum_{\sigma \in S_{\gforest}} \prod_{j=1}^{\bar{n}} (\fforest_{\pi(j)} \star \gforest_{\sigma(j)}).
    \end{equation}
    Utilizing \cref{lemma:permanentmultiset} and inverting the product of symmetries \(\prod_{j=1}^{m} k_j!\) gives the final LHS result.
\end{proof}
It is noted that \cite[Section 5.4]{chabaudQuantuminspiredPermanentIdentities2022} works through formal variables of generating functions on vector products to unitary matrices. Similar results for repeated rows and columns are obtained in this case with trees and outer product.
\begin{remark}
    The purpose of the permanent here is the property of inducing all permutations of terms in the matrix. Thus, if the entries of the matrix are composed of stumped elementary differentials then we assume that \(\circledcirc\) is used over multiplication. This is to allow for all permutations in intermediates form.
\end{remark}
The main point in using the permanent is to describe a simple way to traverse binary operations of a tree recursively over the children. However, in practice there are computational difficulties with the permanent computation \cite{valiantComplexityComputingPermanent1979}. Thus, a more practical method utilizing memory via dynamic programming can improve performance of practical implementations. This leads to a linear form which works recursively by using results for lower cardinality forest. In practice the lower order forest terms are needed for other terms thus, can be stored in memory to reduce computation. The recursive form to reduce, and the linear form to avoid direct permanent computation are given as follows.
\begin{theorem}[Linear Recursive Permanent Form]\label{theorem:generalLinearTreeRecursion}
    Given some binary operation \(\star\) over forest \(\fforest = \forestTQt\) and \(\gforest = \forestTQs\) where the binary operation must return the same result for arbitrary \(\emptyset\) paddings. If \(\star\) is defined for forest as
    \begin{equation}
        \begin{split}
            \fforest \star \gforest &= \frac{1}{(\prod_{i=1}^m k_i !)(p!)}\permanent(\forestTQt \outergen{\star} \forestTQs \emptyset^{p})\\
            \emptyset \star \gforest &= 0
        \end{split}
    \end{equation}
    where \(p = \max(\bar{n} - \bar{m},0)\) then the linear form is given as
    \begin{equation}
        \begin{split}
            \fforest \star \gforest &=  \sum_{j=1}^{m} (\tree_i \star \stree_j)\left(\frac{\fforest}{\tree_i} \star \frac{\gforest}{\stree_j} \right) + (\tree_i \star \emptyset) \, \left(\frac{\fforest}{\tree_i} \star \gforest \right)\, \textit{ for some } i \in \{1 .. n\}.\\
        \end{split}
    \end{equation}
    An alternative formulation based of \(\gforest\) selection gives 
    \begin{equation}
        \begin{split}
            \fforest \star \gforest &= \sum_{i=1}^{n} \frac{q_i}{k_j} (\tree_i \star \stree_j)\left(\frac{\fforest}{\tree_i} \star \frac{\gforest}{\stree_j} \right) + \frac{\tilde{q}}{k_j} (\emptyset \star \stree_j)\left(\fforest \star \frac{\gforest}{\stree_j} \right)\, \textit{if } \bar{n} \ge \bar{m} \quad \textit{ for some } j \in \{1 .. m\}\\
             &= \sum_{i=1}^{n} \frac{q_i}{p} (\tree_i \star \emptyset)\left(\frac{\fforest}{\tree_i} \star \gforest \right) \textit{if } \bar{n} > \bar{m}.\\
        \end{split}
    \end{equation}
\end{theorem}
\begin{proof}
    Given the definition of the permanent \cref{def:permanent} there is the Laplace (co-factor) expansion form. For abuse of notation, division by the tree is removing it from the forest.
    By removing \(\tree_i\) this is removing the \(i\)th row of the matrix permanent for the co-factor expansion. With some arbitrary integer map \(\pi\) the co-factor expansion this gives
    \begin{equation}
        \begin{split}
            \fforest \star \gforest &= \frac{1}{\prod_{i=1}^{m}k_j! p!} (\sum_{j=1}^{\bar{m}} (\tree_i \star \stree_{\pi(j)}) \permanent\left(\frac{\fforest}{\tree_i} \outergen{\star} \frac{\gforest}{\stree_\pi(j)}\right)  \\
            +& \sum_{j=\bar{m} + 1}^{\bar{n}} \tree_i \permanent\left(\frac{\fforest}{\tree_i} \outergen{\star} \gforest\, \emptyset^{p-1} \right) )\textit{ for some } i \in \{1 .. n\}.
        \end{split}
    \end{equation}
    By applying \(\star\) the result is
    \begin{equation}
        \begin{split}
            \fforest \star \gforest &= (\sum_{j=1}^{\bar{m}} \frac{1}{k_{\pi(j)}} (\tree_i \star \stree_{\pi(j)}) \left(\frac{\fforest}{\tree_i} \star \frac{\gforest}{\stree_\pi(j)} \right)  \\
            +& \sum_{j=\bar{m}+1}^{\bar{n}} \frac{1}{p} \tree_i \left(\frac{\fforest}{\tree_i} \star \gforest \right) )\textit{ for some } i \in \{1 .. n\}.
        \end{split}
    \end{equation}
    Given this is simply the sum over the multiplicities then by \cref{lemma:permanentmultiset} this reduces to the statement's result.
    For the alternative definition taking the Laplace expansion, taking the \(j\) column tree, and applying the \(\star\) definition gives
    \begin{equation}
        \begin{split}
            \fforest \star \gforest &= \sum_{i=1}^{\bar{n}} \frac{q_{\pi(i)}}{k_{j}} (\tree_{\pi(i)} \star \stree_j) \left(\frac{\fforest}{\tree_{\pi(i)}} \star \frac{\gforest}{\stree_j} \right)  \\
            & \textit{if } \bar{n} \ge \bar{m} \textit{ for some } j \in \{1 .. m\},\\
            &= \sum_{i=1}^{\bar{n}} \frac{q_{\pi(i)}}{p} (\tree_{\pi(i)} \star \emptyset) \left(\frac{\fforest}{\tree_{\pi(i)}} \star \gforest \right)   \textit{if } \bar{n} > \bar{m}.\\
        \end{split}
    \end{equation}
    Again this is the sum over multiplicities thus, this reduces to the statements result.
\end{proof}
Note that the permanent form can provide some benefits in presentation of binary operations due to the recursive nature. In practice when implementing such operations in software it is preferred to use the form from \cref{theorem:generalLinearTreeRecursion} with memory as all lower order trees are computed already.

Note that \cref{theorem:generalLinearTreeRecursion} could be defined more generally however, for the results to be equivalent for arbitrary padding then a choice of the \(\emptyset\) interaction must be defined to \(0\) for either binary operation input. In our case the natural choice is \(\emptyset \star \gforest = 0\).

%!TEX root = ../main_theory.tex
%%%%%%%%%%%%%%%%%%%%%%%%
\section{\BSeriesText Composition\label{section:Composition}}
%%%%%%%%%%%%%%%%%%%%%%%%

We now formally define the initial conditions of \eqref{eq:ApproximationProblem} as 
\begin{equation}
    \vg\comp{i}(\tau_\initstate) = \NBSeries{\Dt}(b, \y_{\initstate})
\end{equation}
where the labels \(b(\tree) = 0\) if the tree \(\tree\) contains any colors of \(\gun\comp{j}\,\forall\,j\). The series comprising the ODE system needs to be restated in terms of the initial condition \(\y_{\initstate}\). Thus, to handle initial conditions this requires the use of \BSeriesText \textit{Composition} \cite{Chartier_2010_Bseries,Hairer_2006_geometric-book,Butcher_2021_Bseries-book}.

We provide a new definition of tree pruning to support composition by utilizing the permanent to avoid using \textit{ordered} trees. Additionally, triangular matrix analysis is provided in the appendix to support the block matrix formulation of composition. Such analysis provides a formulation for repeated composition, properties of the exact labels, and an alternative way to get the \textit{backward differential equation} labels \cite{Hairer_2006_geometric-book}.

%%%%%%%%%%%%%%%%%%%%%%%%
\subsection{Tree Pruning}
%%%%%%%%%%%%%%%%%%%%%%%%
To start, the concept of tree pruning is introduced. Throughout literature, it is typical to define pruning in terms of ordered trees \cite{Hairer_2006_geometric-book, Chartier_2010_Bseries,Chartier_2005_substitution,Butcher_2021_Bseries-book,hairerOrderConditionsNumerical1981,connesHopfAlgebrasRenormalization1998}.

For this work we define pruning via recursion over only \textit{unordered} trees and extend it to forests.
\begin{definition}[Pruning]\label{def:pruning}
    For the trees $\tree = \tau[\forestTQt] \in \NTTC$ and $\stree = \psi[\forestTQs]  \in \NTTC$ where \(\tau\) and \(\psi\) are nodes tree pruning is defined as
    \begin{align}
        \tau[\forestTQt] \prune \psi[\forestTQs] &= \scalemath{1.0}{\begin{cases}
            \forestTQt \prune \forestTQs \textbf{ if } \tau = \psi, m,n > 0\\
            \forestTQt \prune \emptyset \textbf{ if } \tau = \psi, n > 0, m = 0\\
            0 \textbf{ otherwise}
        \end{cases}}, \\
        \tree \prune \emptyset &= \tree,\\
        \emptyset \prune \stree &= 0,\\
        \emptyset \prune \emptyset &= 1,\\
        \tree \prune \emptyset_{c} &= \begin{cases}
            \tree \quad \text{ if } \treeRoot(\tree) = c\\
            0 \quad \text{ otherwise}
        \end{cases}.
    \end{align}
    For the forests $\forestTQt \in \ForestSpace{\NTTC}$ and $\forestTQs \in \ForestSpace{\NTTC}$  pruning for them is defined as 
    \begin{equation}
        \begin{aligned}
            \forestTQt \prune \forestTQs &= \scalemath{1.0}{\begin{cases}
                \frac{1}{(\prod_{i=1}^m k_i !)p!}\permanent(\forestTQt \outergen{\prune} \forestTQs \emptyset^{p}) \textbf{ if } p \ge 0\\
                0 \textbf{ otherwise}
            \end{cases}},\\
            \forestTQt \prune \emptyset &= \forestTQt
        \end{aligned}
    \end{equation}
    where $p = \sum_{i=1}^n q_i - \sum_{i=1}^m k_i$.
\end{definition}
\begin{theorem}[Pruning Elementary Differentials]\label{theory:pruning}
    Given two trees \(\tree, \stree \in \NTTC\) and a \BSeriesText \(\y^{*} = \NBSeries{\Dt}(a, \y_\initstate)\) then 
    \begin{equation}
        \frac{\BDiffOp(\stree)(\y^{*})}{\sigma(\stree)} = \sum_{\tree \in \NTTC} a(\tree \prune \stree) \frac{\BDiffOp(\tree)(\y_\initstate)}{\sigma(\tree)}.
    \end{equation}
    In different terms, \(\tree \prune \stree\) corresponds to the leading label of \(\BDiffOp(\tree)(\y_\initstate)\) obtained from differentiating \(\BDiffOp(\stree)(\y^{*})\).
\end{theorem}
A proof of this is provided in \cref{proof:pruning}.
\begin{example}
    A simple, yet extensive, example for multicolored trees is
    \begin{equation}
        \btree{[x[x[z[x]]][x[z][x]][x][x]]} \prune \btree{[x[x[z]][x]]} = \btree{[x]}^4\btree{[z]} + 2\btree{[x[z[x]]]}\,\btree{[x]}^2 + 2\btree{[x[z][x]]}\,\btree{[x]}^2 + \btree{[z[x]]}\,\btree{[x]}^3.
    \end{equation}
\end{example}
\begin{remark}[Pruning Linear Recursion]\label{remark:pruninglinear}
    Using the same notations in \cref{theory:pruning} and base conditions the operation can be further reduced from \cref{theorem:generalLinearTreeRecursion}. If \(\bar{n} \ge \bar{m} \ge 1\) then
    \begin{equation}
        \begin{aligned}
            \tree \prune \stree &= \delta_{\tau, \psi} \Big( \sum_{j=1}^{m} (\tree_i \prune \stree_j)\left(\tau\left[\frac{\forestTQt}{\tree_i}\right] \prune \psi\left[\frac{\forestTQs}{\stree_j}\right] \right)  \\
                    &+ \tree_i\left(\tau\left[\frac{\forestTQt}{\tree_i}\right] \prune \stree \right) \Big)\textit{ for some } i \in \{1 .. n\}.\\
        \end{aligned}
    \end{equation}
    If \(\bar{n} < \bar{m}\) then \(\tree \prune \stree = 0\). Thus, pruning of two trees can be constructed by product sum of other pruning of trees. This form can be less memory intensive than \cref{theorem:generalLinearTreeRecursion} as only the possible operations on trees are needed instead of all forest. We exclude adding the alternative form for children of \(\stree\). We note that this form is the same as in \cite{muruaFormalSeriesNumerical1999} via the Butcher product.
\end{remark}
The linear formulation is preferred when utilized for composition. In practice the construction of lower order tree pruning are needed anyway thus, utilizing memory reduces computation down to only a sum over the distinct children of the root of \(\stree\). Additional simplifications can be applied to the formula to give better performance.
\begin{example}
    Taking the example from before and choosing the subtree \(\btreega{[x[z[x]]]}\) gives
    \begin{equation}
        \btreega{[x[x[z[x]]][x[z][x]][x][x]]} \prune \btreega{[x[x[z]][x]]} = \left(\btreega{[x[z[x]]]} \prune \btreega{[x]}\right)\left(\btreega{[x[x][x][x[x][z]]]} \prune \btreega{[x[x[z]]]}\right) + \left(\btreega{[x[z[x]]]} \prune \btreega{[x[z]]}\right)\left(\btreega{[x[x][x][x[x][z]]]} \prune \btreega{[x[x]]}\right) + \btreega{[x[z[x]]]} \left(\btreega{[x[x][x][x[x][z]]]} \prune \btreega{[x[x[z]][x]]}\right).
    \end{equation}
\end{example}
 \begin{remark}[Alternative Pruning Definitions]
    There are alternative definitions in literature.
    \begin{itemize}
        \item \cite{Chartier_2010_Bseries, Hairer_2006_geometric-book} Work with the set of ordered partitions. We avoid such formulations given the usage of ordered trees.
        \item \cite[Section 2.8]{Butcher_2021_Bseries-book} This formulation focuses on forming the full sum via algebraic relations instead of the individual pruning terms.  
        \item \cite{butcherBseriesCompositionTheorem2025} The newest derivation takes counting formulation via sub-graphs.
    \end{itemize}
 \end{remark}
To help simplify results in further sections some additional subsets are provided. The subsets are formally defined by the prune operation.
\begin{definition}[Tree Pruning Subset \cite{Butcher_2021_Bseries-book}]\label{def:prunesub}
    Given a tree \(\tree\) the subset \(\stree < \tree\) is the set of all trees \(\stree\) such that \(\tree \prune \stree \ne 0\) not including \(\tree\). 
\end{definition}
\begin{definition}[Monotonic Leaf Set]\label{def:monoLeafSet}
    Let \(\mathcal{L}(\psi,\,\mathcal{N})\) be the set of trees where all trees of a single color \(\psi\) must have at least one leaf of a different color from the set of colors \(\mathcal{N}\). The formal definition of such set is
    \begin{equation}\label{eq:treeCL}
        \begin{aligned}
            &\mathcal{L}(\psi,\mathcal{N}) = \Big\{\stree\,|\,\stree \in \NTTi{\{\psi\}} \text{ and }\, \exists\, \rtree \in \NTTi{\{\treeRoot(\stree)\}}\quad \\
            &\text{s.t.}\quad \stree \prune \rtree = \prod_{i=1}^{m} c_i\, \qtree_{i}^{k_i},\, \treeRoot{\qtree_i} \ne \psi ,\,  \treeorder{\qtree_i} = 1,\, m > 0,\, k_i > 0,\, c_i \in \mathbb{R},\, \text{and}\, \qtree_i \in \NTTi{\{\treeRoot(\qtree_i)\}} \subset \NTTi{\mathcal{N}}\Big\}.
        \end{aligned}
    \end{equation}
\end{definition}
%
%%%%%%%%%%%%%%%%%%%%%%%%%%%
\subsection{Composition of \BSeriesText}
%%%%%%%%%%%%%%%%%%%%%%%%%%%
%
With the operation of pruning defined this immediately leads to its application in \BSeriesText composition.
\begin{theorem}[\BSeriesText Composition \cite{Chartier_2010_Bseries, Hairer_2006_geometric-book, Butcher_2021_Bseries-book, hairerOrderConditionsNumerical1981, hairer1974butchergruppe}]\label{theory:BSeriesComposition}
    Given $\NBSeries{\Dt}(b,\Tilde{\y})$ where the initial condition is $\Tilde{\y} = \BSeries{\Dt}(a,\y_\initstate)$ then composition reduces the expression to
    \begin{equation}
        \NBSeries{\Dt}(b,\NBSeries{\Dt}(a,\y_\initstate)) = \NBSeries{\Dt}(b \labelComposition a,\y_\initstate).
    \end{equation}
    The labels $(b \labelComposition a)$ are defined as
    \begin{equation}
        (b \labelComposition a)(\tree) = \sum_{\stree \in \NTTC} a(\tree \prune \stree) b(\stree)
    \end{equation}
    where \(b(\emptyset_{c}) = 1\) for all colors.
\end{theorem}
\begin{proof}
    Applying \cref{theory:pruning} and rearranging terms gives
    \begin{equation}
        \begin{split}
            &\NBSeries{\Dt}(b,\Tilde{\y}) = \sum_{\tree \in \AbstractT} b(\tree) \Dt^{\treeorder{\tree}} \frac{\BDiffOp(\tree)(\tilde{\y})}{\sigma(\tree)} = \sum_{\tree \in \AbstractT} \sum_{\stree \in \AbstractT} b(\tree) a(\tree \prune \stree)\Dt^{\treeorder{\stree}}\frac{\BDiffOp(\stree)(\tilde{\y})}{\sigma(\stree)}\\
            &= \sum_{\stree \in \AbstractT} (b \labelComposition a)(\stree) \Dt^{\treeorder{\stree}}\frac{\BDiffOp(\stree)(\tilde{\y})}{\sigma(\stree)} = \NBSeries{\Dt}(b \labelComposition a,\y_\initstate).
        \end{split}
    \end{equation}
\end{proof}
Alternative proofs of composition can be found in \cite{Butcher_2021_Bseries-book,Chartier_2010_Bseries, Hairer_2006_geometric-book}. An alternative perspective of composition is presented in \cite{Butcher_2021_Bseries-book}[Section 2.9] via a linear algebra style for the composition operator.
\begin{definition}[\BSeriesText Composition Matrix]\label{def:pruneMatrix}
    Let  \(\treevecinf\) be an infinite vector where each entry is a tree from \(\NTTC\). For simplicity this vector is assumed to be ordered by tree order with lower tree orders at lower indices. The matrix $\pruneMatrix$ is defined as
    \begin{equation}
        \pruneMatrix =  \treevecinf \outergen{\prune} \treevecinf.
    \end{equation}
    The matrix can be indexed similarly to matrices by trees instead as \(\pruneMatrix_{\tree, \stree}\). If a set of integer indices \(i,j\) is provided instead then block indexing is done over tree order. This gives \(\blockIndex{\pruneMatrix}{i,j} = \blockIndex{\treevecinf}{i} \outergen{\prune} \blockIndex{\treevecinf}{j}\) where the sub-vector \(\blockIndex{\treevecinf}{i}\) is a vector restricted to trees satisfying \(\treeorder{\tree} = i\).

    Using the inner product form then \cref{theory:BSeriesComposition} is rewritten in matrix form as
    \begin{equation}
        \NBSeries{\Dt}(b \labelComposition a,\y_\initstate) = \left( \BDiffOp(\treevecinf, \y_\initstate) \frac{\Dt^{\treeorder{\treevecinf}}}{\sigma(\treevecinf)} \right)^{T}  a(\pruneMatrix)\,  b(\treevecinf).
    \end{equation}
\end{definition}
\begin{example}[Composition]
    Given two colors \(\btree{[gth]}, \btree{[gthr]}\) let the compressed color be \(\btree{[d]}\). Given the labels \(a\) and \(b\) then composition applied to the tree \(\btreega{[x[d][d]]}\) gives 
    \begin{equation}
        (b \labelComposition a)\left(\btreega{[x[d][d]]}\right) = a\left(\btreega{[x[d][d]]}\right) + a\left(\blockIndex{\pruneMatrix}{\btreega{[x[d][d]]},\, \btreega{[x[d]]}}\right) b\left(\btreega{[x[d]]}\right) + a\left(\blockIndex{\pruneMatrix}{\btreega{[x[d][d]]},\, \btreega{[x]}}\right) b\left(\btreega{[x]}\right) + b\left(\btreega{[x[d][d]]}\right),
    \end{equation}
    and in the expanded form
    \begin{equation}
        \scalemath{1.0}{
                \begin{bmatrix}
                    (b \labelComposition a)\left(\btreega{[x[gth][gth]]}\right)\\
                    (b \labelComposition a)\left(\btreega{[x[gth][gthr]]}\right)\\
                    (b \labelComposition a)\left(\btreega{[x[gthr][gthr]]}\right)\\
                \end{bmatrix} = \begin{bmatrix}
                    a\left(\btreega{[x[gth][gth]]}\right)\\
                    a\left(\btreega{[x[gth][gthr]]}\right)\\
                    a\left(\btreega{[x[gthr][gthr]]}\right)\\
                \end{bmatrix} + \begin{bmatrix}
                    2a\left(\btreega{[gth]}\right) & 0 \\
                    a\left(\btreega{[gth]}\right) & a\left(\btreega{[gthr]}\right) \\
                    0 & 2a\left(\btreega{[gthr]}\right)
                \end{bmatrix}\begin{bmatrix}
                    b\left(\btreega{[x[gth]]}\right)\\
                    b\left(\btreega{[x[gthr]]}\right)
                \end{bmatrix} +\begin{bmatrix}
                    a\left(\btreega{[gth]}\right)^2 \\
                    a\left(\btreega{[gth]}\right)a\left(\btreega{[gthr]}\right)\\
                    a\left(\btreega{[gthr]}\right)^2
                \end{bmatrix}b\left(\btreega{[x]}\right) + 
                \begin{bmatrix}
                    b\left(\btreega{[x[gth][gth]]}\right)\\
                    b\left(\btreega{[x[gth][gthr]]}\right)\\
                    b\left(\btreega{[x[gthr][gthr]]}\right)\\
                \end{bmatrix}
        }.
    \end{equation}
\end{example}

%%%%%%%%%%%%%%%%%%%%%%%%%%%
\subsection{Repeated Composition\label{subsection:RepeatedComposition}}
%%%%%%%%%%%%%%%%%%%%%%%%%%%
%
Properties are established utilizing the matrix derivations in \cref{section:MatId} for the composition matrix (\cref{def:pruneMatrix}). For this section let  \(D_{\alpha}\) be defined as follows:
\begin{equation}
    \DiagMat_{\alpha} = \textrm{Diag}([\underbrace{1, \dots, 1}_{|\treevecinf_0|}, \underbrace{\alpha, \dots, \alpha}_{|\treevecinf_1|}, \underbrace{\alpha^2, \dots, \alpha^2}_{|\treevecinf_2|}, \dots, \underbrace{\alpha^i, \dots, \alpha^i}_{|\treevecinf_i|}, \dots]).
\end{equation}

A time scaling of the composition operator is done similarly as \cref{property:timescale}.
\begin{property}[Composition Time Scaling] \label{property:timeScaling}
    Given some label \(a\) and \(\tilde{a}(\tree) = \alpha^{\treeorder{\tree}} a(\tree)\) then
    \begin{equation}
        \tilde{a}(\pruneMatrix) = \DiagMat_{\alpha}\, a(\pruneMatrix)\, \DiagMat_{\alpha}^{-1}
    \end{equation}
    since each block \(\blockIndex{\pruneMatrix}{i,j}\) must contain forest of order \(\max(i - j,\, 0)\).
\end{property}
\begin{theorem}[Pruning Properties] \label{theory:pruneProperties}
    All conditions stated in \cref{theory:TriangleMatrix} hold for the matrix \(\frac{1}{\gammaFun{\Tcolors}}(\pruneMatrix)\) where \(\Tcolors\) is some set of colors. More specifically the following hold
    \begin{equation}
        \begin{aligned}
            \DiagMat_{\Dt}\, \frac{1}{\gammaFun{\Tcolors}}(\pruneMatrix)\, \DiagMat_{\Dt}^{-1} &= \frac{1}{\gammaFun{\Tcolors}}(\pruneMatrix)^{\Dt},\\
            \frac{1}{\gammaFun{\Tcolors}}(\blockIndex{\pruneMatrix}{n,n-k}) &= \frac{1}{k!} \prod_{i=0}^{k-1} \frac{1}{\gammaFun{\Tcolors}}(\blockIndex{\pruneMatrix}{n - i, n - i - 1}) \qquad \forall\, n > k,\\
            \blockIndex{\log(\frac{1}{\gammaFun{\Tcolors}}(\pruneMatrix))}{i,j} &= \begin{cases}
                \frac{1}{\gammaFun{\Tcolors}}(\blockIndex{\pruneMatrix}{i,j}) & \textit{ if } i = j - 1,\\
                \matzero & \text{o.w.}.
            \end{cases}
        \end{aligned}
    \end{equation}
\end{theorem}
A proof is provided in \cref{proof:pruneProperties}.
Additionally, we provide an expansion of the composition operator for when a numerical method satisfies order conditions up to \(p\).
\begin{theorem}[Repeated Steps of an Integration Method]\label{theory:repeatedIntegration} 
    Given \cref{theory:TriangleMatrixClose} and \cref{theory:pruneProperties} and \(m \in \outspace\) then
    \begin{equation}
        (\DiagMat_{m})^{-1}\pruneMatrix^{m} \DiagMat_{m} = \Identity + \CompLogBlock_1 + \sum_{k=1}^{\infty} \frac{\CompLogBlock_{k+1}}{m^k}.
    \end{equation}
    The matrices \(\CompLogBlock_{k+1}\) are defined in \cref{theory:TriangleMatrixClose}.
    Given some \BSeriesText method with labels \(a\) where \(a(\tree) = \frac{1}{\gammaFun{\Tcolors}(\tree)} \quad \forall\, 1 \le \treeorder{\tree} \le p\) then the method repeated \(m\) uniform times over a fixed interval \([t_\initstate, t_\finalstate]\) is represented with labels denoted as \(a^{m}(\tree)\) where their expansion is
    \begin{equation}
            a^{m}(\pruneMatrix) - \frac{1}{\gammaFun{\Tcolors}}(\pruneMatrix) = \sum_{k=p}^{\infty} \frac{a(\CompLogBlock_{k+1})}{m^k}.
    \end{equation}
\end{theorem}
\begin{proof}
    Given \cref{theory:pruneProperties} and \cref{theory:TriangleMatrixClose} this comes directly since all prune operations are directly scaled by \(\DiagMat_{m}\) for the order of the tree.
\end{proof}
This gives a description of how coefficients converge under repeated composition thus, expanding the composition operator by arbitrary powers.
It is noted that work in \cite{shmelev2025explicit} arrived at results for arbitrary real steps \(m \in \outspace\) differently via Hopf algebra.
Written in block matrix form the expansion from \cref{theory:repeatedIntegration} gives
\begin{equation}
    \scalemath{1.0}{
    \begin{aligned}
        (\DiagMat_{m})^{-1}\pruneMatrix^{m} \DiagMat_{m} &= \begin{bmatrix}
            \Identity & \matzero & \matzero & \matzero & \dots\\
            \blockIndex{\pruneMatrix}{1,0} & \Identity & \matzero & \matzero & \dots\\
            \frac{\blockIndex{\pruneMatrix}{2,1}\blockIndex{\pruneMatrix}{1,0}}{2!} & \blockIndex{\pruneMatrix}{2,1} & \Identity & \matzero & \dots\\
            \frac{\blockIndex{\pruneMatrix}{3,2}\blockIndex{\pruneMatrix}{2,1}\blockIndex{\pruneMatrix}{1,0}}{3!} & \frac{\blockIndex{\pruneMatrix}{3,2}\blockIndex{\pruneMatrix}{2,1}}{2!} & \blockIndex{\pruneMatrix}{3,2} & \Identity & \dots\\
            \vdots & \vdots & \vdots & \vdots & \ddots
        \end{bmatrix}\\
        &+ \frac{1}{m}\begin{bmatrix}
            \matzero & \matzero & \matzero & \matzero & \dots\\
            \matzero & \matzero & \matzero & \matzero & \dots\\
            \blockIndex{\pruneMatrix}{2,0}  - \frac{\blockIndex{\pruneMatrix}{2,1}\blockIndex{\pruneMatrix}{1,0}}{2} & \matzero & \matzero & \matzero & \dots\\
            \frac{-\blockIndex{\pruneMatrix}{3,2}\blockIndex{\pruneMatrix}{2,1}\blockIndex{\pruneMatrix}{1,0}}{2} + \frac{1}{2}\left(\blockIndex{\pruneMatrix}{3,2}\blockIndex{\pruneMatrix}{2,0} + \blockIndex{\pruneMatrix}{3,1}\blockIndex{\pruneMatrix}{1,0}\right) & \blockIndex{\pruneMatrix}{3,1}  - \frac{\blockIndex{\pruneMatrix}{3,2}\blockIndex{\pruneMatrix}{2,1}}{2} & \matzero & \matzero & \dots\\
            \vdots & \vdots & \vdots & \vdots & \ddots
        \end{bmatrix} \\
        &+ \frac{1}{m^2}\begin{bmatrix}
            \matzero & \matzero & \matzero & \matzero & \dots\\
            \matzero & \matzero & \matzero & \matzero & \dots\\
            \matzero & \matzero & \matzero & \matzero & \dots\\
            \frac{\blockIndex{\pruneMatrix}{3,2}\blockIndex{\pruneMatrix}{2,1}\blockIndex{\pruneMatrix}{1,0}}{3} - \frac{1}{2}\left(\blockIndex{\pruneMatrix}{3,2}\blockIndex{\pruneMatrix}{2,0} + \blockIndex{\pruneMatrix}{3,1}\blockIndex{\pruneMatrix}{1,0}\right) + \blockIndex{\pruneMatrix}{3,0} & \matzero & \matzero & \matzero & \dots\\
            \vdots & \vdots & \vdots & \vdots & \ddots
        \end{bmatrix} + \dots .
    \end{aligned}
    }
\end{equation}
Further terms can be constructed utilizing the block form definition of \(\CompLogBlock_{k+1}\). Note this definition is independent of the colors used. When some label \(a\) is applied then the matrices of order \(\sfrac{1}{m^p}\) start to become zero matrices depending on the order of the method. Note that the first column gives the error constants of an arbitrary time integration method. Since the coefficients are parameterized by the amount of uniform steps, then by \cref{theory:TriangleMatrixClose} the limit matrix converges to the exact composition operator.

It is noted that the matrix analysis provided has corresponding relationship to that of Pascal's matrix \cite{aceto2001matrices,edelman2004pascal} since the composition matrix \(\frac{1}{\gammaFun{\Tcolors}}(\pruneMatrix)\) is Pascal's triangle when the differential equation is independent of state. As noted before these results are only for uniform steps.
\begin{remark}[Non-uniform Steps]\label{remark:nonuniform}
    Given some method with labels \(a \ne \frac{1}{\gammaFun{\Tcolors}}\) that takes non-uniform steps there is no fractional power of \(a(\pruneMatrix)\) that represents this error. As an example let \(\Dt_1 \ne \Dt_2\) and \(\Dt_1 + \Dt_2 = 1\) then
    \begin{equation}
        \begin{aligned}
            &\forall\, a \ne \frac{1}{\gammaFun{\Tcolors}} \quad \nexists\, m \in \outspace \textit{ s.t. }  \DiagMat_{\Dt_1} \, a(\pruneMatrix)\, \DiagMat_{\Dt_1}^{-1}\, \DiagMat_{\Dt_2}\, a(\pruneMatrix)\, \DiagMat_{\Dt_2}^{-1} = \DiagMat_{m}\, a(\pruneMatrix)^{m}\, \DiagMat_{m}^{-1}.
        \end{aligned}
    \end{equation}
    To get around this utilizing the minimum and maximum steps can provide approximate results. An alternative choice is to use \(m = \left(\sum_{i=1}^{n} \Dt_{i}^{p+1}\right)^{-1/p}\) to handle only the \(p+1\) order terms. Higher terms than \(p+1\) will not be correct necessarily. In summary, uniform stepping can only act as a proxy for non-uniform stepping.
\end{remark}

It is noted if instead the labels \(a = \frac{1}{\gammaFun{\tilde{\Tcolors}}}(\pruneMatrix)\) where \(\tilde{\Tcolors} \subset \Tcolors\) then the results will still converge; however, they converge to solutions of the sub-colors instead. As an example let there be two colors represented by \(\fun\comp{1}\) and \(\fun\comp{2}\). Given the differential equation
\begin{equation}
    \ddifffrac{t}\y = \fun\comp{1} + \fun\comp{2}
\end{equation}
then if Lie-Trotter is used to solve this problem this equates to the following:
\begin{equation}
    \begin{aligned}
        \ddifffrac{t}{\y\comp{1}(t)} &= \fun\comp{1}(\y\comp{1}(t)),\\
        \y\comp{1}(t_\initstate) &= \y_\initstate,\\
    \end{aligned} \qquad 
    \begin{aligned}
        \ddifffrac{t}{\y\comp{2}(t)} &= \fun\comp{2}(\y\comp{2}(t)),\\
        \y\comp{2}(t_\initstate) &= \y\comp{1}(t_\finalstate).\\
    \end{aligned}
\end{equation}
If some sub-methods are used to solve the problems with uniform steps \(m_1\) and \(m_2\) then utilizing \cref{theory:repeatedIntegration} the final formal coefficients of Lie-Trotter with sub-integrators are
\begin{equation}
    \scalemath{0.95}{
    \begin{aligned}
        a^{m_1}_{1}(\pruneMatrix)\,a^{m_2}_{2}(\pruneMatrix) &= \left(\frac{1}{\gammaFun{\{c_1\}}}(\pruneMatrix) + \sum_{k_1=1}^{\infty} \frac{a_1(\CompLogBlock_{k_1+1})}{m_1^k}\right)\left(\frac{1}{\gammaFun{\{c_2\}}}(\pruneMatrix) + \sum_{k_2=1}^{\infty} \frac{a_2(\CompLogBlock_{k_2+1})}{m_2^k}\right)\\
        &= \frac{1}{\gammaFun{\{c_1\}}}(\pruneMatrix)\frac{1}{\gammaFun{\{c_2\}}}(\pruneMatrix) + \frac{a_1(\CompLogBlock_{2})}{m_1} \frac{1}{\gammaFun{\{c_2\}}}(\pruneMatrix) + \frac{1}{\gammaFun{\{c_1\}}}(\pruneMatrix) \frac{a_2(\CompLogBlock_{2})}{m_2} + \textit{ H.O.T} .
    \end{aligned}
    }
\end{equation}
This produces the local label effects under arbitrary amount of sub-steps and order. However, increasing the steps or order has no impact on the coupling error since they only affect monotonic colored trees. The composition operator introduces the coupling trees or trees containing more than one color thus, acting as a coupling operation.
\begin{corollary}[\BSeriesText Composition Superset Representation]\label{corollary:BSeriesCompSuperSet}
    Let \(\NTTC\) be a tree set with colors \(\Tcolors\) where \(\NTTi{\Tcolors_1} \subseteq \NTTC\) and \(\NTTi{\Tcolors_2} \subseteq \NTTC\) are tree subsets where \(\Tcolors_1 \subseteq \Tcolors\) and \(\Tcolors_2 \subseteq \Tcolors\) are their respective colors.\\
    Let \(\tree \in \NTTi{\Tcolors_2 < \Tcolors_1} = (\NTTi{\Tcolors_2} < \NTTi{\Tcolors_1}) \subseteq \NTTC\) be the set of trees such that there is a tree \(\stree \in \NTTi{\Tcolors_2}\) such that \(\tree \prune \stree \in \ForestSpace{\NTTi{\Tcolors_1}}\).
    For the labels \(a_i(\tree) = 0 \quad \, \forall\, \tree \notin \NTTi{\Tcolors_i}\) then applying composition gives 
    \begin{equation}
        (a_2 \labelComposition a_1)(\tree) = 0 \quad \, \forall\, \tree \notin \NTTi{\Tcolors_2 < \Tcolors_1}.
    \end{equation}
\end{corollary}
This demonstrates that decoupled sub-systems need to eventually interact with the other sub-systems. Composition is one of the classical ways to re-couple sub-systems seen in operator splitting methods. However, relying on this alone can cause issues when targeting higher orders such as negative time steps \cite{Mclachlan_2002_splitting,blanes2024splitting}.

\subsection{Backward Error Block Formulation}
We provided an alternative derivation of backward error analysis \cite{Hairer_2006_geometric-book,hairerLifespanBackwardError1997} based on the block matrix formulations of composition. The backward error problem is defined as an exact solution of a differential equation which gives exactly the same numerical solution. Such a problem is given by some arbitrary labels \(a\) as follows (assuming \(a(\emptyset) = 1\)):
\begin{equation}\label{eq:backwardErrorProblem}
    \frac{\ddiff}{\ddiff \tau}\tilde{\y}(\tau) = \sum_{\tree \in \TT} m^{-\treeorder{\tree} + 1} \frac{(\log(a))(\tree) \, \Dt^{\treeorder{\tree} - 1} \, \BDiffOp(\tree)(\tilde{\y}(\tau)) }{\sigma(\tree)}
\end{equation}
where \(m\) is the amount of uniform steps of the method.
\begin{theorem}[Backward Error Labels]\label{theorem:backwardErrorLabels}
    Given a set of labels \(a\) where \(a(\emptyset) = 1\). Then the backward error problem \cref{eq:backwardErrorProblem} has the following coefficients
    \begin{equation}
        \log(a)(\treevecinf) = a(\log(\pruneMatrix))\, \elemVec{\emptyset}
    \end{equation}
    where \(\log(\pruneMatrix)\) is the matrix \(\log\) of the composition matrix \(\pruneMatrix\) and \(\elemVec{\emptyset}\)is the unit vector at the tree \(\emptyset\).
\end{theorem}
\begin{proof}
    Given the \BSeriesText with arbitrary uniform composition of coefficients \(a\) from \cref{theory:repeatedIntegration} is rewritten as 
    \begin{equation}
        \y_\finalstate = \y_\initstate + \sum_{\tree \in \TT} \frac{a(\tree) \, \Dt^{\treeorder{\tree}} \BDiffOp(\tree)(\y_\initstate)}{\sigma(\tree)} =  \y_\initstate + \sum_{\tree \in \TT} \sum_{k=0}^{\infty} m^{-k} \frac{a(\CompLogBlock_{k+1})(\tree)\, \Dt^{\treeorder{\tree}}\, \BDiffOp(\tree)(\y_\initstate) }{\sigma(\tree)}.
    \end{equation}
    Note that though the sum of \(k\) is taken to infinity the coefficient \(a(\CompLogBlock_{k+1})(\tree)\) are \(0\) when \(k > \treeorder{\tree}\).
    Rearranging terms gives the numerical series 
    \begin{equation}
        \y_\finalstate = \y_\initstate + \sum_{k=0}^{\infty} \Dt^{k} \sum_{\tree \in \TT}  m^{-(\treeorder{\tree} - k) + 1} \frac{a(\CompLogBlock_{\treeorder{\tree} - k})(\tree)\, \Dt^{\treeorder{\tree} - k}\, \BDiffOp(\tree)(\y_\initstate) }{\sigma(\tree)}.
    \end{equation}
    The series is now collected in terms such that the Taylor series \(\y_\finalstate = \y_\initstate + \Dt\frac{\ddiff}{\ddiff t}\tilde{\y}(t_\initstate) + \frac{\Dt^2}{2!}\frac{\ddiff^2}{\ddiff t^2}\tilde{\y}(t_\initstate) + \dots\) has been formed. This is given by the increasing \(\log\) powers of \(\CompLogBlock\) and the increasing derivatives of higher terms. Thus, taking the first term gives the backward error differential equation as 
    \begin{equation}
        \frac{\ddiff}{\ddiff\tau}\tilde{\y}(\tau) = \sum_{\tree \in \TT} m^{-\treeorder{\tree} + 1} \frac{a(\CompLogBlock_{|\tree|})(\tree) \, \Dt^{\treeorder{\tree} - 1} \, \BDiffOp(\tree)(\tilde{\y}(\tau)) }{\sigma(\tree)}
    \end{equation}
    where \(\tau \in [t_\initstate, t_\finalstate]\) and \(\tilde{\y}(t_\initstate) = \y_\initstate\). Thus, by the definition of \(\CompLogBlock_{|\tree|}\) then the statement holds.
\end{proof}

Repeated integration is given as extending the time domain of the backward error problem. Thus, the use of the matrix \(\log\) of the prune matrix implies the use of the backward error. This gives an alternative formulation to modified vector fields \cite{Chartier_2010_Bseries,Hairer_2006_geometric-book}. 

\section{Replacement/Substitution}
\label{section:Replacement}

This section will introduce a new binary operation on trees/forests called ``replacement''. The overall goal of replacement is to replace a subset of nodes of colors \(\mathcal{U}\) in a tree such that it matches another tree with no colors of \(\mathcal{U}\). A binary operation is utilized to specify all possible transformations which produce the result; such as with pruning. This provides an operation assuming the elementary differentials associated with \(\mathcal{U}\) can be constructed via a \BSeriesText without colors \(\mathcal{U}\).

For the case of the approximate problem \cref{eq:ApproximationProblem} we now provide the definition of the approximate RHS (without time dependent terms) as
\begin{equation}\label{eq:modifiedAsymptoticConstant}
    \begin{split}
        \gun(\vg(t)) &= \sum_{p \in (\mathbb{Z}^{+})^{\nparts}} \frac{\basisK_{p}}{\prod_{j=1}^{\nparts} p_{j}!} \,\prod_{j=1}^{\nparts} \,(\vg\comp{j}(t) - \vg\comp{j}(\tau_\initstate))^{p_{j}},\\
        \dforest_{p} &= \prod_{j=1}^{\nparts} \xnode_{c_{p_{j}}}^{p_{j}},\\
        \frac{\basisK_{p}}{\prod_{j=1}^{\nparts} p_{j}!} &= \NBSeries{\Dt}^{\dforest_{p}}(a^{[p]}, \vg(\tau_\initstate)) = \frac{\BDiffOp(\psi[\dforest_{p}])(\vg(\tau_\initstate))}{\sigma(\psi[\dforest_{p}])},
    \end{split}
\end{equation}
where \(\psi\) is some node that is to be replaced, and \(a^{[p]}\) is some label which is zero for any trees which contain node \(\psi\). We note that this is not the full asymptotic series we study for \cref{eq:ApproximationProblem} and additional time modifications are reserved for \cref{section:TFDEQ}.

Each approximate field basis \(\basisK_{p}\) is represented by a temporary colored node \(\psi\) in the expansion. To get the full expansion w.r.t. only the elementary differentials \(\fun\comp{i}\) a new operation called \textit{replacement} is introduced. Replacement's goal is to reduce the expression to a different space of elementary differentials.

An example of replacement is provided as follows. Let \(a_1, a_2\) be some set of labels on trees with elementary differentials in terms of \(\gun(\y,\z),\, \fun(\y,\z),\, \basisK_i\).  Let the elementary differentials \(\basisK_i\) have the expansion 
\begin{equation}
    \basisK_i = \begin{bmatrix}
            \gun & \gun_v\fun & \gun_w \gun & \dots
        \end{bmatrix} \begin{bmatrix}
            a_i(\btree{[z]})\\
            a_i(\btree{[z[x]]})\\
            a_i(\btree{[z[z]]})\\
            \vdots\\
        \end{bmatrix}.
\end{equation}
Expanding the elementary differential \(\fun_{\z,\z} (\basisK_1, \basisK_2)\) in terms of only \(\gun\) and \(\fun\) gives 
\begin{equation}
    \scalemath{0.85}{\fun_{\z,\z} (\basisK_1, \basisK_2) = \begin{bmatrix}
        \frac{1}{2}\fun_{\z,\z} (\gun, \gun) & \fun_{\z,\z} (\gun_{\y} \fun, \gun) & \fun_{\z,\z} (\gun_{\z} \gun, \gun) & \frac{1}{2} \fun_{\z,\z} (\gun_{\y} \fun, \gun_{\y} \fun) & \dots
    \end{bmatrix} \begin{bmatrix}
        2 a_1(\btree{[z]}) a_2(\btree{[z]})\\
        \left( a_1(\btree{[z[x]]}) a_2(\btree{[z]}) + a_1(\btree{[z]}) a_2(\btree{[z[x]]}) \right)\\
        \left( a_1(\btree{[z[z]]}) a_2(\btree{[z]}) + a_1(\btree{[z]}) a_2(\btree{[z[x]]}) \right)\\
        2 a_1(\btree{[z[x]]}) a_2(\btree{[z[x]]})\\
        \vdots\\
    \end{bmatrix}.}
\end{equation}
Thus, the action of \textit{replacement} reduced the representation down to only terms of \(\gun\) and \(\fun\).
The replacement operation is defined using the same techniques as composition with the permanent. This operation is similar to \BSeriesText \textit{substitution} \cite{Chartier_2005_substitution,Chartier_2010_Bseries,Hairer_2006_geometric-book}; however, key differences will be explained. 

\subsection{Marked Pruning}
To handle replacement additional operations and structures are to be introduced. To start a tree bracket is introduced to help pair recursive outputs.
\begin{definition}[Tree Bracket]\label{def:treebrack}
    The tree bracket is defined as 
    \begin{equation}
        \begin{aligned}
            &\labelbrack{\fforest}{\gforest}
        \end{aligned}
    \end{equation}
    where the forest \(\fforest\) is paired with forest \(\gforest\). This bracket comes equipped with the properties
    \begin{equation}
        \begin{aligned}
            \labelbrack{a_1\, \fforest_1 + a_2\, \fforest_2}{\gforest} &= a_1\,\labelbrack{\fforest_1}{\gforest} + a_2\,\labelbrack{\fforest_2}{\gforest},\\
            \tau \ButcherProduct (\labelbrack{\fforest}{\gforest} \hforest) &= \labelbrack{\fforest}{\gforest}  (\tau \ButcherProduct \hforest)
        \end{aligned}
    \end{equation}
    where \(a_1,a_2 \in \outspace\), \(\tree_i \in \NTTC\), \(\fforest_1, \fforest_2 \in \ForestSpace{\NTTC}\), and \(\psi\) is a node. The pairings allow for  multi-variable label evaluations and richer recursive structures. Extending the labels to handle pairs \(w: (\NTTC \times \NTTC) \rightarrow \outspace\) then the application to a tree bracket gives 
    \begin{equation}
        \begin{aligned}
            w(\labelbrack{\tree}{\stree}) &= w(\stree)(\tree),\\
            w(a_1\,\labelbrack{\fforest_1}{\gforest} + a_2\,\labelbrack{\fforest_2}{\hforest}) &= a_1\, w(\gforest)(\fforest_1) + a_2\, w(\hforest)(\fforest_2).\\
        \end{aligned}
    \end{equation}
    Lastly, the operations
    \begin{equation}
        \begin{aligned}
            \brackindexfirst{\labelbrack{\fforest}{\gforest}} &= \fforest,\\
            \brackindexsecond{\labelbrack{\fforest}{\gforest}} &= \gforest,\\
        \end{aligned}
    \end{equation}
    provide access to the forests inside the brackets.
\end{definition}
To support the new brackets an extension of Butcher grafting is applied to the right argument.
\begin{definition}[Pairs Right Grafting]\label{def:graftmark}
    The grafting with pairs on the right is defined as
    \begin{equation}
        \begin{aligned}
            \tau[\forestTQt]  \ProductMarked \labelbrack{\fforest}{\gforest} &= \labelbrack{\fforest}{\tau[\forestTQt \gforest]},\\
            \tau[\forestTQt]  \ProductMarked \prod_{i=1}^{n}\labelbrack{\fforest_i}{\gforest_i} &= \labelbrack{\prod_{i=1}^{n}\fforest_i}{\tau\left[\forestTQt \prod_{i=1}^{n}\gforest_i\right]},\\
            \tau[\forestTQt]  \ProductMarked \left(\sum_{i=1}^{n}a_i \labelbrack{\fforest_i}{\gforest_i}\right) &= \sum_{i=1}^{n} a_i \labelbrack{\fforest_i}{\tau[\forestTQt \gforest_i]},\\
            \tau[\forestTQt]  \ProductMarked a &= \labelbrack{a}{\tau[\forestTQt]},\\
            \tau[\forestTQt]  \ProductMarked 0 &= 0\\
        \end{aligned}
    \end{equation}
    where \(a,a_i \in \outspace\), \(\tau\) is a node, and \(\fforest, \gforest, \fforest_i,\gforest_i, \hforest_i \in \ForestSpace{\NTTC}\).
\end{definition}
\begin{example}[Example of Pair Right Grafting]
    An example of pair right grafting is given as  
    \begin{equation}
        \btree{[x]} \ProductMarked \left(\labelbrack{\btree{[x]}}{\btree{[x]}}\labelbrack{\btree{[x]}}{\btree{[e]}} + 2\labelbrack{\btree{[x[x]]}}{\btree{[x]}}+ \labelbrack{1}{\btree{[x[x]]}}\labelbrack{\btree{[x]}}{\btree{[z]}}\right) = \labelbrack{\btree{[x]}^2}{\btree{[x[x][e]]}} + 2\labelbrack{\btree{[x[x]]}}{\btree{[x[x]]}} + \labelbrack{\btree{[x]}}{\btree{[x[z][x[x]]]}}.
    \end{equation}
\end{example}
Utilizing the tree bracket and right grafting gives a new variation of pruning which produces more information about \textit{where} trimmings come from.
\begin{definition}[Marked Pruning]\label{def:pruningMarked}
    Given the trees $\tree = \tau[\forestTQt] \in \NTTC$ and $\stree = \psi[\forestTQs] \in \NTTC$ where \(\tau\) and \(\psi\) are nodes then the prunings between the trees are paired with the stumped locations they were originally from as tree bracket pairs. This is provided by the marked pruning operation defined as
    \begin{align}
        \tau[\forestTQt] \pruneMark \psi[\forestTQs] &= \scalemath{1.0}{\begin{cases}
            \tau \ProductMarked (\forestTQt \pruneMark \forestTQs) \textbf{ if } \tau = \psi, m,n > 0\\
            \tau \ProductMarked (\forestTQt \pruneMark \emptyset) \textbf{ if } \tau = \psi, n > 0, m = 0\\
            0 \textbf{ otherwise}
        \end{cases}}, \\
        \tree \pruneMark \emptyset &= \labelbrack{\tree}{\xnode_{\treeRoot(\tree)}},\\
        \emptyset \pruneMark \stree &= 0,\\
        \emptyset \pruneMark \emptyset &= \labelbrack{1}{\emptyset},\\
        \tree \pruneMark \emptyset_{c} &= \begin{cases}
            \labelbrack{\tree}{\xnode_{\treeRoot(\tree)}} \quad \textbf{ if } \treeRoot(\tree) = c\\
            0 \quad \textbf{ otherwise}
        \end{cases}.
    \end{align}
    Given the forests $\forestTQt \in \ForestSpace{\NTTC}$ and $\forestTQs \in \ForestSpace{\NTTC}$ and letting $p = \sum_{i=1}^n q_i - \sum_{i=1}^m k_i$ then marked pruning for forest is defined as 
    \begin{equation}
        \begin{aligned}
            \forestTQt \pruneMark \forestTQs &= \scalemath{1.0}{\begin{cases}
                \frac{1}{(\prod_{i=1}^m k_i !)p!}\permanent(\forestTQt \outergen{\pruneMark} \forestTQs \emptyset^{p}) \textbf{ if } p \ge 0\\
                0 \textbf{ otherwise}
            \end{cases}},\\
            \forestTQt \pruneMark \emptyset &= \prod_{i=1}^{n} \labelbrack{\tree_{i}}{\xnode_{\treeRoot(\tree_i)}}^{q_i}.
        \end{aligned}
    \end{equation}
\end{definition}
\begin{theorem}[Pruning Marked Elementary Differentials]\label{theory:pruningMarked}
    Given two trees \(\tree, \stree \in \NTTC\) and a \BSeriesText \(\y^{*} = \NBSeries{\Dt}(a, \y_\initstate)\) then for \(\tree \ProductMarked \stree\) the marked pruning, for when \(\stree \ne \emptyset\), gives
    \begin{equation}
        \frac{\BDiffOp(\stree)(\y^{*})}{\sigma(\stree)} = \sum_{\tree \in \NTTC} a(\tree \pruneMark \stree) \frac{\BDiffOp(\tree)(\y_\initstate)}{\sigma(\tree)} = \sum_{\ktree \in \stree\lceil \TStumpcolors \rceil} \stumpCountFactorial(\ktree) \sum_{\hforest \in \ForestSpace{\NTTC}_{\stumpCount(\ktree)}} a(\labelbrack{\hforest}{\ktree}) \frac{\BDiffOp(\ktree)(\y_\initstate)}{\sigma(\ktree)} \symprod \frac{\BDiffOp(\hforest)(\y_\initstate)}{\sigma(\hforest)}
    \end{equation}
    where, for shorthand, when \(\bar{p}=0\) then \(\hforest = 1\). Otherwise, for when \(\stree = \emptyset\) then 
    \begin{equation}
        \BDiffOp(\emptyset)(\y^{*}) = \sum_{\tree \in \NTTC} a(\tree \pruneMark \emptyset) \frac{\BDiffOp(\tree)(\y_\initstate)}{\sigma(\tree)} = \left(\sum_{\ktree \in \NTTC} a(\labelbrack{\ktree}{\xnode_{\treeRoot(\ktree)}}) \frac{\BDiffOp(\ktree)(\y_\initstate)}{\sigma(\ktree)}\right) + a(\labelbrack{1}{\emptyset})\BDiffOp(\emptyset)(\y_\initstate).
    \end{equation}
\end{theorem}
A proof of this is provided in \cref{proof:pruningMarked}. The first part is the same result as \cref{theory:pruning}. Note that though the marked pruning is defined in terms of the permanent it can still be reduced via \cref{theorem:generalLinearTreeRecursion}. We do not explicitly write out the derivation, but it is the method to use for implementations.
\begin{example}[Example of Marked Pruning]
    Some examples of marked pruning are
    \begin{equation}
        \begin{aligned}
            \btree{[x[x[x]][x]]} \pruneMark \btree{[x[x]]} &= \labelbrack{\btree{[x]}^2}{\btree{[x[x[e]][e]]}} + \labelbrack{\btree{[x[x]]}}{\btree{[x[x][e]]}},\\
            \btree{[x[x[x][x]][x[x]]]} \pruneMark \btree{[x[x[x]][x]]} &= \labelbrack{\btree{[x]}^2}{\btree{[x[x[e][e]][x[x]]]}} + 2\labelbrack{\btree{[x]}^2}{\btree{[x[x[x][e]][x[e]]]}}.\\
        \end{aligned}
    \end{equation}
    With the last example it is noted that if standard pruning is utilized, then \(\btree{[x[x[x][x]][x[x]]]} \prune \btree{[x[x[x]][x]]} = 3\,\btree{[x]}^2\). This leaves the locations of the trimmings, and the amount from them, indecipherable.
\end{example}

\subsection{Replacement Operation}
With the tree bracket and marked pruning the general replacement is defined.
\begin{definition}[Replacement]\label{def:treereplacement}
    Let \(\mathcal{U}\) be the set of colored nodes to replace for the trees $\tree = \tau[\fforest] = \tau[\forestTQt]  \in \NTTi{\mathcal{A}}$ and $\stree = \psi[\gforest] = \psi[\forestTQs] \in \NTTC$ where the set \(\NTTi{\mathcal{A}} \subset \NTTC\) under the restriction \(\mathcal{A} \cup \mathcal{U} = \Tcolors\) and \(\mathcal{A} \cap \mathcal{U} = \emptyset\). Then the collection of replacements of nodes \(\mathcal{U}\) on \(\stree\) is formally given by
    \begin{align}
        \scalemath{1.0}{\tree \replaceNode{\mathcal{U}} \stree =  \tau[\forestTQt]  \replaceNode{\mathcal{U}}  \psi[\forestTQs]} & \scalemath{1.0}{ = \begin{cases} 
            \forestTQt \replaceNode{\mathcal{U}} \forestTQs \\
            \qquad \textbf{ if }  \tau = \psi\\
            \sum_{\emptyset < \vtree < \tree}  \left((\tree  \pruneMark \vtree) \replaceNode{\mathcal{U}} \labelbrack{\forestTQs}{\psi}\right) \\
            \qquad \textbf{ if }  \psi \in \mathcal{U}\\
            0 \qquad \textbf{ otherwise or } \tree \in \NTTC \setminus \NTTi{\mathcal{A}} 
        \end{cases}},\\
        \tree \replaceNode{\mathcal{U}} \psi &= \begin{cases}
            \labelbrack{\tree}{\psi} \textbf{ if } \psi \in \mathcal{U}\\
            0 \textbf{ otherwise}
        \end{cases},\\
        \tau \replaceNode{\mathcal{U}} \psi[\forestTQs] &= 0,\\
        \tau \replaceNode{\mathcal{U}} \psi &= \begin{cases}
            1 \textbf{ if } \tau = \psi\\
            \labelbrack{\tau}{\psi} \textbf{ if }  \psi \in \mathcal{U}\\
            0 \textbf{ otherwise}
        \end{cases},\\
        \tree \replaceNode{\mathcal{U}} \emptyset &= 0,\\
        \emptyset \replaceNode{\mathcal{U}} \stree &= 0,\\
        \emptyset \replaceNode{\mathcal{U}} \emptyset &= 1.
    \end{align}
    Forest replacement is given as 
    \begin{equation}
        \scalemath{1.0}{
        \begin{aligned}
            \forestTQt \replaceNode{\mathcal{U}} \forestTQs &= \begin{cases}
                \frac{1}{\prod_{i=1}^m (k_i !) } \permanent(\forestTQt \outergen{\replaceNode{\mathcal{U}}} \forestTQs) \textbf{ if } \bar{n} = \bar{m}\\
                0 \textbf{ otherwise}
            \end{cases},\\
            \sum_{j=1}^{p}c_j \labelbrack{\hforest_j}{\utree_j} \replaceNode{\mathcal{U}} \labelbrack{\forestTQs}{\psi} &= \sum_{j=1}^{p} c_j\labelbrack{\utree_j}{\psi} (\hforest_j \replaceNode{\mathcal{U}} \forestTQs),
        \end{aligned}
        }
    \end{equation}
    where \(c_j \in \outspace\), \(\forestTQs \in \ForestSpace{\NTTC}\), \(\forestTQt, \fforest_j  \in  \ForestSpace{\NTTi{\mathcal{A}}}\).
\end{definition}
\begin{theorem}[Replacement Elementary Differentials]\label{theory:treereplacement}
    Given two trees \(\tree \in \NTTi{\mathcal{A}}\), \(\stree \in \NTTC\) and for all colors \(\psi \in \mathcal{U}\) \(\BDiffOp(\psi)(\y_\initstate) = \NBSeries{\Dt}(a(\psi), \y_\initstate)\) then 
    \begin{equation}
        \frac{\BDiffOp(\stree)(\y_\initstate)}{\sigma(\stree)} = \sum_{\tree \in \NTTi{\mathcal{A}}} a(\tree \replaceNode{\mathcal{U}} \stree) \frac{\BDiffOp(\tree)(\y_\initstate)}{\sigma(\tree)}.
    \end{equation}
    In different terms \(\tree \replaceNode{\mathcal{U}} \stree\) corresponds to the leading label of \(\BDiffOp(\tree)(\y_\initstate)\) obtained from \(\BDiffOp(\stree)(\y_\initstate)\) by replacing nodes \(\psi\) with their series expansion.
\end{theorem}
A proof of this has been provided in \cref{proof:treereplacement}. As before when implementing this it is recommended to use \cref{theorem:generalLinearTreeRecursion} as the number of possible trees can grow quite large for replacement. We omit such derivations for simplicity.
\begin{example}
    As an example given the tree \(\btree{[x[gth][gth][z[gth]]]}\) and the tree \(\btree{[x[z[x]][z[x]][z[x][x]]]}\) then the replacement of node \(\btree{[gth]}\) is given as
    \begin{equation}
        \btree{[x[z[x]][z[x]][z[x][x]]]} \replaceNode{\{\btree{[gth]}\}} \btree{[x[gth][gth][z[gth]]]} = \frac{1}{2!} \permanent\left(\begin{bmatrix}
            \labelbrack{\btree{[x]}}{\btree{[gth]}} & \labelbrack{\btree{[z[x]]}}{\btree{[gth]}} & \labelbrack{\btree{[z[x]]}}{\btree{[gth]}}\\
            \labelbrack{\btree{[x]}}{\btree{[gth]}} & \labelbrack{\btree{[z[x]]}}{\btree{[gth]}} & \labelbrack{\btree{[z[x]]}}{\btree{[gth]}}\\
            0 & \labelbrack{\btree{[z[x][x]]}}{\btree{[gth]}} & \labelbrack{\btree{[z[x][x]]}}{\btree{[gth]}}\\
        \end{bmatrix} \right) = 2! \labelbrack{\btree{[x]}\,\btree{[z[x]]}\, \btree{[z[x][x]]}}{\btree{[gth]}}
    \end{equation}
    Another example, given the tree \(\btree{[x[x[x]][x]]}\) and \(\btree{[gth[gthr]]}\) then replacement applied with the set \(\mathcal{U} = \{\btree{[gth]}, \btree{[gthr]}\}\).
    \begin{equation}
        \btree{[x[x[x]][x]]} \replaceNode{\mathcal{U}} \btree{[gth[gthr]]} = \labelbrack{\btree{[x[x][e]]}}{\btree{[gth]}}\labelbrack{\btree{[x[x]]}}{\btree{[gthr]}} + \labelbrack{\btree{[x[x[x]][e]]}}{\btree{[gth]}}\labelbrack{\btree{[x]}}{\btree{[gthr]}} + \labelbrack{\btree{[x[x][x[e]]]}}{\btree{[gth]}}\labelbrack{\btree{[x]}}{\btree{[gthr]}}
    \end{equation}
\end{example}
Just as with pruning a matrix representation is available.
\begin{definition}[Replacement Matrix Form]\label{theory:replacementMatrix}
    Given the set of nodes to replace \(\mathcal{U}\) let the vectors \(\treevecinf, \treevecinfs\) be comprised of trees \(\blockIndex{\treevecinf}{i} \in \NTTi{\mathcal{A}}\) and \(\blockIndex{\treevecinfs}{j} \in (\NTTC \setminus \NTTi{\mathcal{A}})\). Then the replacement matrix is defined as 
    \begin{equation}
        \replaceMatrixInd{\mathcal{U}} = \begin{bmatrix}
        \treevecinf\\
        \treevecinfs
        \end{bmatrix} \outergen{\replaceNode{\mathcal{U}}} \begin{bmatrix}
            \treevecinf\\
            \treevecinfs
        \end{bmatrix} = \begin{bmatrix}
            \Identity & \widetilde{\replaceMatrixInd{\mathcal{U}}}\\
            \matzero & \matzero
        \end{bmatrix}
    \end{equation}
    where \(\widetilde{\replaceMatrixInd{\mathcal{U}}}\) is the submatrix representing possible non-zero entries.
\end{definition}
With replacement now the action of replacing elementary differentials in \BSeriesText is applied.
\begin{theorem}[Replacement on \BSeriesText]\label{theory:BSeriesReplacement}
    Let \(\psi \in \mathcal{U}\) be the set of nodes to replace and \(a(\psi)(\tree)\) is the label defined pairwise on the tree nodes where \(a(\psi)(\emptyset) = 0\).
    Given some series \(\NBSeries{\Dt}(b, \y_\initstate)\) then replacement of the nodes \(\mathcal{U}\) by labels \(a\) gives a series 
    \begin{equation}
        \NBSeries{\Dt}(a \replaceNodeBSeries{\mathcal{U}} b, \y_\initstate)
    \end{equation}
    where the new labels are defined as 
    \begin{equation}
        (a \replaceNodeBSeries{\mathcal{U}} b)(\tree) = \sum_{\stree \in \NTTC} a(\tree \replaceNode{\mathcal{U}} \stree) b(\stree).
    \end{equation}
\end{theorem}
\begin{proof}
    Letting \(a(\tree) = 0,\,\forall\,\tree \in (\NTTC \setminus \NTTi{\mathcal{A}})\) then by applying \cref{theory:treereplacement} and rearranging terms gives
    \begin{equation}
        \begin{split}
            \NBSeries{\Dt}(b,\y_\initstate) &= \sum_{\tree \in \NTTC} b(\tree) \Dt^{\treeorder{\tree}}\frac{\BDiffOp(\tree)(\y_\initstate)}{\sigma(\tree)} = \sum_{\tree \in \NTTC} \sum_{\stree \in \NTTC} b(\tree) a(\tree \replaceNode{\mathcal{U}}  \stree)\Dt^{\treeorder{\stree}}\frac{\BDiffOp(\stree)(\y_\initstate)}{\sigma(\stree)}\\
            &= \sum_{\stree \in \NTTC} (a \replaceNodeBSeries{\mathcal{U}} b)(\stree) \Dt^{\treeorder{\stree}}\frac{\BDiffOp(\stree)(\y_\initstate)}{\sigma(\stree)} = \NBSeries{\Dt}(a \replaceNodeBSeries{\mathcal{U}} b,\y_\initstate).
        \end{split}
    \end{equation}
\end{proof}
\begin{example}[Replacement]
    As an example let there be two colors \(\btree{[gth]}, \btree{[gthr]}\) and let the compressed color be \(\btree{[d]}\). The associated labels for replacement are \(a\) and \(b\). Focusing on replacing with \(\mathcal{U} = \{\btree{[d]}\}\) on the tree \(\btree{[x[x][x]]}\) gives 
    \begin{equation}
        \begin{aligned}
            (a \replaceNodeBSeries{\mathcal{U}} b)(\btree{[x[x][x]]}) &= b(\btree{[x[x][x]]}) + a\left( \blockIndex{(\replaceMatrixInd{\mathcal{U}})}{\btree{[x[x][x]]}, \btree{[d[d][d]]}}\right)\,b(\btree{[d[d][d]]}) + a\left( \blockIndex{(\replaceMatrixInd{\mathcal{U}})}{\btree{[x[x][x]]}, \btree{[d[d][x]]}}\right)\,b(\btree{[d[d][x]]}) \\
            & + a\left( \blockIndex{(\replaceMatrixInd{\mathcal{U}})}{\btree{[x[x][x]]}, \btree{[d[x][x]]}}\right)\,b(\btree{[d[x][x]]}) + a\left( \blockIndex{(\replaceMatrixInd{\mathcal{U}})}{\btree{[x[x][x]]}, \btree{[x[d][d]]}}\right)\,b(\btree{[x[d][d]]}) + a\left( \blockIndex{(\replaceMatrixInd{\mathcal{U}})}{\btree{[x[x][x]]}, \btree{[d[d][x]]}}\right)\, b(\btree{[x[d][x]]}) \\
            & + a\left( \blockIndex{(\replaceMatrixInd{\mathcal{U}})}{\btree{[x[x][x]]},\, \btree{[d[d]]}}\right)\, b(\btree{[d[d]]}) + a\left( \blockIndex{(\replaceMatrixInd{\mathcal{U}})}{\btree{[x[x][x]]},\, \btree{[d[x]]}}\right)\, b(\btree{[d[x]]})\\
            & + a\left( \blockIndex{(\replaceMatrixInd{\mathcal{U}})}{\btree{[x[x][x]]}, \btree{[d]}}\right)\, b(\btree{[d]}).
        \end{aligned}
    \end{equation}
    Expanding and focusing only on a few terms gives 
    \begin{equation}
        \scalemath{1.0}{
        \begin{aligned}
            a\left( \blockIndex{(\replaceMatrixInd{\mathcal{U}})}{\btree{[x[x][x]]}, \btree{[x[d][d]]}}\right)\,b(\btree{[x[d][d]]}) &= a(\btree{[gth]})(\btree{[x]})^2 b(\btree{[x[gth][gth]]}) + 2 a(\btree{[gth]})(\btree{[x]})a(\btree{[gthr]})(\btree{[x]}) b(\btree{[x[gth][gthr]]}) + a(\btree{[gthr]})(\btree{[x]})^2 b(\btree{[x[gthr][gthr]]}),\\
            a\left( \blockIndex{(\replaceMatrixInd{\mathcal{U}})}{\btree{[x[x][x]]},\, \btree{[d[d]]}}\right)\, b(\btree{[d[d]]}) &= a(\btree{[gth]})(\btree{[x[x][e]]})a(\btree{[gth]})(\btree{[x]})b(\btree{[gth[gth]]}) + a(\btree{[gthr]})(\btree{[x[x][e]]})a(\btree{[gth]})(\btree{[x]})b(\btree{[gthr[gth]]})\\
            &+ a(\btree{[gth]})(\btree{[x[x][e]]})a(\btree{[gthr]})(\btree{[x]})b(\btree{[gth[gthr]]}) + a(\btree{[gthr]})(\btree{[x[x][e]]})a(\btree{[gthr]})(\btree{[x]})b(\btree{[gthr[gthr]]}),\\
            a\left( \blockIndex{(\replaceMatrixInd{\mathcal{U}})}{\btree{[x[x][x]]}, \btree{[d]}}\right)\, b(\btree{[d]}) &= a(\btree{[gth]})(\btree{[x[x][x]]}) b(\btree{[gth]}) + a(\btree{[gthr]})(\btree{[x[x][x]]}) b(\btree{[gthr]}).
        \end{aligned}
        }
    \end{equation}
\end{example}

\begin{remark}[Case of Leaves Replacement Only]\label{remark:leavesReplace}
    It is noted that in the case of leaves only being replaced this does \textbf{not} indicate only scaling by symmetry. A simple counter example to the claim is given by two trees with the same symmetry factor as
    \begin{equation}
        \btree{[x[z[z]][z][z]]} \replaceNode{\mathcal{U}} \btree{[x[gth][gth][z]]} = 2\labelbrack{\btree{[z[z]]} \btree{[z]}}{\btree{[gth]}}.
    \end{equation}
\end{remark}

\subsection{Tree Cutting}
There are some relations to \BSeriesText substitution \cite{Chartier_2005_substitution} from replacement such as the vector fields are changed. However, it is important to distinguish between the two different operations. Before describing their differences, substitution will be defined and thus, tree cutting as well. For this text tree cutting is redefined from literature utilizing \cref{lemma:permanentTrees} to explain similarities in operations.
\begin{definition}[Tree Cut]\label{def:treecut}
    Given the trees \(\tree = \tau[\forestTQt]\), \( \stree = \psi[\forestTQs] \in \NTTC\) then the recursive formulation of tree cutting \(\treecut\) is
    \begin{align}
        \scalemath{1.0}{\tree \treecut \stree = \tau[\forestTQt]  \treecut  \psi[\forestTQs]} & \scalemath{0.9}{ = \begin{cases} 
            \sum_{\emptyset < \vtree < \tree} \vtree\, ((\tree \prune \vtree) \treecut \forestTQs) \\
            \qquad \textbf{ if }  \tau = \psi\\
            0 \textbf{ otherwise}
        \end{cases}},\\
        \tree \treecut \psi = \tau[\forestTQt]  \treecut \psi &= \begin{cases}
            t \textbf{ if } \tau = \psi\\
            0 \textbf{ otherwise}
        \end{cases},\\
        \tau \treecut \psi[\forestTQs] &= 0,\\
        \tau \treecut \psi &= \begin{cases}
            \tau \textbf{ if } \tau = \psi\\
            0 \textbf{ otherwise}
        \end{cases},\\
        \tree \treecut \emptyset &= 0,\\
        \emptyset \treecut \stree &= 0,\\
        \emptyset \treecut \emptyset &= 1.
    \end{align}
    For forest cutting then 
    \begin{equation}
        \scalemath{1.0}{
            \begin{aligned}
                \forestTQt \treecut \forestTQs &= \begin{cases}
                    \frac{1}{\prod_{i=1}^m (k_i !) } \permanent(\forestTQt \outergen{\treecut} \forestTQs) \textbf{ if } \bar{n} = \bar{m}\\
                    0 \textbf{ otherwise}
                \end{cases},\\
                \left(\sum_{j=1}^{p} c_j \, \fforest_j \right)\treecut \forestTQs &= \sum_{j=1}^{p} c_j (\fforest_j \treecut \forestTQs),
            \end{aligned}
        }
    \end{equation}
    where \(c_j \in \outspace\) and \(\forestTQt, \forestTQs, \fforest_j  \in  \ForestSpace{\NTTC}\).
\end{definition}
\begin{theorem}[Tree Cutting Elementary Differentials]\label{theory:treecut}
    Given two trees \(\tree,\stree \in \NTTC\) and for all colors \(\psi \in \Tcolors\) the \BSeriesText elementary differentials \(\BDiffOp(\psi)(\y) = \NBSeries{\Dt}(a, \y)\) then 
    \begin{equation}
        \frac{\BDiffOp(\stree)(\y_\initstate)}{\sigma(\stree)} = \sum_{\tree \in \NTTC} a(\tree \treecut \stree) \frac{\BDiffOp(\tree)(\y_\initstate)}{\sigma(\tree)}.
    \end{equation}
    The definition of each colors elementary differential after the cut is redefined to the classical partitioned elements \(\BDiffOp(\psi)(\y_\initstate) = \fun\comp{\psi}(\y_\initstate)\). Thus, \(\tree \treecut \stree\) corresponds to the leading label of \(\BDiffOp(\tree)(\y_\initstate)\) obtained from \(\BDiffOp(\stree)(\y_\initstate)\) by expanding the nodes by the series \(\NBSeries{\Dt}(a, \y)\).
\end{theorem}
\begin{proof}
    Given that tree cutting represents a replacement of vector fields by a \BSeriesText then using the replacement operation (\cref{theory:treereplacement}) can be substituted for this. Given that replacement reduces the space of nodes we consider two sets of colors \(\vnode_{i} \in \mathcal{V}\) and \(\unode_{i} \in \mathcal{U}\) such that every color has an \(i\)th index correspondence. Let the label detection parity function be 
    \begin{equation}
        \kappa(\labelbrack{\tau[\fforest]}{\psi}) = \begin{cases}
            \tau[\fforest] & \quad \tau = \vnode_{i} \textit{ and } \psi = \unode_{i}\\
            0& \quad \textit{o.w.}
        \end{cases}.
    \end{equation}
    Then showing how trees of \(\NTTi{\mathcal{U}}\) transform via replacement to \(\NTTi{\mathcal{V}}\) demonstrates the change in vector field. Thus, tree cutting is replicated by 
    \begin{equation}
        \mathcal{E}(\kappa(\tree \replaceNode{\mathcal{U}} \stree) \quad \forall\, \tree \in \NTTi{\mathcal{V}}\, \stree \in \NTTi{\mathcal{U}})
    \end{equation}
    where the operation \(\mathcal{E}\) removes the stumps from the trees as tree cuts are not concerned with differential directions.
    To match the definition \cref{def:treecut} then converting the parity sets down to one set by removing the tree brackets, reverting the marked pruning operation to normal pruning, and collapse the colors to replace case into the nodes equal case gives the result.
\end{proof}
As in this case it is recommended to use the linear form of \cref{theorem:generalLinearTreeRecursion} as the sum over lower cardinality prunings increases the cost. We omit such formula for simplicity.
\begin{definition}[Substitution Matrix Form]\label{theory:subsitutionMatrix}
    Let  \(\treevecinf\) be an infinite vector where each entry is a tree from \(\NTTC\). The cutting matrix $\cutMatrix$ is defined as
    \begin{equation}
        \cutMatrix =  \treevecinf \outergen{\treecut} \treevecinf.
    \end{equation}
\end{definition}

\subsection{\BSeriesText Substitution}
\begin{definition}[\BSeriesText Substitution \cite{Chartier_2010_Bseries, Hairer_2006_geometric-book}]\label{def:BSeriesSubsitution}
    Given the \BSeriesText with elementary derivative information 
    \begin{itemize}
        \item \(\gun\comp{\treeRoot}(\y) = \NBSeries{\Dt}(a,\Tilde{\y})\) where \(\BDiffOp(\tree)(\y) = \fun\comp{\treeRoot(\tree)}(\y) \quad \forall\, \tree \in \NTTC\),
        \item \(\NBSeries{\Dt}(b,\y_\initstate)\) where \(\BDiffOp(\tree)(\y) = \gun\comp{\treeRoot(\tree)}(\y) \quad \forall\, \tree \in \NTTC\).
    \end{itemize}
    Then using substitution to reduce the \BSeriesText to only elementary differentials of \(\fun\comp{\treeRoot}\) gives
    \begin{equation}
        \NBSeries{\Dt}(a \labelSubstitution b, \y_\initstate) \text{, changing }\, \BDiffOp(\tree)(\y) = \gun\comp{\treeRoot(\tree)}(\y) \rightarrow \fun\comp{\treeRoot(\tree)}(\y).
    \end{equation}
    The labels $(a \labelSubstitution b)$ are defined as
    \begin{equation}
        (a \labelSubstitution b)(\tree) = \sum_{\stree \in \NTTC} a(\tree \treecut \stree) b(\stree)
    \end{equation}
    where \(a(\emptyset_{c}) = 0\) for all colors.
\end{definition}
\begin{proof}
    Applying \cref{theory:treecut} and rearranging terms gives
    \begin{equation}
        \begin{split}
            \NBSeries{\Dt}(b,\y_\initstate) &= \sum_{\tree \in \NTTC} b(\tree) \Dt^{\treeorder{\tree}}\frac{\BDiffOp(\tree)(\y_\initstate)}{\sigma(\tree)} = \sum_{\tree \in \NTTC} \sum_{\stree \in \NTTC} b(\tree) a(\tree \treecut \stree)\Dt^{\treeorder{\stree}}\frac{\BDiffOp(\stree)(\y_\initstate)}{\sigma(\stree)}\\
            &= \sum_{\stree \in \NTTC} (a \labelSubstitution b)(\stree) \Dt^{\treeorder{\stree}}\frac{\BDiffOp(\stree)(\y_\initstate)}{\sigma(\stree)} = \NBSeries{\Dt}(a \labelSubstitution b,\y_\initstate).
        \end{split}
    \end{equation}
\end{proof}
An alternative proof of this result are found in \cite{Chartier_2005_substitution} utilizing tree cutting partitions.

\subsection{Differences Between Replacement and Substitution}
When looking at both \cref{theory:treecut} and \cref{theory:treereplacement} they look fairly similar with some main differences. One difference is the node checking in that tree cutting requires equality of root nodes to return non-zero results \(\tau = \psi\). Whereas replacement has the same condition plus the additional recursion check if the nodes are in the set of nodes to replace \(\psi \in \mathcal{U}\). This gives the action only on a subset of nodes whereas tree cutting must always act on the full set of nodes.

One of the key differences is that substitution makes the assumption that the elementary function being replaced is expanded as
\begin{equation}
    \begin{aligned}
        \gun(\y) &= b(\btree{[x]}) \Dt\,\fun(\y) + b(\btree{[x[x]]}) \Dt^2\,\fun_{\y}(\y)\fun(\y)\\
            &+ b(\btree{[x[x][x]]}) \frac{\Dt^3}{2}\,\fun_{\y,\y}(\y)(\fun(\y), \fun(\y))  + b(\btree{[x[x[x]]]}) \Dt^3\,\fun_{\y}(\y)\fun_{\y}(\y)\fun(\y) + \dots,
    \end{aligned}
\end{equation}
where there is a state dependent input \(\y\) for each elementary differential. On the other hand, replacement makes zero assumptions as to what elementary differential is being replaced. This gives a more general expansion
\begin{equation}
    \begin{aligned}
        \gun(\y) &= \left[ a(\btree{[x]}) \Dt\,\fun + a(\btree{[x[x]]}) \Dt^2\,\fun_{\y}\fun + a(\btree{[x[x][x]]}) \frac{\Dt^3}{2}\,\fun_{\y,\y}(\fun, \fun)  + b(\btree{[x[x[x]]]}) \Dt^3\,\fun_{\y}\fun_{\y}\fun  \right]\\
        &+ \left[ a(\btree{[x[e]]}) \Dt \,\fun_{\y} + a(\btree{[x[x][e]]}) \Dt^2 \, \fun_{\y,\y}(\fun, \cdot) + a(\btree{[x[x[e]]]}) \Dt^2 \, \fun_{\y}\fun_{\y} + \dots \right](\y - \y_\initstate)\\
        &+ \dots.
    \end{aligned}
\end{equation}

For example the function could be a constant one such as  
\begin{equation}
    \begin{aligned}
        \gun(\y) &= \frac{1}{2}\fun(\y_\initstate) + \frac{1}{2}\fun(\y_{\initstate + 1/2}) = \fun(\y_\initstate) + \frac{1}{8} \Dt\,\fun_{\y}(\y_\initstate)\fun(\y_\initstate) + \dots.
    \end{aligned}
\end{equation}
\BSeriesText substitution would need higher state derivatives with non-constant inputs, but those must always be zero since the operation is constant. For dynamics which require state dependency replacement handles this via stumps. An example is
\begin{equation}
    \gun(\y) = \fun_{\y}(\y_\initstate)(\y - \y_\initstate) + \fun(\y_\initstate)
\end{equation} 
where for \cref{theory:treereplacement} the label is \(a(\psi)(\btree{[x]}) = 1\), \(a(\psi)(\btree{[x[e]]}) = 1\), and all other terms are zero.

Even though replacement doesn't require linear spans of elementary differentials the operations evaluated at initial conditions look \textit{similar}. We follow the derivation from \cite{Chartier_2005_substitution} to expand substitution and replacement. For simplicity of presentation for replacement we assume \(a(\labelbrack{\tree}{\psi}) = a(\tree)\). From substitution the Taylor expansion requires
\begin{equation}
    \frac{\ddiff}{\ddiff t}\y = \gun(\y) = \NBSeries{\Dt}\left(b, \y\right) = \gun + \gun_{\y}(\y - \y_\initstate) + \frac{1}{2!}\,\gun_{\y,\y}(\y - \y_\initstate)^2 + \dots
\end{equation}
whereas replacement relaxes the condition to give 
\begin{equation}
    \frac{\ddiff}{\ddiff t}\y = \mathbf{k} + \mathbf{A}(\y - \y_\initstate) + \frac{1}{2!}\mathcal{T}(\y - \y_\initstate)^2 + \dots
\end{equation}
where \(\mathbf{k}, \mathbf{A}, \mathcal{T}\) are constant series evaluated at the initial condition with number of stumps being \(0,1,2\) respectively and individually have arbitrary labels. For example focusing on a few terms of the expansion are
\begin{equation}
    \gun_{\y}(\y_\initstate)\, \gun(\y_\initstate) \quad \text{vs.} \quad \mathbf{A}\, \mathbf{k}
\end{equation}
then applying substitution gives
\begin{equation}
    \begin{aligned}
        \gun(\y_\initstate) &= b(\btree{[x]}) \Dt\,\fun + b(\btree{[x[x]]}) \Dt^2\,\fun_{\y}\,\fun + \dots,\\
        \gun_{\y}(\y_\initstate) & = \frac{\ddiff}{\ddiff\y}(b(\btree{[x]}) \Dt\, \fun(\y) + b(\btree{[x[x]]}) \Dt^2\, \fun_{\y}(\y)\,\fun(\y) + \dots)(\y_\initstate)\\
        &= b(\btree{[x]}) \Dt \,\fun_{\y} + b(\btree{[x[x]]}) \Dt^2 \,\fun_{\y,\y}(\fun(\y), \cdot) + b(\btree{[x[x]]}) \Dt^2\, \fun_{\y}\,\fun_{\y}  + \dots,\\
        \gun_{\y}(\y_\initstate)\, \gun(\y_\initstate) &= b(\btree{[x]})^2 \Dt^2 \,\fun_{\y}\,\fun +  2b\left(\btree{[x]}\,\btree{[x[x]]}\right)\frac{\Dt^3}{2} \,\fun_{\y,\y}\left(\fun, \fun\right) + 2b\left(\btree{[x]}\,\btree{[x[x]]}\right)\Dt^3 \,\fun_{\y}\,\fun_{\y}\,\fun + \dots,
    \end{aligned}
\end{equation}
whereas in contrast with replacement the terms are 
\begin{equation}
    \begin{aligned}
        \mathbf{k} &= a(\btree{[x]}) \Dt \,\fun + a(\btree{[x[x]]}) \Dt^2 \,\fun_{\y}\,\fun + \dots,\\
        \mathbf{A} & = a(\btree{[x[e]]}) \Dt \,\fun_{\y} + a(\btree{[x[x][e]]}) \Dt^2 \, \fun_{\y,\y}(\fun, \cdot) + a(\btree{[x[x[e]]]}) \Dt^2 \, \fun_{\y}\,\fun_{\y} + \dots ,\\
        \mathbf{A} \mathbf{k} & = a\left(\btree{[x[e]]}\btree{[x]}\right) \Dt^2 \,\fun_{\y}\,\fun + 2a\left(\btree{[x[x][e]]}\btree{[x]}\right) \frac{\Dt^3}{2} \, \fun_{\y,\y}(\fun, \fun) + a\left(\btree{[x[x[e]]]}\btree{[x]} + \btree{[x[e]]}\btree{[x[x]]}\right) \Dt^3 \, \fun_{\y}\,\fun_{\y}\,\fun + \dots ,\\
    \end{aligned}
\end{equation}
thus allowing for more general representations.

If the focus is only on vector fields where only \(\mathbf{k}\) has non-zero labels then, under some reasonable label transformations and zeroing, substitution gives the same results as replacement. Such labels become dependent on the structure of the problem. Though for such a simple case then alternatively composition could be used to replicate the same results as well.

Even so, substitution still treats it as a state dependent function whereas replacement treats it as a constant one. Additionally, substitution will generate additional terms which become zeroed by the labels afterwords. Thus, replacement meets the expected formulation of the results.

Another key difference is that replacement does not change the elementary differentials representations to a new field whereas substitution does as
\begin{equation}
    \begin{aligned}
        \BDiffOp(\btree{[x]})(\y_\initstate) = \gun(\y_\initstate) \quad &\underset{Substitution}{\rightarrow} \quad \BDiffOp(\btree{[x]})(\y_\initstate) = \fun(\y_\initstate),\\
        Replacement: \quad &\BDiffOp(\btree{[gth]})(\y_\initstate) = \gun(\y_\initstate) \quad \BDiffOp(\btree{[x]})(\y_\initstate) = \fun(\y_\initstate).
    \end{aligned}
\end{equation}
This allows for mixtures of terms for replacement such as \(\btree{[x[gth]]}\) whereas substitution fails to encapsulate such terms. Such terms can arise from method structure and/or composition.
Overall, substitution shares concepts with replacement, but there are distinct differences between them.

\section{\BSeriesText Time Labels\label{section:TFDEQ}}

We extend the expansion \eqref{eq:modifiedAsymptoticConstant} to include scalar time dependent functions on the base elementary differentials. This gives \BSeriesText labels with time parameterized labels with an underlying representation. We now provide the full definition of the approximate problem \cref{eq:ApproximationProblem} as
\begin{equation}\label{eq:modifiedAsymptotic}
    \begin{split}
        \gun(t,\vg(t)) &= \sum_{p \in (\mathbb{Z}^{+})^{\nparts}} \frac{\gun_{p}(t)}{\prod_{j=1}^{\nparts} p_{j}!} \,\prod_{j=1}^{\nparts} \,(\vg\comp{j}(t) - \vg\comp{j}(\tau_\initstate))^{p_{j}},\\
        \dforest_{p} &= \prod_{j=1}^{\nparts} \xnode_{c_{p_{j}}}^{p_{j}},\\
        \gun_{p}(t) &= \sum^{r}_{q=1} \phi_{p,q}(t)\,\basisK_{p,q} = \NBSeries{\Dt}^{\dforest_{p}}(\timeLabel{a^{[p]}}{\tau_\initstate}{t}, \vg(\tau_\initstate)),\\
        \frac{\basisK_{p,q}}{\prod_{j=1}^{\nparts} p_{j}!} &= \NBSeries{\Dt}^{\dforest_{p}}(a^{[p,q]}, \vg(\tau_\initstate)) = \frac{\BDiffOp(\psi[\dforest_{p}])(\vg(\tau_\initstate))}{\sigma(\psi[\dforest_{p}])},
    \end{split}
\end{equation}
where \(\phi_{p,q}(t): \mathbb{C} \rightarrow \mathbb{C}\) are smooth functions, \(r\) is either some positive integer or infinite, \(\psi\) is some node that is to be replaced, and \(a^{[p,q]}\) is some label which is zero for any trees which contain node \(\psi\). The new type of label \(\timeLabel{a^{[p]}}{\tau_\initstate}{t}\) is a new time parameterized form of \BSeriesText. 

In regard to the expansion \eqref{eq:modifiedAsymptotic}, we do not expand the function \(\phi_{p,q}(t)\) out in series. Instead, we embed the time dependency in the \BSeriesText labels themselves. This allows for a more general expansion of numerical methods which are continuous. By doing such we gain alternative views of the expansion and its associated properties. It is noted that utilizing time parameterized (and general parameterizations) labels are not new \cite{butcherAlgebraicTheoryIntegration,brouder2000runge,miyatakeCharacterizationEnergyPreservingMethods2016,hairer2010energy}.

We provide formal definitions of such time dependent series and how they relate to the underlying system expansion. Their expansions are simplified with a compact reduction such that the sub-sequences are reduced to single nodes as time dependent functions. We introduce the representation for time parameterization with its underlying classical series that is approximated by some method. An example applied to Runge-Kutta methods is provided.

\subsection{Tree Time Nodes}
To motivate the derivation we discuss a simplified differential equation and generalize back to \eqref{eq:modifiedAsymptotic}. With the simplified problem
\begin{equation}\label{eq:TFProblem}
    \begin{split}
        \ddifffrac{t} \vg &= \fun(t,\,\vg,\,\wg),\\
        \ddifffrac{t} \wg &= \gun(t,\wg,\vg) = \sum_{p_1,p_2 = 0} \phi\comp{p_1,p_2}(t)\, \basisK\comp{p_1,p_2} (\vg^{p_1}\,\wg^{p_2}),\\
        t &\in [t_\initstate,\, t_\finalstate], \quad \vg(t_\initstate) = \vzero, \quad \wg(t_\initstate) = \vzero
    \end{split}
\end{equation}
the focus is only on \(2\) partitions. We use this as an example to motivate the underlying representation of the time dependent terms. Extensions come naturally from multiple colors applied to additional terms in the  sums and partitions.

To describe \eqref{eq:TFProblem} classically there are \(3\) colors one for the state \(\vg\), time \(t\), and the approximate state \(\wg\). Focusing on the case when no stump nodes are present and there are only two colors \(\BDiffOp(\btree{[x]})(\y) = \fun(\y)\), \(\BDiffOp(\btree{[b]})(\y) = \phi\comp{0,0} \basisK\) and \(\phi\comp{0,0}\) is the leading time function. 
These trees can be visualized with the list of orders $1$ to $3$ with only one color of each node type. For simplicity, we write \(\phi = \phi\comp{0,0}\). The non-zero elementary differential representations are as follows:
\begin{equation}
    \scalemath{0.9}{\begin{array}{c|c|c}
        \text{Order} & \text{Trees} & \text{Trees Elementary Differentials} \\
        \hline
        1 & \btree{[x]} \quad \btree{[b]} & \fun \quad \phi \, \basisK  \\
        \hline
        2 & \btree{[x[x]]} \quad \btree{[x[t]]} \quad \btree{[x[b]]} \quad \btree{[b[t]]} & \fun_v\, \fun \quad \fun_t \quad \fun_w\, \phi\, \basisK \quad \phi_t\, \basisK \\
        \hline
        3 & \btree{[x[x][x]]} \quad \btree{[x[t][x]]} \quad \btree{[x[t][t]]} \quad \btree{[x[b][x]]} \quad \btree{[x[b][b]]} & \fun_{v,v}(\fun, \fun) \quad \fun_{v,t}\,\fun \quad \fun_{t,t} \quad \fun_{v,w}(\fun, \phi\, \basisK) \quad \fun_{w,w}(\phi\, \basisK, \phi\, \basisK)\\
         & \quad \btree{[x[b][t]]} \quad \btree{[b[t][t]]} \quad \btree{[x[x[x]]]} \quad \btree{[x[x[t]]]} \quad \btree{[x[x[b]]]} \quad \btree{[x[b[t]]]} & \quad \fun_{w,t}(\phi\, \basisK) \quad \phi_{t,t} \basisK \quad \fun_v\,\fun_v\,\fun \quad \fun_v\,\fun_t \quad \fun_v\,\fun_w \,\phi\, \basisK \quad \fun_v\, \phi_t\, \basisK \\
        \hline
    \end{array}.}
\end{equation}
Next we compactify the representation along the time direction of the node \(\btree{[b]}\). As an example taking the tree \(\btree{[x[b]]}\) and looking at its elementary differential expansion gives
\begin{equation}
    a_1 \Dt^2\, \fun_w \,\phi\, \basisK + a_2 \Dt^3\, \fun_w\, \phi_t\, \basisK + \frac{a_3}{2!} \Dt^4\,\fun_w\, \phi_{t,t}\, \basisK + \dots = (a_1 \phi + a_2 \Dt \phi_t + \frac{a_3}{2!} \Dt \phi_{t,t} + \dots) \Dt^2\, \fun_w \,\basisK = a(t)\, \fun_w \,\basisK.
\end{equation}
Thus, we have pushed the derivatives of the function \(\phi\) into the new reduced label as a scalar function \(a(t)\).

To display the compression we let \(\btree{[gth]}\) represent the time derivatives summed on \(\phi\) thus, replacing the color \(\btree{[b]}\) with a compressed one. This can be visualized by trees orders $1$ to $3$ for one color of each node. Again, focusing on the case when no stump nodes are present the colors are \(\BDiffOp(\btree{[x]})(\y) = \fun(\y)\) and the reduced node \(\BDiffOp(\btree{[gth]}) = \basisK\). For the non-zero elementary differentials this gives:
\begin{equation}
    \begin{array}{c|c|c}
        \text{Order} & \text{Trees} & \text{Trees Elementary Differentials}\\
        \hline
        1 & \btree{[x]} \quad \btree{[gth]} & \fun \quad \basisK \\
        \hline
        2 & \btree{[x[x]]} \quad \btree{[x[t]]} \quad \btree{[x[gth]]} & \fun_v\,\fun \quad \fun_t \quad \fun_w\, \basisK \\
        \hline
        3 & \btree{[x[x][x]]} \quad \btree{[x[t][x]]} \quad \btree{[x[t][t]]} \quad \btree{[x[gth][x]]} \quad \btree{[x[gth][gth]]} & \fun_{v,v}(\fun, \fun) \quad \fun_{v,t}\,\fun \quad \fun_{t,t} \quad \fun_{v,w}(\fun, \basisK)\\
         & \btree{[x[gth][t]]} \quad \btree{[x[x[x]]]} \quad \btree{[x[x[t]]]} \quad \btree{[x[x[gth]]]} & \fun_{w,t}\,\basisK \quad \fun_v\,\fun_v\,\fun \quad \fun_v\,\fun_t \quad \fun_v\,\fun_w\, \basisK \\
        \hline
    \end{array}
\end{equation}
Each new node is assigned a smooth scalar function which becomes a part of the labels including all of its derivatives. Thus, returning to the general case \eqref{eq:modifiedAsymptotic} where every leading time node has a pairing with an underlying classical node representation. Thus, every label becomes dependent on the scalars being pushed down. To get the new label a time elementary operator is defined.
\begin{definition}[Nodes Elementary Function]\label{def:reductionNodeFunctions}
    Let there be two sets of colored nodes \(\unode_{k} \in \mathcal{U}\) and \(\vnode_{k} \in \mathcal{V}\) where the index \(k\) corresponds to the comparable colors between nodes. Let \(\psi \in \mathcal{A} = \Tcolors \setminus (\mathcal{U} \cup \mathcal{V})\) denote the set of remaining colors utilized with non-reduced nodes.
    The nodes elementary time operator \(\reductionFun : \NTTC \rightarrow (\outspace \rightarrow \outspace)\) is defined as
    \begin{equation}
        \begin{split}
            \reductionFun(\vnode_{k},\dforest)(x) &= \phi_{\vnode_{k},\dforest}(x),\\
            \reductionFun(\vnode_{k},\dforest)^{(i)}(x) &= \phi^{(i)}_{\vnode_{k},\dforest}(x),\\
            \reductionFun(\unode_{k},\dforest)(x) = 0,\\
            \reductionFun(\psi, \dforest)(x) = 1,
        \end{split}
    \end{equation}
    where \(x \in \outspace\) and \(\dforest = \xnode_{1}^{b_1}\,\xnode_{2}^{b_2}\, \dots\, \xnode_{p}^{b_{p}}\) represents the colored stumped nodes for colors \((\mathcal{A} \cup \mathcal{V}) \setminus \{\btree{[t]}\}\) plus the empty set. Any node not defined is \(0\). A shorthand notation is used for when all stumps of the same node contain the same function as
    \begin{equation}
        \reductionFun(\unode_{k})(x) = \phi_{\vnode_{k}} = \reductionFun(\vnode_{k},\dforest) \quad \forall \, \dforest.
    \end{equation}
    For simplicity throughout let
    \begin{equation}
        \timeElemOp{x}^{\reductionFun}(\tree) = \begin{cases}
            \reductionFun(\vnode_{k}, \dforest)^{(j)}(x)  \textbf{ if } \tree = \labelbrack{ (\unode_{k} \ButcherProduct \btree{[t]}^{j}) \ButcherProduct \dforest}{\vnode_{k}},\\
            1 \textbf{ if } \tree = 1,\\
            0 \textbf{ otherwise},
        \end{cases}
    \end{equation}
    where \(j \ge 0\).
\end{definition}
With the elementary operation for the time dependent terms defined then the reduction nodes have a defined relation to the unreduced nodes.
\begin{theorem}[Reduction Node Representation]\label{theorem:reductionForceRepresentation}
    Utilizing the set of colors and functions from \cref{def:reductionNodeFunctions}. For a tree \(\stree \in \NTTC\) and a time label \(\timeLabel{a}{t_\initstate}{t_\finalstate}^{\reductionFun}\) over the time domain \([t_\initstate, t_\finalstate]\) the underlying label is \(a\), with \(a(\tree) = 0 \quad \forall\, \tree \notin \NTTi{\Tcolors \setminus \mathcal{V}}\), then the relationship to the time dependent label is
    \begin{equation}
        \timeLabel{a}{t_\initstate}{t_\finalstate}^{\reductionFun}(\stree) = \sigma(\stree) \sum_{\rtree \in \NTTC}  \Dt^{\treeorder{\rtree} - \treeorder{\stree}} \timeElemOp{t_\initstate}^{\reductionFun}(\rtree \replaceNode{\mathcal{U}} \stree) \frac{a(\rtree)}{\sigma(\rtree)}.
    \end{equation}
\end{theorem}
\begin{proof}
    For a tree \(\stree \in \NTTi{\Tcolors \setminus \mathcal{V}}\) the elementary operator \(\timeElemOp{x}^{\reductionFun}\) is only one when given one. Then over the sum of all trees must give the exact same tree. Thus, giving the result of \(\timeLabel{a}{t_\initstate}{t_\finalstate}^{\reductionFun}(\stree) = a(\stree)\) which is the expected result for trees with time nodes.

    For the case when there are time nodes the operator \(\timeElemOp{t_\initstate}^{\reductionFun}(\rtree \replaceNode{\mathcal{U}} \stree)\) will only be non-zero when the reduced nodes overlap, and they have time nodes.  Thus, giving products of time elementary differentials with no additional factors. Lastly to add the effects of factorials in expansion the symmetries of both trees are divided. Lack of symmetry in the tree gives additions of alternative elementary time derivative selections and additional symmetry in the time derivatives is captured.
\end{proof}
This gives a reduced \BSeriesText representation for any elementary differential equation of the form \(a(t)\fun(\y)\) where \(a(t)\) is a scalar function. Looking back at \eqref{eq:TFProblem} the elementary nodes become
\begin{equation}
    \begin{split}
        \reductionFun(\btree{[x]})(x) &= 1\\
        \reductionFun(\btree{[gth]}, \xnode_1^{p_1}\,\xnode_2^{p_1})(x) &= \phi\comp{p_1,p_2}(x)\\
    \end{split}
\end{equation}
Multiple examples of using the reduction nodes and their underlying representations are provided in \cref{table:reductionExamples}.
However, there are cases where it can become difficult to resolve what nodal definitions are to be used. A definition is provided such that the nodal locations are not necessary.

\begin{definition}[Reudction Node Non-Nodal]\label{def:reductionForceRepresentationNonNodal}
    For a tree \(\tree \in \NTTC\) and the time label \(\timeLabel{a}{t_\initstate}{t_\finalstate}\) over the time domain \([t_\initstate, t_\finalstate]\) and \(a\) is the underlying label. Let the set 
    \begin{equation}
        \mathcal{Q}_{\Tcolors}(\tree) = \{\stree\, |\, \stree \in \NTTi{\Tcolors},\, \stree \prune \tree = c\, \btree{[t]}^{q},\, c \ne 0,\, q > 0\}
    \end{equation}
    where the set produces all trees which only have prunings of time nodes thus, collecting possible time function locations.
    The relationship between the parameterized label and \(a\) is 
    \begin{equation}
        \timeLabel{a}{t_\initstate}{t_\finalstate}(\tree) = \sigma(\tree) \sum_{\stree \in \mathcal{Q}_{\Tcolors}(\tree)}  \Dt^{\treeorder{\stree} - \treeorder{\tree}} \timeElemOp{t_\initstate}(\tree, \stree) \frac{a(\stree)}{\sigma(\stree)}
    \end{equation}
    where \(\timeElemOp{t_\initstate}(\tree, \stree)\) returns a scalar function evaluated at \(t_{\initstate}\) based off the two trees.
\end{definition}
This leaves the time dependent functions arbitrary per-tree instead of nodally defined. All series produced by nodal operations can be represented by this. Such representation will be used if the locations of the time nodes are difficult to resolve with nodal based definitions. An example is when replacement \cref{theory:BSeriesReplacement} is applied. Replacement can give two trees of different orders with the same time dependent function only with a slightly different constant scaling. Note this does not mean \textit{any} choice of locations for time nodes is correct; just the structure of the tree and additional context are needed to put them in correct locations.
\begin{remark}[Stumps in Reduction Nodes]
    Note that in \cref{def:reductionNodeFunctions} stumps sharing the same colors as not stumps or higher ones are not required. One could simply split the elementary differential into separate colors to represent the same effect and have multiple zero elementary differentials. However, so many zeros can make the expansion have too many wasteful terms.
\end{remark}
\begin{remark}[Composition and Reduction Nodes]
    Note that there is nothing preventing the usage of composition (\cref{theory:BSeriesComposition}) with the reduced nodes. However, additional care is required as the time dependent effects must be composed upon as well. Applying repeated composition on reduced labels cannot use \cref{theory:repeatedIntegration} \textit{directly} since time dependent prune matrices will be different at each step.
\end{remark}
\begin{sidewaystable}
    \[
    \scalemath{0.85}{
    \begin{array}{l|l|l}
        \text{Reduction Tree Label} & \text{Represented Sequence} & \text{Infinite Forcing Expansion}\\[5pt]
        \hline
        \Dt \, \timeLabel{a}{t_\initstate}{t_\finalstate}^{\reductionFun}\left(\btreega{[gth]}\right) \BDiffOp\left(\btreega{[gth]}\right)(\y) & \Dt \, \Big( a(\btreega{[b]})\phi(t_\initstate) + \Dt \, a\left(\btreega{[b[t]]}\right)\phi^{(1)}(t_\initstate) + \Dt^2 \, \cfrac{a\left(\btreega{[b[t][t]]}\right)}{2!}\phi^{(2)}(t_\initstate)  & \sum\limits_{i=0}^{\infty} \Dt^{i+1} \,\cfrac{a\left(\btreega{[b]} \ButcherProduct \btreega{[t]}^{i}\right)}{i!} \BDiffOp\left(\btreega{[b]} \ButcherProduct \btreega{[t]}^{i}\right)(\y)\\[5pt]
        & \quad + \Dt^3 \,\cfrac{a\left(\btreega{[b[t][t][t]]}\right)}{3!} \phi^{(3)}(t_0)  + \dots \Big)\BDiffOp\left(\btreega{[gth]}\right)(\y) & \\[5pt]
        \hline
        \Dt^2 \, \timeLabel{a}{t_\initstate}{t_\finalstate}^{\reductionFun}\left(\btreega{[x[gth]]}\right) \BDiffOp\left(\btreega{[x[gth]]}\right)(\y) & \Dt^2 \, \Big( a\left(\btreega{[x[b]]}\right)\phi(t_\initstate) + \Dt \, a\left(\btreega{[x[b[t]]]}\right)\phi^{(1)}(t_\initstate) + \Dt^2 \, \cfrac{a\left(\btreega{[x[b[t][t]]]}\right)}{2!}\phi^{(2)}(t_\initstate) & \sum\limits_{i=0}^{\infty}  \Dt^{i+2} \, \cfrac{a\left(\btreega{[x]} \ButcherProduct \left(\btreega{[b]} \ButcherProduct \btreega{[t]}^{i}\right)\right)}{i!}  \\[5pt]
        & \quad + \Dt^3 \,\cfrac{a\left(\btreega{[x[b[t][t][t]]]}\right)}{3!} \phi^{(3)}(t_\initstate)  + \dots \Big)\BDiffOp\left(\btreega{[x[gth]]}\right)(\y) &  \BDiffOp\left(\btreega{[x]} \ButcherProduct \left(\btreega{[b]} \ButcherProduct \btreega{[t]}^{i}\right)  \right)(\y)\\[5pt]
        \hline
        \Dt^3 \, \timeLabel{a}{t_\initstate}{t_\finalstate}^{\reductionFun}\left(\btreega{[x[gth][x]]}\right) \BDiffOp\left(\btreega{[x[gth][x]]}\right)(\y) & \Dt^3 \, \Big( a\left(\btreega{[x[b][x]]}\right)\phi(t_\initstate) + \Dt \, a\left(\btreega{[x[b[t]][x]]}\right)\phi^{(1)}(t_\initstate) + \Dt^2 \, \cfrac{a\left(\btreega{[x[b[t][t]][x]]}\right)}{2!} \phi^{(2)}(t_\initstate) & \sum\limits_{i=0}^{\infty}  \Dt^{i+3} \, \cfrac{a\left( \btreega{[x[x]]} \ButcherProduct \left(\btreega{[b]} \ButcherProduct \btreega{[t]}^{i}\right)\right)}{i!} \\[5pt]
        & \quad + \Dt^3 \, \cfrac{a\left(\btreega{[x[b[t][t][t]][x]]}\right)}{3!} \phi^{(3)}(t_\initstate)  + \dots \Big)\BDiffOp\left(\btreega{[x[gth][x]]}\right)(\y) & \BDiffOp\left( \btreega{[x[x]]} \ButcherProduct \left(\btreega{[b]} \ButcherProduct \btreega{[t]}^{i}\right)\right)(\y) \\[5pt]
        \hline
        \Dt^3 \, \timeLabel{a}{t_\initstate}{t_\finalstate}^{\reductionFun}\left(\btreega{[x[gth][t]]}\right) \BDiffOp\left(\btreega{[x[gth][t]]}\right)(\y) & \Dt^3 \, \Big( a\left(\btreega{[x[b][t]]}\right)\phi(t_\initstate) + \Dt \, a\left(\btreega{[x[b[t]][t]]}\right) \phi^{(1)}(t_\initstate) + \Dt^2 \, \cfrac{a\left(\btreega{[x[b[t][t]][t]]}\right)}{2!} \phi^{(2)}(t_\initstate) & \sum\limits_{i=0}^{\infty}  \Dt^{i+3} \, \cfrac{a\left( \btreega{[x[t]]} \ButcherProduct \left(\btreega{[b]} \ButcherProduct \btreega{[t]}^{i}\right)\right)}{i!}\\[5pt]
        & \quad + \Dt^3 \, \cfrac{a\left(\btreega{[x[b[t][t][t]][t]]}\right)}{3!} \phi^{(3)}(t_\initstate)  + \dots \Big)\BDiffOp\left(\btreega{[x[gth][t]]}\right)(\y) &  \BDiffOp\left( \btreega{[x[t]]} \ButcherProduct \left(\btreega{[b]} \ButcherProduct \btreega{[t]}^{i}\right) \right)(\y)\\[5pt]
        \hline
        \Dt^3 \, \cfrac{\timeLabel{a}{t_\initstate}{t_\finalstate}^{\reductionFun}\left(\btreega{[x[gth][gth]]}\right)}{2!} \BDiffOp\left(\btreega{[x[gth][gth]]}\right)(\y) &  \Dt^3 \, \Big(\cfrac{a\left(\btreega{[x[b][b]]}\right)}{2!} \left(\phi(t_\initstate)\right)^2 + \Dt \, a\left(\btreega{[x[b[t]][b]]}\right) \phi^{(1)}(t_\initstate)\phi(t_\initstate) & \frac{1}{2!} \sum\limits_{i=0}^{\infty}\sum\limits_{j=0}^{\infty}  \Dt^{i+j+3}\\[5pt]
        & \quad + \Dt^2 \, \cfrac{a\left(\btreega{[x[b[t][t]][b]]}\right)}{2!} \phi^{(2)}(t_\initstate)\phi(t_\initstate) + \Dt^2 \, \cfrac{a\left(\btreega{[x[b[t]][b[t]]]}\right)}{2!} \left(\phi^{(1)}(t_\initstate)\right)^2 + \Dt^3 \, \cfrac{a\left(\btreega{[x[b[t][t]][b[t]]]}\right)}{2!}  \phi^{(2)}(t_\initstate)\phi^{(1)}(t_\initstate) &  \cfrac{a\left(\left(\btreega{[x]} \ButcherProduct \left(\btreega{[b]} \ButcherProduct \btreega{[t]}^{i}\right)\right) \ButcherProduct \left(\btreega{[b]}\ButcherProduct \btreega{[t]}^{j}\right) \right)}{i!j!} \\[5pt]
        & \quad + \Dt^4 \, \cfrac{a\left(\btreega{[x[b[t][t]][b[t][t]]]}\right)}{2!(2!)^2} \left(\phi^{(2)}(t_\initstate)\right)^2 + \Dt^4 \, \cfrac{a\left(\btreega{[x[b[t][t][t]][b[t]]]}\right)}{3!}  \phi^{(3)}(t_\initstate) \phi^{(1)}(t_\initstate) + \dots \Big)\BDiffOp\left(\btreega{[x[gth][gth]]}\right)(\y) & \BDiffOp\left(\left(\btreega{[x]} \ButcherProduct \left(\btreega{[b]} \ButcherProduct \btreega{[t]}^{i}\right)\right) \ButcherProduct \left(\btreega{[b]}\ButcherProduct \btreega{[t]}^{j}\right)\right)(\y)\\[5pt]
        \hline
        \Dt^3 \, \timeLabel{a}{t_\initstate}{t_\finalstate}^{\reductionFun}\left(\btreega{[x[x[gth]]]}\right) \BDiffOp\left(\btreega{[x[x[gth]]]}\right)(\y) & \Dt^3 \, \Big( a\left(\btreega{[x[x[b]]]}\right) \phi(t_\initstate) + \Dt \, a\left(\btreega{[x[x[b[t]]]]}\right) \phi^{(1)}(t_\initstate) + \Dt^2 \, \cfrac{a\left(\btreega{[x[x[b[t][t]]]]}\right)}{2!} \phi^{(2)}(t_\initstate) & \sum\limits_{i=0}^{\infty}  \Dt^{i+3} \, \cfrac{a\left(\btreega{[x]} \ButcherProduct \left(\btreega{[x]} \ButcherProduct \left(\btreega{[b]} \ButcherProduct \btreega{[t]}^{i}\right)\right) \right)}{i!} \\[5pt]
        & \quad + \Dt^3 \, \cfrac{a\left(\btreega{[x[x[b[t][t][t]]]]}\right)}{3!} \phi^{(3)}(t_\initstate)  + \dots \Big)\BDiffOp\left(\btreega{[x[x[gth]]]}\right)(\y) & \BDiffOp\left(\btreega{[x]} \ButcherProduct \left(\btreega{[x]} \ButcherProduct \left(\btreega{[b]} \ButcherProduct \btreega{[t]}^{i}\right)\right)  \right)(\y)\\[5pt]
        \hline
    \end{array}}
    \]
    \caption{Examples of the reduction node representation \cref{theorem:reductionForceRepresentation}.}
    \label{table:reductionExamples}
\end{sidewaystable}

\subsection{Applying Time Labels to Runge-Kutta Methods}
We describe an example of using the reduced labels for a Runge-Kutta method with \(2\) partitions. Specifically, we look at the simplified problem
\begin{equation}\label{eq:diffTimeSum}
    \frac{\ddiff}{\ddiff t}\y = \fun\left(\y, \z_0 + \sum_{i=1}^{m} g_i(t)\,\basisK_i\right), \quad t \in [t_\initstate, t_\finalstate],
\end{equation}
where the terms \(\BDiffOp(\lambda_i) = \basisK_i\) and \(\reductionFun(\lambda_i, \emptyset) = g_i\) are the new reduced nodes. When \(g_i(t_\initstate) \ne 0\) the constant terms are pushed into \(\z_0\) and treated via composition if preferable.

Taking the conditions to an arbitrary staged Runge-Kutta gives the new conditions.
\begin{theorem}[Runge-Kutta Time Expansion]\label{theorem:RKTime}
    Given some Runge-Kutta method scheme coefficients \((A,b,c)\) and the differential equation \cref{eq:diffTimeSum}. The standard order conditions \cite{Hairer_book_I,Butcher_2021_Bseries-book,Butcher_2016_book} are augmented to represent the time dependent nodes \(\reductionFun(\lambda_i ,\emptyset)(x) = g_i(x)\) and \(\reductionFun(\tau)(x) = 1\). The augmented order conditions are given by the functions
    \begin{equation}
        \begin{split}
            \timeLabel{\Phi}{t_\initstate}{t_\finalstate}^{\reductionFun}(\tau[\forestTQt]) &= b^T\, \left(\prod_{i=1}^{n} \timeLabel{\phi}{t_\initstate}{t_\finalstate}^{\reductionFun}(\tree_i)^{\elemtimes q_i}\right),\\
            \timeLabel{\Phi}{t_\initstate}{t_\finalstate}^{\reductionFun}(\lambda_i[\forestTQt]) &= 0,\\
            \timeLabel{\Phi}{t_\initstate}{t_\finalstate}^{\reductionFun}(\tau) &= b^T \vone,\\
            \timeLabel{\Phi}{t_\initstate}{t_\finalstate}^{\reductionFun}(\lambda) &= 0,\\
            \timeLabel{\phi}{t_\initstate}{t_\finalstate}^{\reductionFun}(\tau[\forestTQt]) &= A\, \left(\prod_{i=1}^{n} \timeLabel{\phi}{t_\initstate}{t_\finalstate}^{\reductionFun}(\tree_i)^{\elemtimes q_i}\right),\\
            \timeLabel{\phi}{t_\initstate}{t_\finalstate}^{\reductionFun}(\lambda_i[\forestTQt]) &= \vzero,\\
            \timeLabel{\phi}{t_\initstate}{t_\finalstate}^{\reductionFun}(\tau) &= c,\\
            \timeLabel{\phi}{t_\initstate}{t_\finalstate}^{\reductionFun}(\lambda_i) &= g_{i}(t_\initstate + \Dt\, c),\\
        \end{split}
    \end{equation}
    where \(\timeLabel{\Phi}{t_\initstate}{t_\finalstate}^{\reductionFun}\) is the base order condition as in literature and \(g_{i}(t_\initstate + \Dt\, c)\) is element wise evaluated over the stages of \(c\).
\end{theorem}
A proof of this is provided in \cref{proof:RKTime}. Note that we focus on the results for the simplified problem over the full problem \eqref{eq:ApproximationProblem}. Results here can be expanded to the more general case.
As an example utilizing a \(3\) stage explicit Runge-Kutta method gives 
\begin{equation}
    \begin{split}
        \YY_{1} &= \y_\initstate,\\
        \YY_{2} &= \y_\initstate + \Dt\, a_{2,1}\, \fun\left(\YY_{1}, \z_0\right),\\
        \YY_{3} &= \y_\initstate + \Dt\, a_{3,1}\, \fun\left(\YY_{1}, \z_0\right) + \Dt\, a_{3,2}\, \fun\left(\YY_{2}, \z_0 + \sum_{i=1}^{m} g_i(t_\initstate + \Dt\, a_{2,1})\,\basisK_i\right),\\
        \y_{\finalstate} &= \y_\initstate + \Dt\, b_1\, \fun\left(\YY_{1}, \z_0\right) + \Dt\, b_2\, \fun\left(\YY_{2}, \z_0 + \sum_{i=1}^{m} g_i(t_\initstate + \Dt\, c_{2})\,\basisK_i\right) \\
        &+ \Dt\, b_3\, \fun\left(\YY_{3}, \z_0 + \sum_{i=1}^{m} g_i(t_\initstate + \Dt\, c_{3})\,\basisK_i\right).
    \end{split}
\end{equation}
Some of the labels with at least one reduced node are  
\begin{equation}
    \begin{split}
        \timeLabel{\Phi}{t_\initstate}{t_\finalstate}^{\reductionFun}(\tau[\lambda_i]) = \sum_{j=1}^{3} b_j\, g_i(t_{\initstate} + \Dt\, c_j),\\
        \timeLabel{\Phi}{t_\initstate}{t_\finalstate}^{\reductionFun}(\tau[\tau[\lambda_i]]) = \sum_{j_1 = 1}^{3} \sum_{j_2 = 1}^{3} b_{j_1}\, a_{j_1, j_2}\, g_i(t_{\initstate} + \Dt\, c_{j_2}),\\
        \timeLabel{\Phi}{t_\initstate}{t_\finalstate}^{\reductionFun}(\tau[\tau \lambda_i]) = \sum_{j = 1}^{3} b_{j}\left(c_{j}\, g_i(t_{\initstate} + \Dt\, c_{j})\right),\\
        \timeLabel{\Phi}{t_\initstate}{t_\finalstate}^{\reductionFun}(\tau[\lambda_{i_1} \lambda_{i_2}]) = \sum_{j = 1}^{3} b_{j}\,\left(g_{i_1}(t_{\initstate} + \Dt\, c_{j})\, g_{i_2}(t_{\initstate} + \Dt\, c_{j})\right).\\
    \end{split}
\end{equation}

\subsection{Partial Exact Solution}
In some cases the true solution of \textit{part} of the differential equation may be known and accessible. This leaves a time dependent operator parameterized on the other states of the system. As an example given the system
\begin{equation}\label{eq:diffRed}
    \begin{split}
        \frac{\ddiff}{\ddiff t}\y\comp{1} &= \fun\comp{1}(\y\comp{1}, \y\comp{2})\\
        \frac{\ddiff}{\ddiff t}\y\comp{2} &= \fun\comp{2}(\y\comp{1}, \y\comp{2})\\
    \end{split} \quad \underbrace{\rightarrow}_{\text{ Solve } \fun\comp{2}} \quad \frac{\ddiff}{\ddiff t}\y\comp{1} = \fun\comp{1}(\y\comp{1}, \tilde{\y}\comp{2}\left(t; \y\comp{1}\right))
\end{equation}
where the solution is given by inserting Picard iteration in as 
\begin{equation}
    \tilde{\y}\comp{2}\left(t; \y\comp{1}\right) = \y_{0} + \int_{0}^{t} \fun\comp{2}(\y\comp{1}(\tau), \tilde{\y}\comp{2}\left(t; \y\comp{1}\right)) \dint \tau.
\end{equation}

Now the \BSeriesText of the simplified system has only \(1\) color in the time reduced case. Thus, some method applied to the simplified problem only needs that many colors. However, the full system has \(2\) colors. Thus, to describe the error of the original full system a transformation between them is needed.

Such simplifications allow for a transformation of the reduction node labels for when partial system is exact and the other is approximated. A label transformation is provided as follows.
\begin{theorem}[Partial Exact]\label{theorem:partialExact}
    Given the label \(\timeLabel{b}{0}{x}\) for \(x \in \outspace\) and a differential equation with \(\nparts\) partitions that are represented as colors \(\Tcolors\). Given the sub-set of colors such that \(\Tcolors_1 \cup \Tcolors_2 = \Tcolors\) and \(\Tcolors_1 \cap \Tcolors_2 = \emptyset\). Suppose the differential equation for colors \(\Tcolors_2\) is solved exactly but not the other. Let the coefficients \(\timeLabel{a}{0}{x}\comp{k}\) be for the simplified differential equation solved by some method. An assumption of \(\timeLabel{a}{0}{0}\comp{k}(\tree) = 0\) holds for all trees except the \(\emptyset\). Each integer \(k\) for the labels represents time dependent labels from higher iterations. Given the functions 
    \begin{equation}
        \begin{split}
            \varkappa\comp{k,a}(\tree) &= \varkappa\comp{k,a}(\tau[\forestTQt]) = 
            \begin{cases}
                \tau \ButcherProduct \left(\prod_{i=1}^{n}\varkappa\comp{k,a}(\tree_i)^{q_i}\right) \text{ if } \tau \in \Tcolors_1,\\
                \btree{[t]}  \, \left(\prod_{i=1}^{n}\beta\comp{k,a}(\tree_i)^{q_i}\right) \text{ if } \tau \in \Tcolors_2,
            \end{cases},\\
            \varkappa\comp{k,a}(\tau) &= \begin{cases}
                \tau \text{ if } \tau \in \Tcolors_1 \text{ or } \tau = \btree{[t]}\\
                \btree{[t]} \text{ if } \tau \in \Tcolors_2
            \end{cases},\\
            \beta\comp{k,a}(\tree) &= \beta\comp{k,a}(\tau[\forestTQt]) = 
            \begin{cases}
                \sigma(\tree) \sum_{i=1}^{\infty} \btree{[t]}^{i} \frac{\ddiff^{i}}{\ddiff x^{i}} \left( \frac{\timeLabel{a}{0}{x}\comp{k}}{\sigma}(\varkappa\comp{k+1,a}(\tree))\right)\Big|_{x=0}  \text{ if } \tau \in \Tcolors_1,\\
                \btree{[t]}  \, \left(\prod_{i=1}^{n}\beta\comp{k,a}(\tree_i)^{q_i}\right) \text{ if } \tau \in \Tcolors_2,
            \end{cases},\\
            \beta\comp{k,a}(\tau) &= \begin{cases}
                \sum_{i=1}^{\infty} \btree{[t]}^{i} \frac{\ddiff^{i}}{\ddiff x^{i}} \left(\timeLabel{a}{0}{x}\comp{k}(\tau)\right)\Big|_{x=0}   \text{ if } \tau \in \Tcolors_1 \text{ or } \tau = \btree{[t]}\\
                \btree{[t]} \text{ if } \tau \in \Tcolors_2
            \end{cases},
        \end{split}
    \end{equation}
    the transformation of the reduced labels back to the full system labels is
    \begin{equation}
        \timeLabel{b}{0}{x}(\tree) = \sigma(\tree) \frac{\timeLabel{a}{0}{x}\comp{0}}{\sigma}(\varkappa\comp{1,a}(\tree)).
    \end{equation}
\end{theorem}
A proof of this is provided in \cref{proof:partialExact}.

\begin{remark}[Simplified Case of \cref{theorem:partialExact}]\label{remark:reducMono}
    In the case that the coupling sub-system was solved exactly and has no dependency on the other state then the higher iterations are removed. This gives that the higher labels are simplified to 
    \begin{equation}
        \timeLabel{a}{0}{x}\comp{k}(\tree) = 0 \quad \forall\, k > 0, x \in \outspace, \tree \in \NTTC.
    \end{equation}
\end{remark}

\section{Tree Iterated Integrals (TII) for Exact Solutions\label{section:TII}}

In this section we provide a systematic view of the exact solution via the \textit{Tree Iterated Integrals (TII)} \cite{kreimer1999chen, butcher1987numerical}. This allows for a simple representation of the exact solution of time parameterized series. We use this formulation to describe the exact solution of \cref{eq:modifiedAsymptotic} to avoid relying on the derivatives of the labels. Additionally, we analyze the case of \textit{step size invariance} given as
\begin{equation}
    \phi_{p,q}(t) = \tilde{\phi}_{p,q}\left(\frac{t - t_\initstate}{\Dt}\right)
\end{equation}
where the step size \(\Dt = t_\finalstate - t_\initstate\). Such functions are present in some continuous numerical methods, such as MRI and CRK. This is so the method doesn't change if a different starting time point \(t_\initstate\) or step length \(\Dt\) is used. In the step size invariance case as \(t \rightarrow t_\finalstate\) then the expansions higher order terms start to have stronger effects. This can be seen by
\begin{equation}\label{eq:TimeInverseScale}
    \begin{split}
        \frac{\ddiff^p}{\ddiff t^p}(\phi_{p,q}) &= \frac{1}{\Dt^p} \frac{\ddiff^p}{\ddiff t^p}(\tilde{\phi}_{p,q}),\\
        \phi_{p,q}(\epsilon + t_\initstate) &= \sum_{q=0}^{\infty} \frac{\epsilon^q}{q!}\frac{\ddiff^q}{\ddiff t^q}(\phi_{p,q}(t_\initstate)) = \sum_{q=0}^{\infty} \frac{1}{q!} \frac{\epsilon^q}{\Dt^q}\frac{\ddiff^q}{\ddiff t^q}(\tilde{\phi}_{p,q}(0)),\\
        \epsilon &\in [0, \Dt].
    \end{split}
\end{equation}
Typically, the orders of a numerical method leave the higher order terms \(\mathcal{O}(\Dt^p)\) as error. However, as \(\epsilon \rightarrow \Dt\) every term in the expansion reaches magnitudes of \(\mathcal{O}(1)\) and becomes independent of \(\Dt\). Any decay in higher terms must come from the elementary differentials value and the leading factorials \(1/q!\). Thus, a \BSeriesText which incorporates these terms must handle their labels carefully since as \(\Dt \rightarrow 0\) the approximation does not necessarily improve. Applying TII to such labels allows one to handle all infinite terms and to describe the effect of approximations to them.

\subsection{Tree Iterated Integrals (TII)}
To obtain the exact solution in terms of reduction nodes (\cref{theorem:reductionForceRepresentation}) we use integrals iterated over the trees. To motivate this a brief example is provided. Given the scalar function $g(\tau): \outspace \rightarrow \outspace$ and expanding out into Taylor series around $\tau = 0$ gives
\begin{equation}
    g(\tau) = g(0) + \tau g'(0) + \frac{\tau^2}{2!} g^{(2)}(0) + \frac{\tau^3}{3!} g^{(3)}(0) + \dots  = \sum^{\infty}_{i=0} \frac{\tau^{i}}{i!} g^{(i)}(0).
\end{equation}
Integrating \(g\) w.r.t. \(\tau\) is equivalent to integrating the Taylor series giving
\begin{equation}
    \int_0^{\tau} g(\sigma) \dint\sigma = \sum^{\infty}_{i=0} \frac{\tau^{i+1}}{(i+1)!} g^{(i)}(0).
\end{equation}
Looking back at the Taylor expansion of an ODE it is known that the labels for the exact solution are \(a = \frac{1}{\gamma}\). Taking the reduced expansion \cref{theorem:reductionForceRepresentation} for the tree \(\btree{[gth]}\) and inserting \(\frac{1}{\gamma}\) as the underlying label gives 
\begin{equation}
    \timeLabel{a}{0}{\tau}^{\reductionFun}(\btree{[gth]}) = \sum_{i=0} \frac{\tau^{|\btree{[b]} \ButcherProduct \btree{[t]}^{i}|}}{\gamma(\btree{[b]} \ButcherProduct \btree{[t]}^{i})\sigma(\btree{[b]} \ButcherProduct \btree{[t]}^{i})} F(\btree{[b]} \ButcherProduct \btree{[t]}^{i})(0) = \sum_{i=0} \frac{ \tau^{i+1}}{(i+1)(i!)} g^{(i)}(0) =  \int_0^{\tau} g(\sigma) \dint \sigma
\end{equation}
representing the exact solution.

The concept of integrals representing the exact solution for reduced nodes will be extended for arbitrary trees and functions. This leads to the concept of Tree Iterated Integrals (TII). For clarity a \textbf{Tree Iterative Integral} is defined as an integral which repeats integrations over the nodes of a tree. This is considered a generalization of \textbf{Iterated Integrals}.
\begin{definition}[Tree Iterated Integral (TII) \cite{kreimer1999chen}]\label{def:treeIterativeIntegral}
    Let the function $\treeIterInt{s}{x}(\tree) : (\outspace \times \outspace \times \NTTC) \rightarrow \outspace$ be the TII. Given a set of assignments of functions to nodes \(\reductionFun(\psi)(\tree)\) the Tree Iterative Integral $\treeIterInt{s}{x}^{\reductionFun}(\tree)$ is defined recursively on a tree $\tree = \tau[\tree_1^{k_1}\tree_2^{k_2} \dots \tree_n^{k_n}] \in \NTTC$ as
    \begin{align}
        &\treeIterInt{s}{x}^{\reductionFun}(\tau) = \int_s^{x}\reductionFun(\tau, \emptyset)(u_{\tau}) \dint u_{\tau},\\
        &\treeIterInt{s}{x}^{\reductionFun}(\tree) = \treeIterInt{s}{x}^{\reductionFun}(\tau[\forestTQt]) = \int_{s}^{x} \reductionFun(\tau,\dforest)(u_{\tree}) \prod_{i=1}^{n} (\treeIterInt{s}{u_{\tree}}^{\reductionFun}(\tree_i))^{q_i} \dint u_{\tree},\\
        &\treeIterInt{s}{x}^{\reductionFun}(\emptyset) = 1,
    \end{align}
    where if no function is assigned to a node for \(\reductionFun\) it is assumed to return \(0\). 
    For the case of color compression element-wise evaluation is done.
\end{definition}
The concepts are similar to those of  \cite{brouder2000runge, gubinelli2010ramification, butcher1987numerical}. 

For color compression an alternative formulation is provided. First, we define the symmetric Kronecker power as 
\begin{definition}[Symmetric Kronecker Power]\label{def:symkronprod}
    Given some positive integer \(p > 0\), the sets from \cref{def:IntegerPartitions}, and a vector \(\vec{a} \in \outspace^{n}\) the symmetric Kronecker product is defined as 
    \begin{equation}
        (\vec{a})^{\outersym p} = v \quad \text{ where } \quad  v_{\pi(\sigma)} =  \prod_{j = 1}^{n} (\vec{a}_{j})^{\sigma(j)}, \quad \sigma \in \mathfrak{V}_{p,n}
    \end{equation}
    where \(\pi\) is some arbitrary bijection for the elements \(\sigma\) to the integers \(\{1,\, \dots\, |\mathfrak{V}_{p,n}|\}\).
\end{definition}
\begin{example}[Example: Symmetric Kronecker Power]
    Some examples of the symmetric Kronecker power are 
    \begin{equation}
        \begin{bmatrix}
            a_1 \\
            a_2
        \end{bmatrix}^{\outersym 3} = \begin{bmatrix}
            a_1^3\\
            a_1^2 a_2\\
            a_1 a_2^2\\
            a_2^3
        \end{bmatrix}, \qquad \begin{bmatrix}
            a_1 \\
            a_2 \\
            a_3
        \end{bmatrix}^{\outersym 2} = \begin{bmatrix}
            a_1^2\\
            a_1 a_2\\
            a_1 a_3\\
            a_2^2\\
            a_2 a_3\\
            a_3^2
        \end{bmatrix}.
    \end{equation}
\end{example}
With the new symmetric Kronecker product the compressed TII is given.
\begin{corollary}[Compressed TII]
    An alternative to element-wise TII is the compressed form
    \begin{equation}
        \treeIterInt{s}{x}^{\reductionFun}(\tree) = \treeIterInt{s}{x}^{\reductionFun}(\tau[\forestTQt]) = \int_{s}^{x} \reductionFun(\tau ,\dforest)(u_{t}) \bigotimes_{i=1}^{n} (\treeIterInt{s}{u_{\tree}}^{\reductionFun}(\tree_i))^{\outersym q_i} \dint u_{\tree}\\
    \end{equation}
    where \(\bigotimes\) is the Kronecker product applied repeatedly by index.
    Just as in \cref{def:symkronprod} the order of trees can be permutated as long as indexing is done by trees, thus it is agnostic to the ordering of the results from \(\bigotimes\) and \(\outersym\).
\end{corollary}
\begin{example}
    Multiple examples of \cref{def:treeIterativeIntegral} are as follows.
    Given \(\reductionFun(\btree{[x]})(x) = 1\) and \(\reductionFun(\btree{[gth]})(x) = \phi(x)\) the TII of the tree \(\btree{[x[x][x[gth]]]}\) is 
    \begin{equation}
        \treeIterInt{0}{\Dt}^{\reductionFun}(\btree{[x[x][x[gth]]]}) = \int_0^\Dt u_1 \left(\int_0^{u_1} \left(\int_0^{u_2} \phi(u_3) \dint u_3\right) \dint u_2\right)  \dint u_1.
    \end{equation}
    A more complicated example utilizing compressions is given as follows. Let \(\varrho(\btree{[x]})(x) = 1\) and \(\varrho(\btree{[gth]})(x) = [\varrho(\btree{[b]})(x), \varrho(\btree{[r]})(x)]^T = [\phi_1(x)\,\phi_2(x)]^{T}\). For the tree \(\btree{[x[gth][gth]]}\) the TII is 
    \begin{equation}
        \begin{split}
        \treeIterInt{0}{\Dt}^{\varrho}(\btree{[x[gth][gth]]}) &= \int_0^\Dt \left(\int_{0}^{u_1} \begin{bmatrix}
                    \phi_1(u_2)\\
                    \phi_2(u_2)
                \end{bmatrix} \dint u_2\right)^{\outersym 2} \dint u_1 \\
                &= \int_0^\Dt \begin{bmatrix}
                    (\int_{0}^{u_1} \phi_1(u_2) \dint u_2)^2\\
                    (\int_{0}^{u_1} \phi_1(u_2) \dint u_2)(\int_0^{u_1} \phi_2(u_2) \dint u_2)\\
                    (\int_{0}^{u_1} \phi_2(u_2) \dint u_2)^2\\
                \end{bmatrix} du_1 = \begin{bmatrix}
                    \treeIterInt{0}{\Dt}^{\varrho}(\btree{[x[b][b]]})\\
                    \treeIterInt{0}{\Dt}^{\varrho}(\btree{[x[b][r]]})\\
                    \treeIterInt{0}{\Dt}^{\varrho}(\btree{[x[r][r]]})\\
                \end{bmatrix}.
        \end{split}
    \end{equation}
\end{example}
To further apply the TII some properties are given.
\begin{property}[Iterated Tree Integral Linear Substitution] \label{property:subitt}
    For a tree \(\tree\) and integration bounds \([a, b]\) let each node \(u\) be defined as \(\reductionFun(u)(x)\) over that interval and let \(\varrho(u)(x) = \reductionFun(u)((b-a)x + a)\) be the normalized variation to \([0,1]\). Then by integral substitution
    \begin{equation}
        \treeIterInt{a}{b}^{\reductionFun}(\tree) = (b - a)^{\treeorder{\tree}}\treeIterInt{0}{1}^{\varrho}(\tree).
    \end{equation}
\end{property}
\begin{proof}
    Given the variable transformation \(k(z) = (b-a)z + a\) and \(k^{-1}(z) = \frac{z - a}{b-a}\) and some value \(\epsilon \in [0,1]\) the base integral of TII gives 
    \begin{equation}
        \int_{a}^{k(\epsilon)}  \reductionFun(\tau,\dforest)(u_{\tree}) \prod_{i=1}^{n} (\treeIterInt{a}{u_{\tree}}^{\reductionFun}(\tree_i))^{q_i} \dint u_{\tree} = (b-a)\,\int_{0}^{\epsilon}  \reductionFun(\tau,\dforest)(k(u_{\tree})) \prod_{i=1}^{n} (\treeIterInt{a}{k(u_{\tree})}^{\reductionFun}(\tree_i))^{q_i} \dint u_{\tree}.
    \end{equation}
    Recursively applying this relationship to each node gives the factor \((b-a)\) thus, filtering out gives it raised to the cardinality of the tree.
\end{proof}
By \cref{property:subitt} then all further TII will be normalized around the \([0,1]\) and \([0,\Dt]\) domain for simplicity.
\begin{property}[Linearity and Multi-Colored] \label{property:linitt}
    Let there be a set of nodes \(\unode \in \mathcal{U}\) such that \(\reductionFun(\unode)(x) = \phi_{\unode}(x)\) and a single node \(\vnode \notin \mathcal{U}\) such that \(\reductionFun(\vnode)(x) = \sum_{\unode \in \mathcal{U}}\reductionFun(\unode)(x)\). Let \(\mathcal{S}_{\mathcal{U}, \vnode}(\tree)\) be a subset of trees such that all nodes of \(\vnode\) in \(\tree\) are replaced with all unique permutations of nodes in set \(\mathcal{U}\). The TII of \(\tree\) is changed to 
    \begin{equation}
        \treeIterInt{0}{1}^{\reductionFun}(\tree) = \sum_{\stree \in \mathcal{S}_{\mathcal{U}, v}(\tree)} \frac{\sigma(\tree)}{\sigma(\stree)}\treeIterInt{0}{1}^{\reductionFun}(\stree).
    \end{equation}
\end{property}
\begin{proof}
    For the case of ladder trees this is due to linearity of TII at each level. Thus, only color permutations are necessary to add with no inclusion of symmetry. For the case of tall bush trees \(\TallBushTrees{\NTTC}\) then to expand the tree with nodes of \(\vnode\) each term must be linearly expanded as well. However, this implies expanding \((\sum_{\unode \in \mathcal{U}}\reductionFun(\unode)(x))^{q}\) where \(q\) is the number of duplicate trees in the forest. Thus, expanding the terms out leads to the lower symmetry factor for when different colors are given. This gives the leading symmetry factor ratio. Extending this recursively over the tall bush trees gives the result for all trees.
\end{proof}
\begin{property}[Leading Scalars]\label{property:scalaritt}
    Given some node \(\unode\) the nodal function \(\reductionFun(\unode)(x) = a_{\unode}\phi_{\unode}(x)\) such that \(a_{\unode} \in \outspace\) and a set of functions \(\varrho(\unode)(x) = \phi_{\unode}(x)\). Let \(q(\tree,\unode)\) be a function that returns the number of nodes \(\unode\) that are present in the tree \(\tree\). The TII of \(\tree\) is expanded as 
    \begin{equation}
        \treeIterInt{0}{1}^{\reductionFun}(\tree) = \left[\prod_{\unode \in \mathcal{U}} a_{u}^{q(\tree, u)} \right] \treeIterInt{0}{1}^{\varrho}(\tree).
    \end{equation}
\end{property}
This is given by the integrals pulling constant scalars out. Recursively this pulls all the nodal scalars out and scales by the number in the tree.

Next a property on the domain of integration for TII is provided as follows.
\begin{property}[Domain Splitting TII]\label{property:domainSplitTII}
    Suppose the domain of integration is \(a\) to \(\epsilon\) where \(\epsilon \ge a\) and let \(b > a\) be some other point. Given \\\(\tree = \tau[\tree_1^{q_1}\tree_2^{q_2} \dots \tree_n^{q_n}] \in \NTTC\) the TII over the time domain has cases
    \begin{equation}
        \treeIterInt{a}{\epsilon}^{\reductionFun}(\tree) = \left\{\begin{aligned}
            &\treeIterInt{a}{\epsilon}^{\reductionFun}(\tree) &\textit{ if } b \ge \epsilon,\\
            &\treeIterInt{a}{b}^{\reductionFun}(\tree) + \int_{b}^{\epsilon} \reductionFun(\tau, \dforest)(\sigma) \prod_{i=1}^{n} \treeIterInt{a}{\sigma}^{\reductionFun}(\tree_i)^{q_i} \dint \sigma  &\textit{ if } b < \epsilon \textit{ and } \treeorder{\tree} > 1,\\
            &\treeIterInt{a}{b}^{\reductionFun}(\tau) + \treeIterInt{b}{\epsilon}^{\reductionFun}(\tau) &\textit{ if } b \le \epsilon \textit{ and } \treeorder{\tree} = 1.\\
        \end{aligned}\right.
    \end{equation}
\end{property}
This is given by \cref{def:treeIterativeIntegral} and that integrals can split their integration domain.

Lastly, a relationship is formed from Cauchy's repeated integration formula and the TII.
\begin{lemma}[TII in Cauchy's Repeated Integration]\label{lemma:TIICauchys} 
    Given the nodal function \(\reductionFun'(u)(x) = 1,\, \forall\, u \in \mathcal{U} \subset \mathcal{D}\) and some tree \(\tree_{n}[\fforest] = \psi_{1}[\psi_{2}[\dots \psi_{n}[\fforest]]]\) where it has \(n\) nodes of \(\mathcal{U}\) linearly used up to where the final node has children of  \(\fforest \in \ForestSpace{\NTTC}\). For simplicity let \(\tree_{k}\) be the sequence of nodes up to position \(k\) where \(\tree_{0} = \emptyset\), \(\tree'_{k}\) is the sequence of nodes start at position \(k\) going up to \(n\), and let \(\psi_{k}\) be the node at the \(k\)th position in the chain. By Cauchy's repeated integration a summed relationship is formed as 
    \begin{equation}
        \treeIterInt{0}{\tau}^{\reductionFun'}(\tree_{n}[\fforest]) = \sum_{k=0}^{n-1} \frac{-1^{k}}{k!} \treeIterInt{0}{\tau}^{\reductionFun'}(\tree_{n-1-k})\,\treeIterInt{0}{\tau}^{\reductionFun'}\left(\psi_{1}\left[\left(\prod_{i=2}^{k} \psi_{i}\right) \fforest\right]\right)
    \end{equation}
    or alternatively 
    \begin{equation}
        \treeIterInt{0}{\tau}^{\reductionFun'}(\tree_{n}[\fforest]) = \sum_{k=0}^{n-1} -1^{k} \treeIterInt{0}{\tau}^{\reductionFun'}(\tree_{n-1-k})\, \treeIterInt{0}{\tau}^{\reductionFun'}(\tree'_{k} \ButcherProduct \fforest)
    \end{equation}
    or some construction of scaling such that the product matches \(\frac{\tau^{k}}{k!}\).
\end{lemma}
\begin{proof}
    For simplicity of the forest writing  
    \begin{equation}
        \treeIterInt{0}{\tau}^{\reductionFun'}(\fforest) = G(\tau).
    \end{equation}
    Given the TII of the original tree is 
    \begin{equation}
        \treeIterInt{0}{\tau}^{\reductionFun'}(\tree_{n}[\fforest]) = \int_{0}^{\tau} \dots \int_{0}^{\sigma_{n-1}} G(\sigma_n) \dint \sigma_{n} \dots \dint  \sigma_{1}.
    \end{equation}
    By Cauchy's repeated integration rule 
    \begin{equation}
        = \frac{1}{(n-1)!} \int_{0}^{\tau} (\tau - \sigma_1)^{n-1}\, G(\sigma_1)\dint \sigma_1.
    \end{equation}
    Applying the binomial expansion and simplifying terms gives
    \begin{equation}
        = \sum_{k=0}^{n-1} \frac{-1^{k}}{k!} \frac{\tau^{n-1-k}}{(n-1-k)!} \int_{0}^{\tau} \sigma_1^{k}\, G(\sigma_1)\dint \sigma_1.
    \end{equation}
    Given that any linear tree of size \(k\) has the resulting TII of \(\frac{\tau^k}{k!}\) the equation reduces to the statement of the lemma.
\end{proof}
Returning to the idea of describing the exact solution of \eqref{eq:TFProblem}. The utilization of \cref{def:treeIterativeIntegral} has been done before to represent the exact solution of the non-autonomous problem in \cite{butcher1987numerical}. This was shown in regard to Picard iteration or the analytical solution to an ODE \(\y(t_\finalstate) = \y(t_\initstate) + \int_{t_\initstate}^{t_\finalstate} \fun(\y(\sigma))\, \dint \sigma\) and was mapped directly to the \(\gamma\) function and trees.
\begin{theorem}[Iterated Integral Gamma] \label{theory:iterlabelgamma}
    Given \(\reductionFun(\unode)(x) = 1\) for all nodes \(\unode \in \Tcolors\) and for a tree  $\tree \in \NTTC$ then 
    \begin{equation}
        \frac{(\Dt)^{\treeorder{\tree}}}{\gamma(\tree)} = \treeIterInt{t_\initstate}{t_\finalstate}^{\reductionFun}(\tree).
    \end{equation}
\end{theorem}
\begin{proof}
    Since \cref{property:subitt} for integral transformation and in \cite[Pg. 326]{butcher1987numerical} it was proven this holds for mono-colored trees. For the case of multicolored trees there is no difference as multi-partition problems encapsulate their effects in the \(\frac{1}{\sigma}\) term.
\end{proof}
For each node which has a non-constant labels a relationship between exact and TII is provided.
\begin{theorem}[TII Labels and Exact Solution]\label{theory:ittLabelExact}
    Utilizing the set of colors and functions from \cref{def:reductionNodeFunctions}. Given some tree \(\stree \in \NTTC\) then from \cref{theorem:reductionForceRepresentation} the exact solution via TII is defined as 
    \begin{equation}
        \treeIterInt{t_\initstate}{t_\finalstate}^{\reductionFun}(\stree) = \sigma(\stree) \sum_{\rtree \in \NTTC}  \Dt^{\treeorder{\rtree} - \treeorder{\stree}} \timeElemOp{t_\initstate}^{\reductionFun}(\rtree \replaceNode{\mathcal{U}} \stree) \frac{1}{\gammaFun{\Tcolors \setminus \mathcal{V}}(\rtree)\sigma(\rtree)}.
    \end{equation}
\end{theorem}
The proof is provided in the appendix \cref{proof:ittLabelExact}. 

Lastly taking the reduced expansion representation into the \BSeriesText.
\begin{theorem}[Exact \BSeriesText Flow Reduction] \label{theory:ittflow}
    Taking the set of colors, sets and functions from \ref{def:reductionNodeFunctions} the \BSeriesText formulation with reduced labels is 
    \begin{equation}
        \NBSeries{\Dt}\left(\frac{1}{\gammaFun{\Tcolors \setminus \mathcal{V}}}, \y_\initstate\right) =  \NBSeries{1}(\treeIterInt{t_\initstate}{t_\finalstate}^{\reductionFun}, \y_\initstate).
    \end{equation}
\end{theorem}
\begin{proof}
    This is given by \cref{theorem:reductionForceRepresentation} and the zero entries in \(\reductionFun\) for terms in expansion.
\end{proof}

As an example take Picard iteration for the solution of the ODE 
\begin{equation}
    \begin{split}
        \y' &= \alpha_1(t) \, \fun\comp{1}(\y) + \alpha_2(t) \, \fun\comp{2}(\y), \quad t \in [t_\initstate,\, t_\finalstate],\\
        \y(t) &= \y(t_\initstate) + \int_{t_\initstate}^{t_\finalstate} \alpha_1(\sigma) \, \fun\comp{1}(\y(\sigma)) + \alpha_2(\sigma) \, \fun\comp{2}(\y(\sigma)) \dint\sigma.
    \end{split}
\end{equation}
Repeating the Picard iteration by infinitely substituting into itself \(\y(t)\) then differentiating \(\fun\comp{1}\) and \(\fun\comp{2}\) around \(\y_\initstate\) gives
\begin{equation}
    \scalemath{0.95}{
    \begin{split}
        \y(t) &= \y(t_\initstate) + \int_{t_\initstate}^{t_\finalstate} \alpha_1(\sigma_1) \dint\sigma_1 \, \fun\comp{1} + \int_{t_\initstate}^{t_\finalstate} \alpha_2(\sigma_1) \dint \sigma_1 \, \fun\comp{2} \\
        &+ \quad \int_{t_\initstate}^{t_\finalstate} \alpha_1(\sigma_1) \int_{t_\initstate}^{\sigma_1} \alpha_1(\sigma_2) \dint\sigma_2 \, \dint\sigma_1 \, \fun_{\y}\comp{1}\,\fun\comp{1} + \int_{t_\initstate}^{t_\finalstate} \alpha_1(\sigma_1) \int_{t_0}^{\sigma_1} \alpha_2(\sigma_2) \dint\sigma_2 \, \dint\sigma_1 \, \fun\comp{1}_{\y}\,\fun\comp{2} \\
        &+ \quad \int_{t_\initstate}^{t_\finalstate} \alpha_2(\sigma_1) \int_{t_0}^{\sigma_1} \alpha_1(\sigma_2) \dint\sigma_2 \, \dint\sigma_1 \, \fun_{\y}\comp{2}\fun\comp{1} + \int_{t_\initstate}^{t_\finalstate} \alpha_1(\sigma_1) \int_{t_\initstate}^{\sigma_1} \alpha_1(\sigma_2) \dint\sigma_2 \, \dint\sigma_1 \, \fun\comp{2}_{\y}\,\fun\comp{2} + \dots,\\
    \end{split}
    }
\end{equation}
where \(\reductionFun(\btree{[x]})(x) = \alpha_1(x)\) and  \(\reductionFun(\btree{[z]})(x) = \alpha_2(x)\) giving the reduced \BSeriesText 
\begin{equation}
    \y(t) = \NBSeries{1}(\treeIterInt{t_\initstate}{t}^{\reductionFun}, \y_\initstate).
\end{equation}
A separate example is formed by Runge-Kutta methods with reduced nodes. Looking back at \cref{theorem:RKTime} then this implies the method approximates quadratures and higher quadratures. As an example looking at order conditions up to order \(3\) with time dependent nodes gives
\begin{equation}
    \begin{split}
        b^{T} g_i(t_n + \Dt\, c)  &= \int_{0}^{1} g_i(t_n + \Dt\, \sigma_1)\, \dint \sigma_1 + \bigOh(\Dt^{p - 1}), \\
        b^T A\, g_i(t_n + \Dt\, c)  &= \int_{0}^{1} \int_{0}^{\sigma_1} g_i(t_n + \Dt\, \sigma_2)\, \dint \sigma_2\, \dint \sigma_1 + \bigOh(\Dt^{p - 2}),\\
        b^T (c \Elemtimes g_i(t_n + \Dt\, c))  &= \int_{0}^{1} \sigma_1\, g_i(t_n + \Dt\, \sigma_1)\, \dint \sigma_1 + \bigOh(\Dt^{p - 2}),\\
        b^T  (g_i(t_n + \Dt c)\,g_j(t_n + \Dt\, c))   &= \int_{0}^{1} g_i(t_n + \Dt\, \sigma_1)\,g_j(t_n + \Dt\, \sigma_1)\, \dint \sigma_1 + \bigOh(\Dt^{p - 2})\\
        \forall\, i,\, j \in \{1 .. m\}.
    \end{split}
\end{equation}

\subsection{Composition Utilizing Tree Iterated Integrals (TII)}

The new TII are utilized to expand the definition of composition of exact labels.
\begin{theorem}[Integral Analytical Composition]\label{theory:analyticalComposition}
    Given the time points \(a, b, c \in \outspace\), trees \(\tree, \utree \in \NTTC\), and some arbitrary \(\reductionFun(\psi)(x)\) then
    \begin{equation}
        \begin{split}
            \treeIterInt{a}{c}^{\reductionFun}(\tree \prune \utree) &= \sum_{\stree \in \NTTC} \treeIterInt{a}{b}^{\reductionFun}(\tree \prune \stree) \,\treeIterInt{b}{c}^{\reductionFun}(\stree \prune \utree),\\
            \treeIterInt{a}{c}^{\reductionFun}(\pruneMatrix) &= \treeIterInt{a}{b}^{\reductionFun}(\pruneMatrix)\, \treeIterInt{b}{c}^{\reductionFun}(\pruneMatrix).
        \end{split}
    \end{equation}
\end{theorem}
A proof of this is provided in \cref{proof:analyticalComposition}.
\begin{remark}[Discontinuous Basis Functions]
    The proof of \cref{theory:analyticalComposition} relies only on TII and decomposing its integration domain. This leaves basis functions \(\reductionFun(\psi)(x)\) that can be pointwise discontinuous. As long as the piecewise functions are defined as \(t_{i-1} \le \tau \le t_{i}\) for each interval and not pointwise defined.
\end{remark}
Looking back at the properties derived for \(\frac{1}{\gamma}(\pruneMatrix)\) in \cref{theory:pruneProperties} and \cref{theory:repeatedIntegration}. Such properties derived from matrices are not as easily applied to TII due to there not being powers of TII formed. For example let \(\reductionFun(\psi)(x) = \phi(x)\), \(\varrho_{a}(\psi)(x) = \phi(\Dt\,x + a)\) and suppose there are \(\{t_\initstate,\,t_{\initstate + 1},\, \dots \,t_{\initstate + m}\}\) time points evenly subdivided by \(\Dt\) into \(m\) sub-intervals then
\begin{equation}
    \begin{split}
        \treeIterInt{t_\initstate}{t_\initstate + \Dt m}^{\reductionFun}(\pruneMatrix) &= \DiagMat_{\Dt}\, \treeIterInt{0}{m}^{\varrho_{t_\initstate}}(\pruneMatrix)\, \DiagMat_{\Dt}^{-1} = \DiagMat_{\Dt}\, \left[\prod_{i=0}^{m-1} \treeIterInt{i}{i+1}^{\varrho_{t_\initstate}}(\pruneMatrix)\right]\, \DiagMat_{\Dt}^{-1}\\
        &= \DiagMat_{\Dt}\, \left[\prod_{i=0}^{m-1} \treeIterInt{0}{1}^{\varrho_{t_n + \Dt\, i}}(\pruneMatrix)\right]\, \DiagMat_{\Dt}^{-1} = \prod_{i=0}^{m-1} \treeIterInt{t_n + \Dt\, i}{t_n + \Dt\, (i+1)}^{\reductionFun}(\pruneMatrix).\\
    \end{split}
\end{equation}
This leaves no powers of the matrix being formed since at least either the bounds of the integrals change or the function definitions do. Assuming at least one node of the function \(\reductionFun\) is not the constant then \cref{theory:pruneProperties} is \textbf{not directly} applicable on TII form. However, using the underlying representation \cref{theory:ittLabelExact} provides a way around this. All sub-block terms relating to the reduced labels are some form of Pascal's triangle.
Given such a view the error accumulated in each compressed function over repeated compositions is provided.
\begin{theorem}[Repeated Composition Integral Form]\label{theory:compIntegral}
    Taking the set of colors, sets and functions from \cref{def:reductionNodeFunctions} and definitions from \cref{theory:TriangleMatrixClose}.
    Given the label \(a(\tree) = \frac{1}{\gammaFun{\Tcolors \setminus \mathcal{V}}(\tree)} \quad \forall \, 0 \le \treeorder{\tree} \le p\) over the time interval \([t_\initstate,t_\finalstate]\) and composing the method \(m\) times then the composed reduced label \(a^{m, \reductionFun}_{t_\initstate, t_\finalstate}\) has the difference from the exact solution as 
    \begin{equation}
        \begin{split}
            a^{m, \reductionFun}_{t_\initstate, t_\finalstate}(\stree \prune \utree) - \treeIterInt{t_\initstate}{t_\finalstate}^{\reductionFun}(\stree \prune \utree) &= \sigma(\stree \prune \utree) \sum_{\rtree,\vtree \in \NTTC}  \Dt^{\treeorder{\rtree \prune \vtree} - \treeorder{\stree \prune \utree}}\, \timeElemOp{t_\initstate}^{\reductionFun}((\rtree \replaceNode{\mathcal{U}} \stree) \prune (\vtree \replaceNode{\mathcal{U}} \utree))\\
            &\frac{\delta_{\rtree, \Tcolors \setminus \mathcal{V}}\,\delta_{\vtree, \Tcolors \setminus \mathcal{V}}}{\sigma(\rtree \prune \vtree)} \left[\sum_{k=p}\frac{a((\CompLogBlock_{k+1})_{\rtree,\vtree})}{m^k}\right]
        \end{split}
    \end{equation}
    where \(\delta_{\tree, \mathcal{C}}\) is one if the tree \(\tree\) has only colors in the set \(\mathcal{C}\) and otherwise zero. The pruning applied to replacement results only act on the two results if they share the same right tree bracket entry.
\end{theorem}
\begin{proof}
    The underlying representations are given by \cref{theory:ittLabelExact} and \cref{theorem:reductionForceRepresentation} where a difference is formed of each label. The result is given by inserting into \cref{theory:repeatedIntegration}.
\end{proof}
This gives the error generated in each coefficient label as the method is repeated for multiple steps. Such results can be utilized to create the underlying repeated composition matrix results.

\section{Conclusion}
\label{section:Conclusion}

In this paper, we have introduced and redefined a multitude of operations on \BSeriesText in order to solve the approximate problem \eqref{eq:ApproximationProblem}. By introducing the matrix permanent we redefined tree pruning and cutting giving new formulations for composition and substitution while avoiding ordered trees. With the matrix analysis provided in \cref{sec:MatrixAnalysis} we established new properties of composition. Such properties give an expansion in terms of uniform steps and a new way to derive the labels of a backward differential equation.

We show that the permanent in tandem with the generalized outer product can be used to describe multiple tree recursive binary operations. With these results, we were able to introduce a new operation called replacement for contracting the colors of nodes into others. Thus, giving a systematic way to expand out elementary differential operators which are composed of \BSeriesText themselves. This allows integration methods to use approximate flows constructed by constant frozen vector fields summed and scaled with time as shown in \eqref{eq:modifiedAsymptoticConstant}.

To push the new results to a continuous framework we introduce a variation of \BSeriesText where the labels are time dependent. We utilize this series to analyze the general approximate functions \eqref{eq:modifiedAsymptotic}. This allows the filtering of time dependent scalar terms in elementary differentials to be represented in the labels instead. To support this variation the exact solution is given as a Tree Iterated Integral with multiple associated properties.

The authors plan to apply the operations defined in this paper to derive new methods beyond the ones that are mentioned here.

\section*{Acknowledgement}
\small
This work was funded by the National Science Foundation (NSF) award DMS-2411069 and by the Computational Science Laboratory. AI was used to assist with grammar, spell checking, a latex command for math text displays, and catching minor formatting issues in the text such as broken equation references. The authors have checked and edited all suggestions from the AI for correctness and take full responsibility for the content.

\normalsize

\bibliographystyle{plainnat}
\bibliography{./main_theory}

\appendix

\section{Block Triangular Matrices\label{sec:MatrixAnalysis}}
In this section we provide results on block triangular matrices to support results in the main body of the paper. All analysis provided here is self-contained and independent of the main body of the paper. We only target block triangular matrices with restrictions on the diagonal blocks.
Let \(\Lmat\) be some block lower-triangular matrix such that it has the structure
\begin{equation}
    \Lmat = \begin{bmatrix}
        \blockIndex{\Lmat}{0,0} & \matzero & \matzero & \matzero & \matzero & \dots\\
        \blockIndex{\Lmat}{1,0} & \blockIndex{\Lmat}{1,1} & \matzero & \matzero & \matzero & \dots\\
        \blockIndex{\Lmat}{2,0} & \blockIndex{\Lmat}{2,1} & \blockIndex{\Lmat}{2,2} & \matzero & \matzero & \dots\\
        \blockIndex{\Lmat}{3,0} & \blockIndex{\Lmat}{3,1} & \blockIndex{\Lmat}{3,2} & \blockIndex{\Lmat}{3,3} & \matzero & \dots\\
        \blockIndex{\Lmat}{4,0} & \blockIndex{\Lmat}{4,1} & \blockIndex{\Lmat}{4,2} & \blockIndex{\Lmat}{4,3} & \blockIndex{\Lmat}{4,4} & \dots\\
        \vdots & \vdots & \vdots & \vdots & \vdots & \ddots\\
    \end{bmatrix}
\end{equation}
where \(\blockIndex{\Lmat}{i,j}\) is a rectangular matrix of size \(n_i \times n_j\). Every element of the matrix must be from the same commutative ring.
Note that throughout this section (unless specified) it will be assumed that the eigenvalues of this \(\Lmat\) have the properties: \(\Re(\lambda) > 0\) and the spectral radius \(\rho(\Lmat - \Identity) < 1\).

The beginning of this section will leave the diagonal matrices arbitrary. However, the second section will proceed with assuming \(\blockIndex{\Lmat}{i,i} = \Identity \quad \forall\, i\). The main focus will be on identity diagonal blocks with supporting results derived for arbitrary diagonals.
Another matrix that is used is the block diagonal scaling matrix defined as
\begin{equation}
    \DiagMat_{\alpha} = \textrm{Diag}([\underbrace{1,\, \dots,\, 1}_{n_0},\, \underbrace{\alpha,\, \dots,\, \alpha}_{n_1},\, \underbrace{\alpha^2,\, \dots,\, \alpha^2}_{n_2},\, \dots,\, \underbrace{\alpha^i,\, \dots,\, \alpha^i}_{n_i},\, \dots]).
\end{equation}
Utilized throughout are the following sets of integer partitioning.
\begin{definition}[Integer Partitions]\label{def:IntegerPartitions}
    The set of all partitions of an integer \(p\) into exactly \(n\) integers has the following constrained sets
    \begin{equation}
        \begin{split}
            & \textbf{Set of partitions where order matters: } \\
            & \sigma \in \IntPartOrder{p}{n} \quad \text{s.t.} \quad \sum_{k=1}^{n} \sigma(k) = p, \quad \sigma(k) \ge 1.\\
            & \textbf{Set of partitions where order does not matters: }\\
            & \tree_{1}^{q_1} \tree_{2}^{q_2} \dots \tree_{m}^{q_m} \in \IntPartUnorder{p}{n} \quad \text{s.t.} \quad \sum_{i=1}^{m} i\, q_i = p, \qquad \sum_{i=1}^{m} q_i = n.\\
            & \textbf{Set of partitions where order matters + zeros: } \\
            & \sigma \in \IntPartOrderZero{p}{n} \quad \text{s.t.} \quad \sum_{k=1}^{n} \sigma(k) = p, \quad \sigma(k) \ge 0.\\
        \end{split}
    \end{equation}
    Where \(\sigma\) is a mapping of an index to the integer at that position. \(\tree_i\) are ladder trees of size \(i\) and \(q_i\) are positive integers. Note we use trees here but, other characters which are commutative are valid.
    An additional note is
    \begin{equation}
        \begin{split}
            \forall\, n > p &\quad \IntPartOrder{p}{n} = \IntPartUnorder{p}{n} = \emptyset,\\
            \forall\, n,\,p & \quad \IntPartOrder{p}{n} \subseteq \IntPartOrderZero{p}{n}.
        \end{split}
    \end{equation}
    For simplicity throughout let \(\overline {\sigma(i)} = \sum_{j=1}^{i} \sigma(j)\).
\end{definition}

As an example of all the constrained integer partition sets the following is given.
\begin{example}[Examples of Sets]
    For \(p = 3, n = 2\) the sets \(\IntPartOrder{p}{n}, \IntPartUnorder{p}{n}, \IntPartOrderZero{p}{n}\) are
    \begin{equation}
        \begin{split}
            \IntPartOrder{p}{n} = & \{\{1,\,2\},\, \{2,\,1\}\},\\
            \IntPartUnorder{p}{n} =& \{\{t_1\, t_2\}\},\\
            \IntPartOrderZero{p}{n} =& \{\{1,\,2\},\, \{2,\,1\},\, \{0,\, 3\}, \{3,\, 0\}\}.
        \end{split}
    \end{equation}
\end{example}
For simplicity going forward the following is shown.
\begin{corollary}
    Suppose there is some function \(f\) which takes as input an integer partitioning and suppose the arguments of the function are commutative. Given that multi-choose is \({n \choose q_1, q_2, \dots q_m} = \frac{n!}{\prod_{i=1}^{m}q_i!}\) then the function is written over the sum of inputs as
    \begin{equation}
        \sum_{\sigma \in \IntPartOrder{p}{n}} f(\sigma) = \sum_{t_1^{q_1} t_2^{q_2} \dots t_m^{q_m} \in \IntPartUnorder{p}{n}}  {n \choose q_1, q_2, \dots q_m} \,f(t_1^{q_1} t_2^{q_2} \,\dots t_n^{q_m})
    \end{equation}
    where \(f\) evaluated at the forest is some choice of integer partitioning.
\end{corollary}
The Stirling numbers of the first and second kind will be utilized throughout.
\begin{definition}[Stirling Numbers] \label{def:Stirling} \cite{comtet2012advanced,wilf2005generatingfunctionology}
    The Stirling numbers of the first and second kind are
    \begin{equation}
            \textbf{First Kind:} \quad s(n,k), \quad \textbf{Second Kind:} \quad S(n,k).
    \end{equation}
    Properties for an arbitrary positive integer \(p\) are 
    \begin{equation}
        \begin{split}
            S(p,1) &= 1,\\
            s(p,p) &= 1,\\
            S(p,n) &= p!  \sum_{t_{1}^{q_1} t_{2}^{q_2} \dots t_{m}^{q_{m}} \in \IntPartUnorder{p}{n}} \prod_{i=1}^{m} \frac{1}{q_i! (i !)^{q_i}},\\
            \delta_{p,k} &= \sum_{n=k}^{p} s(n,k)  S(p,n).\\
        \end{split}
    \end{equation}
\end{definition}

\subsection{Block Triangular Matrix Lemmas}
We provide some results on triangular matrices to support the main results.
First a lemma on summation and Stirling numbers.
\begin{lemma}\label{lemma:combStirlingLemma}
   For \(k \le p \in \mathbb{Z}^+\) then
    \begin{equation}\label{eq:combstir}
        \sum_{n=k}^{p} \frac{s(n,k)}{n!} \left[ \sum_{t_{1}^{q_1} t_{2}^{q_2} \dots t_{m}^{q_m} \in \IntPartUnorder{p}{n}}  {n \choose q_1, q_2, \dots q_m} \prod_{i=1}^{m}\frac{1}{(i!)^{q_i}} \right] = \delta_{p,k}\, \frac{1}{p!}.
    \end{equation}
    If the first \(k < r \le p\) terms of the expansion are removed then
    \begin{equation}
        \begin{split}
            &\sum_{n=k + r}^{p} \frac{s(n,k)}{n!} \left[ \sum_{t_{1}^{q_1} t_{2}^{q_2} \dots t_{m}^{q_m} \in \IntPartUnorder{p}{n}}  {n \choose q_1, q_2, \dots q_m} \prod_{i=1}^{m}\frac{1}{(i!)^{q_i}} \right] \\
            &=\frac{1}{p!} \left(\delta_{p,k} - \sum_{n=k}^{k+r} s(n,k)\,S(p,n) \right).
        \end{split}
    \end{equation}
\end{lemma}
\begin{proof}
    Given the properties of Stirling numbers of the second kind from \cref{def:Stirling} then
    \begin{equation}
        S(p,n) =  p!\,  \sum_{t_{1}^{q_1} t_{2}^{q_2} \dots t_{p - n + 1}^{q_{p - n + 1}} \in \IntPartUnorder{p}{n}} \prod_{i=1}^{p-n+1} \frac{1}{q_i!\, (i!)^{q_i}}.
    \end{equation}
    Scaling both sides by \(\frac{1}{p!}\), re-writing as a multi-choose, and inserting into \eqref{eq:combstir} gives
    \begin{equation}
        \frac{1}{p!} \left(\sum_{n=k}^{p} s(n,k)\, S(p,n)\right).
    \end{equation}
    Given the two kinds of Stirling numbers are inverses of each other (\cref{def:Stirling}) then
    \begin{equation}
        \delta_{p,k}\, \frac{1}{p!}.
    \end{equation}
    For the case of when \(r\) terms are omitted. This is shown by adding the remaining terms to both sides of the equation. Applying the Stirling numbers equality gives a subtraction on the RHS of \(-\sum_{n=k}^{k+r} s(n,k)S(p,n)\).
\end{proof}
Next a lemma on the multiplication powers of the triangular matrix provided.
\begin{lemma}\label{lemma:TriangleMatrixLemmaMult}
    For all \(a > b \in \mathbb{Z}^+\) and \(n \in \mathbb{Z}^+ \setminus \{0\}\) then the following holds for the block triangle matrix powers of \(\Lmat\):
    \begin{equation}
        \blockIndex{(\Lmat^n)}{a,b} = \sum_{\sigma \in \mathfrak{E}_{a - b,n}} \prod_{i = 1}^{n} \blockIndex{\Lmat}{a - \overline{\sigma}(i-1),a - \overline{\sigma}(i)}
    \end{equation}
    where 
    \[\mathfrak{E}_{p,n} = \begin{cases}
        \IntPartOrder{p}{n} \text{ if } \blockIndex{\Lmat}{i,i} = \matzero \quad  \forall\, i,\\
        \IntPartOrderZero{p}{n} \text{ otherwise.}
    \end{cases}\]
\end{lemma}
\begin{proof}
    Let there be some other matrix with the same block structure and size as matrix \(\Lmat\) given as
    \begin{equation}
        \mathbf{B} = \begin{bmatrix}
            \blockIndex{\mathbf{B}}{0,0} & \matzero & \matzero & \dots\\
            \blockIndex{\mathbf{B}}{1,0} & \blockIndex{\mathbf{B}}{1,1} & \matzero & \dots\\
            \blockIndex{\mathbf{B}}{2,0} & \blockIndex{\mathbf{B}}{2,1} & \blockIndex{\mathbf{B}}{2,2} & \dots\\
            \vdots & \vdots & \vdots & \ddots\\
        \end{bmatrix}.
    \end{equation}
    The blockwise entry multiplication is given as 
    \begin{equation*}
        (\mathbf{B}\,\Lmat)_{a,b} = \sum_{i_1 = b}^{a} \blockIndex{\mathbf{B}}{a, i_1} \blockIndex{\Lmat}{i_1, b}.
    \end{equation*}
    Suppose that the matrix \(\mathbf{B} = \mathbf{C} \Lmat\) then inserting back into the block sum gives
    \begin{equation*}
        (\mathbf{B}\Lmat)_{a,b} = \sum_{i_1 = b}^{a} (\mathbf{C}\Lmat)_{a, i_1} \blockIndex{\Lmat}{i_1, b} = \sum_{b \le i_1 \le i_2 \le a} \blockIndex{\mathbf{C}}{a, i_2} \blockIndex{\Lmat}{i_2, i_1} \blockIndex{\Lmat}{i_1, b}.
    \end{equation*}
    Repeating for integer powers of \(n\) gives
    \begin{equation*}
        \blockIndex{(\Lmat^{n})}{a,b} = \sum_{b \le i_1 \le \dots \le i_{n-1} \le a} \blockIndex{\Lmat}{a,i_{n-1}} \left[\prod_{j=1}^{n-2} \blockIndex{\Lmat}{i_{n - j}, i_{n - j - 1}}\right] \blockIndex{\Lmat}{i_{1}, b}.
    \end{equation*}
    Next instead of indexing via increasing indices they can be instead shifted by a repeated subtraction of previous entries. This gives
    \begin{equation*}
        i_{j} = a - \sum_{k=1}^{n - j} s_k
    \end{equation*}
    where \(s_k \ge 0\) are nonnegative integers where the total sum \(\sum_{k=1}^{n} s_{k} = b\). Giving an integer partitioning of \(a - b\) into \(n\) integers of \(s_k = \sigma(k)\) and \(i_{j} = a - \overline{\sigma}(n-j)\) where \(\sigma \in \mathfrak{V}_{a-b,n}\).
    Now suppose that \(\blockIndex{\Lmat}{i,i} = \matzero\, \forall\, i\) then any equality between successive indices \(i_{j}\) gives a zero matrix. Thus, simplifying the equation to strictly increasing indices gives integer partitions with no zeros giving the set \(\mathfrak{W}_{a-b,n}\).
\end{proof}
As a demonstration for the powers of \(2\) and \(3\) where it is assumed the diagonals \(\blockIndex{\Lmat}{i,i} = \Identity \, \forall\, i\) then
\begin{equation}
    \scalemath{0.9}{
    (\Lmat - \Identity)^2 = \begin{bmatrix}
        \matzero & \matzero & \matzero & \matzero & \matzero & \dots\\
        \matzero & \matzero & \matzero & \matzero & \matzero & \dots\\
        \blockIndex{\Lmat}{2,1}\blockIndex{\Lmat}{1,0} & \matzero & \matzero & \matzero & \matzero & \dots\\
        \blockIndex{\Lmat}{3,1}\blockIndex{\Lmat}{1,0} + \blockIndex{\Lmat}{3,2}\blockIndex{\Lmat}{2,0} & \blockIndex{\Lmat}{3,2}\blockIndex{\Lmat}{2,1} & \matzero & \matzero & \matzero & \dots\\
        \blockIndex{\Lmat}{4,1}\blockIndex{\Lmat}{1,0} + \blockIndex{\Lmat}{4,2}\blockIndex{\Lmat}{2,0} + \blockIndex{\Lmat}{4,3}\blockIndex{\Lmat}{3,0} & \blockIndex{\Lmat}{4,2}\blockIndex{\Lmat}{2,1} + \blockIndex{\Lmat}{4,3}\blockIndex{\Lmat}{3,1} & \blockIndex{\Lmat}{4,3}\blockIndex{\Lmat}{3,2} & \matzero & \matzero & \dots\\
        \vdots & \vdots & \vdots & \vdots & \vdots & \ddots\\
    \end{bmatrix},
    }
\end{equation}
\begin{equation}
    \scalemath{0.9}{
    (\Lmat - \Identity)^3 = \begin{bmatrix}
        \matzero & \matzero & \matzero & \matzero & \matzero & \dots\\
        \matzero & \matzero & \matzero & \matzero & \matzero & \dots\\
        \matzero & \matzero & \matzero & \matzero & \matzero & \dots\\
        \blockIndex{\Lmat}{3,2}\blockIndex{\Lmat}{2,1}\blockIndex{\Lmat}{1,0} & \matzero & \matzero & \matzero & \matzero & \dots\\
        \blockIndex{\Lmat}{4,2}\blockIndex{\Lmat}{2,1}\blockIndex{\Lmat}{1,0} + \blockIndex{\Lmat}{4,3}\blockIndex{\Lmat}{3,1}\blockIndex{\Lmat}{1,0} + \blockIndex{\Lmat}{4,3}\blockIndex{\Lmat}{3,2}\blockIndex{\Lmat}{2,0} & \blockIndex{\Lmat}{4,3}\blockIndex{\Lmat}{3,2}\blockIndex{\Lmat}{2,1} & \matzero & \matzero & \matzero & \dots\\
        \vdots & \vdots & \vdots & \vdots & \vdots & \ddots\\
    \end{bmatrix}.
    }
\end{equation}

An alternative view of \cref{lemma:TriangleMatrixLemmaMult} is the integer partition representations of powers. For simplicity the matrices are filled with their integer partition representation. Given the sub-matrix indices \(\blockIndex{\Lmat}{i,j}\) and filling the entries with \(i-j\) gives
\begin{equation}
    \scalemath{0.8}{
    \text{pow-1}: \begin{bmatrix}
       0 & 0 & 0 & 0 & 0 & \dots\\
       \{1\} & 0 & 0 & 0 & 0 & \dots\\
       \{2\} & \{1\} & 0 & 0 & 0 & \dots\\
       \{3\} & \{2\} & \{1\} & 0 & 0 & 0 & \dots\\
       \{4\} & \{3\} & \{2\} & \{1\} & 0 & \dots\\
       \vdots & \vdots & \vdots & \vdots & \vdots & \ddots\\
   \end{bmatrix}
    },
\end{equation}
\begin{equation}
    \scalemath{0.8}{
     \text{pow-2}: \begin{bmatrix}
        0 & 0 & 0 & 0 & 0 & \dots\\
        0 & 0 & 0 & 0 & 0 & \dots\\
        \{1,1\} & 0 & 0 & 0 & 0 & \dots\\
        \{2,1\} + \{1,2\} & \{1,1\} & 0 & 0 & 0 & \dots\\
        \{3,1\} + \{2,2\} + \{1,3\} & \{2,1\} + \{1,2\} & \{1,1\} & 0 & 0 & \dots\\
        \vdots & \vdots & \vdots & \vdots & \vdots & \ddots\\
    \end{bmatrix}
    },
\end{equation}
\begin{equation}
    \scalemath{0.8}{
    \text{pow-3}: \begin{bmatrix}
        0 & 0 & 0 & 0 & 0 & \dots\\
        0 & 0 & 0 & 0 & 0 & \dots\\
        0 & 0 & 0 & 0 & 0 & \dots\\
        \{1,1,1\} & 0 & 0 & 0 & 0 & \dots\\
        \{2,1,1\} + \{1,2,1\} + \{1,1,2\} & \{1,1,1\} & 0 & 0 & 0 & \dots\\
        \vdots & \vdots & \vdots & \vdots & \vdots & \ddots\\
    \end{bmatrix} 
    }.
\end{equation}
Every integer power \(m\) must contain a set size of \(m\) in each entry or \(0\) (null set) demonstrating \(\mathfrak{W}\) sets. If the case of \(\blockIndex{\Lmat}{i,i} \ne \Identity \quad \forall\, i\) is given then none of the sets below and including the diagonal resolve to null sets. The integer partitions would still apply; however, this would constitute sets of \(\mathfrak{V}\) as they are never terminating.

Lastly, a lemma on the matrix \(\log\) and fractional power for the given structure.
\begin{lemma}\label{lemma:MatrixBlockLogPower}
    Given some real number \(m \in \mathbb{R}\) and for all \(a > b \in \mathbb{Z}^+\) then
    \begin{equation}
        \frac{\log(\Lmat)^k}{k!} = \sum_{n=k}^{\infty} \frac{s(n,k)}{n!} (\Lmat - \Identity)^{n}.
    \end{equation}
    The blockwise matrix log is defined as
    \begin{equation}
        \blockIndex{\left(\frac{\log(\Lmat)^k}{k!}\right)}{a,b} = \sum_{n=k}^{\infty} \frac{s(n,k)}{n!} \left[ \sum_{\sigma \in \mathfrak{E}_{a - b,n}} \prod_{i = 1}^{n} (\blockIndex{\Lmat}{a - \overline {\sigma}(i-1), a - \overline {\sigma}(i)} - \delta_{0, \overline {\sigma}(i)}\, \Identity) \right]
    \end{equation}
    where
    \[\mathfrak{E}_{p,n} = \begin{cases}
        \IntPartOrder{p}{n} \text{ if } \blockIndex{\Lmat}{i,i} = \Identity \quad \forall\, i,\\
        \IntPartOrderZero{p}{n} \text{ otherwise}.
    \end{cases}\]
    This gives that any block triangular matrix raised to a real power has a block entry defined as 
    \begin{equation}
        \blockIndex{(\Lmat^{m})}{a,b} = \sum_{k=0}^{\infty} m^k \left[\sum_{n=k}^{\infty} \frac{s(n,k)}{n!}\left[ \sum_{\sigma \in \mathfrak{E}_{a - b,n}} \prod_{i = 1}^{n} (\blockIndex{\Lmat}{a - \overline {\sigma}(i-1), a - \overline {\sigma}(i)} - \delta_{0, \overline {\sigma}(i)} I) \right]\right].
    \end{equation}
\end{lemma}
\begin{proof}
    The matrix \(\log\) power derivation here is similar to the one found in \cite{comtet2012advanced,wilf2005generatingfunctionology} except they focus on the scalar case.
    To start note that the expansion for the powers is 
    \begin{equation}
        \Lmat^{m} = e^{m \log(\Lmat)} = \sum_{k=0}^{\infty} m^{k}\, \frac{\log(\Lmat)^{k}}{k!}.
    \end{equation}
    Given the binomial expansion and utilizing the Stirling numbers \cref{def:Stirling} gives
    \begin{equation}
        \begin{split}
            \Lmat^{m} = (\Identity + (\Lmat - \Identity))^m &= \sum_{n=0}^{\infty} {m \choose n} (\Lmat - \Identity)^{n}\\
             &= \sum_{n=0}^{\infty} \frac{(\Lmat - \Identity)^n}{n!} \sum_{k=0}^{n} s(n,k) \,m^{k} \\
             &= \sum_{k=0}^{\infty} m^{k} \sum_{n=k}^{\infty} \frac{s(n,k)}{n!} (\Lmat - \Identity)^{n}.
        \end{split}
    \end{equation}
    Since both series expansions have the same leading sum of \(\sum_{k=0}^{\infty} m^{k}\) then taking the insides gives
    \begin{equation}\label{eq:submercator}
        \frac{\log(\Lmat)^k}{k!} = \sum_{n=k}^{\infty} \frac{s(n,k)}{n!} (\Lmat - \Identity)^{n}.
    \end{equation}
    For the scalar case this is a Mercator series giving a radius of convergence of  \(|x - 1| < 1\). Given the series \eqref{eq:submercator} is alternating via the Stirling numbers signs then the alternating series test gives the scalar function converges to the same radius. Therefore, since the assumption is that the matrix satisfies \(\rho(\Lmat - \Identity) < 1\) and all eigenvalues have positive real parts then the matrix form converges \cite[Theorem 4.7]{higham2008functions}. Since the diagonals of the matrix \(\rho(\blockIndex{\Lmat}{i,i} - \Identity) < 1\) then all block formulations are convergent.\par
    For the block matrix entry \(\blockIndex{\left(\frac{\log(\Lmat)^k}{k!}\right)}{a,b}\) given \cref{lemma:TriangleMatrixLemmaMult} applied to the matrix \((\Lmat - \Identity)^{n}\) then the product inside must give terms of \(\blockIndex{\Lmat}{i,i} - \Identity\) for when a zero integer is given in the partition sequence. Thus, giving the term \(-\delta_{0, \overline {\sigma}(i)}\, \Identity\) added. Thus, inserting the result into the equation \(\log\) power gives the block formulation.
    The real power block form is given by the case of inserting the block entry of \(\log\) into the matrix exponential.
\end{proof}
We note that the same results are obtained from the analysis in \cite{kurbatovGreensFunctionProblem2018} via a series expansions of the diagonal terms.

\subsection{Main Results for Block Triangular Matrix (Identity Diagonals)\label{section:MatId}}
Utilizing the previously defined lemmas, the following theorems are provided for the matrix \(\Lmat\) under the condition \(\blockIndex{\Lmat}{i,i} = \Identity \quad \forall\, i\).
\begin{theorem}\label{theory:TriangleMatrix}
    Let \(m \in \mathbb{R}\), \(a > b \in \mathbb{Z}^+\) and \(\blockIndex{\Lmat}{i,i} = \Identity \quad \forall\, i\) then the matrix \(\Lmat\) satisfying one of the following properties are equivalent to the others
    \begin{equation}
        \begin{split}
            (\mathcal{A})\text{ : } \qquad& \DiagMat_{m}\, \Lmat\, \DiagMat^{-1}_{m} = \Lmat^{m},\\
            (\mathcal{B})\text{ : } \qquad& \blockIndex{\Lmat}{a,b} = \frac{1}{(a-b)!} \prod_{i=0}^{a-b-1} \blockIndex{\Lmat}{a - i, a - i - 1} \quad \forall\, a > b,\\
            (\mathcal{C})\text{ : } \qquad& \blockIndex{\log(\Lmat)}{i,j} = \begin{cases}
                \blockIndex{\Lmat}{i,j} &\quad \textit{ if } i = j - 1\\
                \matzero &\quad \textit{o.w.}
            \end{cases},\\
            &(\mathcal{A}) \equiv (\mathcal{B}) \equiv (\mathcal{C}).
        \end{split}
    \end{equation}
    Furthermore, given some \(\delta \in \mathbb{R} \setminus \{0,1,-1\}\), if the condition
    \begin{equation}
        \DiagMat_{\delta}\, \Lmat\, \DiagMat^{-1}_{\delta} = \Lmat^{\delta}
    \end{equation}
    holds this implies that \((\mathcal{A})\) must be true for the given matrix.
\end{theorem}
\begin{proof}
    To start taking the real matrix power equation from \cref{lemma:MatrixBlockLogPower}. Simplifying since the assumption of the diagonals being identity gives finite terminating \(\log\) power and terminating matrix exponential thus, giving
    \begin{equation}\label{eq:tricondsum}
        \blockIndex{(\Lmat^{m})}{a,b} = \sum_{k=1}^{a-b} m^k \left[\sum_{n=k}^{a-b} \frac{s(n,k)}{n!}\left[ \sum_{\sigma \in \IntPartOrder{a - b}{n}} \prod_{i = 1}^{n} \blockIndex{\Lmat}{a - \overline {\sigma}(i-1), a - \overline {\sigma}(i)} \right]\right].
    \end{equation}
    This is the main form which will be utilized throughout each part of the proof. The case of \((\mathcal{A}) \equiv (\mathcal{B})\) will be shown with the conditional case shown in conjunction. Next it will be shown that \((\mathcal{B}) \equiv (\mathcal{C})\) thus, by equivalency giving the main result. 

    \textbf{Case} \((\mathcal{A}) \equiv (\mathcal{B})\):\\
    This case will be proven via induction. To start the base case is the first lower sub-diagonal giving \(b = a-1\) applied to \eqref{eq:tricondsum} \(\blockIndex{(\Lmat^{m})}{a,a-1} = m \blockIndex{\Lmat}{a,a-1}\) thus, holds true for \((\mathcal{A})\). For the condition \((\mathcal{B})\) this trivializes to the simple case of \(\blockIndex{\Lmat}{a,a-1} = \frac{1}{1!} \blockIndex{\Lmat}{a,a-1}\) which holds giving the base case to be true.

    The induction case will be shown as follows. Suppose that for some positive integer \(p-1\) that
    \begin{equation}\label{eq:inductionTriCond}
        \blockIndex{\Lmat}{a,b} = \frac{1}{(a-b)!} \prod_{i=0}^{(a-b)-1} \blockIndex{\Lmat}{a - i, a - i - 1} \quad \forall\, a - b < p,\, a > b.
    \end{equation}
    Looking at the block elementwise equation \eqref{eq:tricondsum} it can be seen that the only term that will not hold is the ordered partition set \(p\) thus, giving the following
    \begin{equation}
        \scalemath{0.9}{ m \, \blockIndex{\Lmat}{a,b} + \sum_{k=1}^{a - b} m^k \left[\sum_{n = k + \delta_{k,1}}^{a - b} \frac{s(n,k)}{n!} \left[ \sum_{\sigma \in \IntPartOrder{a - b}{n}} \prod_{i = 1}^{n} \blockIndex{\Lmat}{a - \overline {\sigma}(i-1), a - \overline {\sigma}(i)} \right]\right]}.
    \end{equation}
    Utilizing the induction hypothesis they all resolve to the same product of sub-diagonal blocks. Summed over the total number of orderings of the same partition giving the commutative case of \cref{def:IntegerPartitions} thus, simplifying the equation down to
    \begin{equation}
        \scalemath{0.85}{m\, \blockIndex{\Lmat}{a,b} + \sum_{k=1}^{a - b} m^k \left[\sum_{n=k + \delta_{k,1}}^{a - b} \frac{s(n,k)}{n!} \left[ \sum_{t_{1}^{q_1} t_{2}^{q_2} \dots t_{u}^{q_u} \in \IntPartUnorder{a - b}{n}}  {n \choose q_1, q_2, \dots q_u} \prod_{j=1}^{u}\frac{1}{(i!)^{q_j}} \right] \right] \prod_{i = 0}^{a - b - 1} \blockIndex{\Lmat}{a - i, a -  i - 1} }.
    \end{equation}
    Simplifying by applying \cref{lemma:combStirlingLemma} and since \(s(k,k) = 1\) and \(S(a-b, 1) = 1\) then
    \begin{equation}
        m\, \blockIndex{\Lmat}{a,b} + (m^{a-b} - m) \frac{1}{(a-b)!} \prod_{i = 0}^{a - b - 1} \blockIndex{\Lmat}{a - i, a -  i - 1}.
    \end{equation}
    Now since this is the blockwise entry of the power \(\blockIndex{(\Lmat^{m})}{a,b}\) and \((\mathcal{A})\) condition holds then shifting the term \(m \blockIndex{\Lmat}{a,b}\) over to LHS gives
    \begin{equation}
        (m^{a-b} - m) \blockIndex{\Lmat}{a,b} =  (m^{a-b} - m) \frac{1}{(a-b)!} \prod_{i = 0}^{a - b - 1} \blockIndex{\Lmat}{a - i, a -  i - 1}.
    \end{equation}
    Both sides of the equation have the same leading factor \(m^{a-b} - m\). This polynomial has roots of \(0,1\) and additionally \(-1\) when \(a-b\) is odd and greater than one. Given this holds for all \(m\) then it must hold for non-roots as well. Dividing both sides gives
    \begin{equation}
        \blockIndex{\Lmat}{a,b} =  \frac{1}{(a-b)!} \prod_{i = 0}^{a - b - 1} \blockIndex{\Lmat}{a - i, a -  i - 1} .
    \end{equation}
    Thus, the induction hypothesis holds for \(p\) as well giving \((\mathcal{B})\) is true. Given this holds by equality then \((\mathcal{A}) \equiv (\mathcal{B})\).

    Note that both sides have the same factors and no previous steps involved manipulating \(m\). Given the matrix power \((\mathcal{A})\) is known to hold for some value \(\delta \in \mathbb{R} \setminus \{1,0,-1\}\) then the conditions of both \((\mathcal{A})\) and \((\mathcal{B})\) must hold as well.

    \textbf{Case} \((\mathcal{B}) \equiv (\mathcal{C})\):\\
    Taking the block index form of the matrix log \cref{lemma:MatrixBlockLogPower} and removing the first term gives
    \begin{equation}\label{eq:logindpull}
        \blockIndex{\log(\Lmat)}{a, b} = \blockIndex{\Lmat}{a, b} + \sum_{n=2}^{a-b} \frac{s(n,1)}{n!}\left[ \sum_{\sigma \in \IntPartOrder{a - b}{n}} \prod_{i = 1}^{n} \blockIndex{\Lmat}{a - \overline {\sigma}(i-1), a - \overline {\sigma}(i)} \right].
    \end{equation}
    For the case of \(b = a - 1\) the summation cancels to zero thus, giving the sub-diagonal case of the log. For the next lower diagonal the log gives
    \begin{equation}
        \blockIndex{\log(\Lmat)}{a, a - 2} = \blockIndex{\Lmat}{a, a - 2} - \frac{1}{2} \blockIndex{\Lmat}{a, a - 1}\blockIndex{\Lmat}{a - 1, a - 2}.
    \end{equation}
    Thus, for this term to be zero then it must satisfy the case of \(\blockIndex{\Lmat}{a, a - 2} = \frac{1}{2} \blockIndex{\Lmat}{a, a - 1}\blockIndex{\Lmat}{a - 1, a - 2}\) which gives condition \((\mathcal{B})\) for \(p=2\). Applying the same induction step as \eqref{eq:inductionTriCond} up to \(p-1\) and applying \cref{lemma:combStirlingLemma} to \eqref{eq:logindpull} gives
    \begin{equation}
        \blockIndex{\log(\Lmat)}{a, b} = \blockIndex{\Lmat}{a, b} - \frac{1}{(a-b)!} \prod_{i=0}^{a-b-1} \blockIndex{\Lmat}{a - i, a - i - 1} \quad \forall\, a-b < p.
    \end{equation}
    Thus, for the term to be zero then it must satisfy the next higher condition for \((\mathcal{B})\). Therefore, by induction \((\mathcal{B})\) must be true for all cases to get the expected matrix structure giving \((\mathcal{B}) \equiv (\mathcal{C})\).
\end{proof}

It should be noted that such results encompass known properties for Pascal's matrix \cite{edelman2004pascal, aceto2001matrices}. Such results produced here are applied to the prune matrix of unordered trees of arbitrary colors. When looking at the general case of matrices, it is not always the case that one of the properties is satisfied. To handle this we introduce results for when the matrix is ``close'' to satisfying the properties.
\begin{theorem}\label{theory:TriangleMatrixClose}
    Given the result from \cref{theory:TriangleMatrix} suppose that instead condition \((\mathcal{B})\) is satisfied up to \(p \ge 1 \in \mathbb{Z}^+\) terms. Let the matrix \(\Identity + \CompLogBlock_1\) be the matrix where every entry satisfies the condition. This is given formally as
    \begin{equation}
        \begin{split}
            \blockIndex{\Lmat}{a,b} &= \frac{1}{(a-b)!} \prod_{i=0}^{a-b-1} \blockIndex{\Lmat}{a - i, a - i - 1} \quad \forall\, 0 < a-b \le p, \\
            \blockIndex{(\Identity + \CompLogBlock_1)}{a,b} &= \frac{1}{(a-b)!} \prod_{i=0}^{a-b-1} \blockIndex{\Lmat}{a - i, a - i - 1} \quad \forall\, 0 < a-b.
        \end{split}
    \end{equation}
    Condition \((\mathcal{A})\) is transformed to
    \begin{equation}
        \Lmat^m = \DiagMat_m \left(\Identity + \CompLogBlock_1 + \sum_{k = p}^{\infty} \frac{\CompLogBlock_{k + 1}}{m^{k}}\right) \DiagMat_m^{-1}
    \end{equation}
    where the matrices \(\CompLogBlock_{k+1}\) are nilpotent starting at the block sub-diagonal \(k+1\). The infinite sum is bounded by the number of blocks in the matrix. \(\CompLogBlock_{k+1}\) are characterized by interlaced powers of \(\log(\Lmat)\) and defined blockwise as
    \begin{equation}
        \blockIndex{(\CompLogBlock_{k+1})}{a,b} = \begin{cases}
            \blockIndex{\left(\frac{\log(\Lmat)^{(a-b) - k}}{((a-b) - k)!}\right)}{a,b} \qquad &\text{ if } a - b > k, a > b,\\
            0 \qquad &\text{ o.w.}
        \end{cases}
    \end{equation}
    Given the block \(\log\) power \cref{lemma:MatrixBlockLogPower} expansion then \(\forall\, p > j > 1\) gives \(\CompLogBlock_{j+1} = 0\). In the matrix \(\log(\Lmat)\) there will be the main block sub-diagonal and \(p-1\) repeated zero block sub-diagonals till the next non-zero sub-diagonal.
    Furthermore, at the limit, as long as \(p \ge 1\) then
    \begin{equation}
        \lim_{m \rightarrow \infty} \DiagMat_m^{-1}\, \Lmat^m\, \DiagMat_m = \Identity + \CompLogBlock_1.
    \end{equation}
\end{theorem}
\begin{proof}
    To start the proof the same expansion from the proof of \cref{theory:TriangleMatrix} will be used
    \begin{equation}
        \scalemath{0.9}{\blockIndex{(\Lmat^{m})}{a,b} = \sum_{k=1}^{a-b} m^k \left[\sum_{n=k}^{a-b} \frac{s(n,k)}{n!}\left[ \sum_{\sigma \in \IntPartOrder{a - b}{n}} \prod_{i = 1}^{n} \blockIndex{\Lmat}{a - \overline {\sigma}(i-1), a - \overline {\sigma}(i)} \right]\right]}.
    \end{equation}
    Given \cref{lemma:combStirlingLemma} the summed entry where \(k > a - b - p\) will cancel out into \(\frac{1}{p!}\delta_{k,p}\). Since the matrix block \((\CompLogBlock_1)_{a,b}\) contains the \(k=p\) term then
    \begin{equation}
        \blockIndex{(\Lmat^{m})}{a,b} = m^{a-b}\blockIndex{(\CompLogBlock_1)}{a,b} + \sum_{k = 1}^{a - b - p} m^k \left[\sum_{n=k}^{a - b} \frac{s(n,k)}{n!} \left[ \sum_{\sigma \in \IntPartOrder{a - b}{n}} \prod_{i = 1}^{n} \blockIndex{\Lmat}{a - \overline {\sigma}(i-1), a - \overline {\sigma}(i)} \right]\right].
    \end{equation}
    Applying the following sequence of summation identities \(\sum_{k=1}^{a-b-p}f(k) = \sum_{k=0}^{a-b-p-1}f(a-b-p-k) = \sum_{k=p}^{a-b-1}f(a-b-k)\) and expanding the scaling \(m^{a}\) and \(m^{-b}\) onto both sides of the equation gives
    \begin{equation}
        \scalemath{0.8}{\blockIndex{(\Lmat^{m})}{a,b} = m^{a}\left[ \blockIndex{(\CompLogBlock_1)}{a,b} + \sum_{k= p}^{a - b - 1} \frac{1}{m^k} \left[\sum_{n= a-b - k}^{a - b} \frac{s(n,a-b - k)}{n!} \left[ \sum_{\sigma \in \IntPartOrder{a - b}{n}} \prod_{i = 1}^{n} \blockIndex{\Lmat}{a - \overline {\sigma}(i-1), a - \overline {\sigma}(i)} \right]\right]\right]m^{-b}}.
    \end{equation}
    Given \cref{lemma:MatrixBlockLogPower} the inside terms match the block log term then
    \begin{equation}
        \blockIndex{(\Lmat^{m})}{a,b} = m^{a}\left[ \blockIndex{(\CompLogBlock_1)}{a,b} + \sum_{k= p}^{a - b - 1} \frac{1}{m^k} \blockIndex{\left(\frac{\log(\Lmat)^{a-b-k}}{(a-b-k)!}\right)}{a,b} \right]m^{-b}.
    \end{equation}
    Thus, giving the definition of \(\blockIndex{(\CompLogBlock_{k + 1})}{a,b}\). To show the limit inputting the equation derived into the hypothesis gives
    \begin{equation}
        \lim_{m \rightarrow \infty} \Identity + \CompLogBlock_1 + \sum_{i=p}^{\infty} \frac{1}{m^{i}}\, \CompLogBlock_{i+1} = \Identity + \CompLogBlock_1.
    \end{equation}
\end{proof}

This result provides insight into all block triangular matrices with diagonal identities such as the full space of composition, cutting, and replacement matrices formed by \BSeriesText or other binary operations. It is further noted the class of Heisenberg matrices \cite{hallLieGroupsLie2015a} are also covered (though the transposed results). Some examples of the \(\log\)s with specific diagonal block indices \(i \ge 0\) are given as 
\begin{equation}
    \scalemath{0.9}{
        \begin{split}
            \blockIndex{\log(\Lmat)}{1 + i,i} &= \blockIndex{(\CompLogBlock_{1})}{1 + i,i} = \blockIndex{\Lmat}{1 + i,i},\\
            \blockIndex{\log(\Lmat)}{2 + i,i} &= \blockIndex{(\CompLogBlock_{2})}{2 + i,i} = \blockIndex{\Lmat}{2 + i,i} - \frac{\blockIndex{\Lmat}{2 + i,1 + i}\blockIndex{\Lmat}{1 + i,i}}{2},\\
            \blockIndex{\log(\Lmat)}{3 + i,i} &= \blockIndex{(\CompLogBlock_{3})}{3 + i,i} = \blockIndex{\Lmat}{3 + i,i} - \frac{\blockIndex{\Lmat}{3 + i,2 + i}\blockIndex{\Lmat}{2 + i,i} + \blockIndex{\Lmat}{3 + i,1 + i}\blockIndex{\Lmat}{1 + i,i}}{2} \\
            &+ \frac{\blockIndex{\Lmat}{3 + i,2 + i}\blockIndex{\Lmat}{2 + i,1 + i}\blockIndex{\Lmat}{1 + i,i}}{3},\\
            \blockIndex{\log(\Lmat)}{4 + i,i} &= \blockIndex{(\CompLogBlock_{4})}{4 + i,i} = \blockIndex{\Lmat}{4 + i,i} - \frac{\blockIndex{\Lmat}{4 + i,1 + i}\blockIndex{\Lmat}{1 + i,i} + \blockIndex{\Lmat}{4 + i,2 + i}\blockIndex{\Lmat}{2 + i,i} + \blockIndex{\Lmat}{4 + i,3 + i}\blockIndex{\Lmat}{3 + i,i}}{2}\\
            &+ \frac{\blockIndex{\Lmat}{4 + i,2 + i}\blockIndex{\Lmat}{2 + i,1 + i}\blockIndex{\Lmat}{1 + i,i} + \blockIndex{\Lmat}{4 + i,3 + i}\blockIndex{\Lmat}{3 + i,1 + i}\blockIndex{\Lmat}{1 + i,i} + \blockIndex{\Lmat}{4 + i,3 + i}\blockIndex{\Lmat}{3 + i,2 + i}\blockIndex{\Lmat}{2 + i,i}}{3}\\
            & - \frac{\blockIndex{\Lmat}{4 + i,3 + i}\blockIndex{\Lmat}{3 + i,2 + i}\blockIndex{\Lmat}{2 + i,1 + i}\blockIndex{\Lmat}{1 + i,i}}{4}.\\
        \end{split}
    }
\end{equation}
For the second power these are
\begin{equation}
    \scalemath{0.85}{
        \begin{split}
            \blockIndex{\left(\frac{\log(\Lmat)^2}{2}\right)}{2 + i,i} &= \blockIndex{(\CompLogBlock_{1})}{2 + i,i} = \frac{\blockIndex{\Lmat}{2 + i,1 + i}\blockIndex{\Lmat}{1 + i,i}}{2},\\
            \blockIndex{\left(\frac{\log(\Lmat)^2}{2}\right)}{3 + i,i} &= \blockIndex{(\CompLogBlock_{2})}{3 + i,i} = \frac{\blockIndex{\Lmat}{3 + i,2 + i}\blockIndex{\Lmat}{2 + i,i} + \blockIndex{\Lmat}{3 + i,1 + i}\blockIndex{\Lmat}{1 + i,i}}{2} \\
            &- \frac{\blockIndex{\Lmat}{3 + i,2 + i}\blockIndex{\Lmat}{2 + i,1 + i}\blockIndex{\Lmat}{1 + i,i}}{2},\\
            \blockIndex{\left(\frac{\log(\Lmat)^2}{2}\right)}{4 + i,i} &= \blockIndex{(\CompLogBlock_{3})}{4 + i,i} = \frac{\blockIndex{\Lmat}{4 + i,1 + i}\blockIndex{\Lmat}{1 + i,i} + \blockIndex{\Lmat}{4 + i,2 + i}\blockIndex{\Lmat}{2 + i,i} + \blockIndex{\Lmat}{4 + i,3 + i}\blockIndex{\Lmat}{3 + i,i}}{2} \\
            &- \frac{\blockIndex{\Lmat}{4 + i,2 + i}\blockIndex{\Lmat}{2 + i, 1 + i}\blockIndex{\Lmat}{1 + i,i} + \blockIndex{\Lmat}{4 + i,3 + i}\blockIndex{\Lmat}{3 + i,1 + i}\blockIndex{\Lmat}{1 + i, i} + \blockIndex{\Lmat}{4 + i,3 + i}\blockIndex{\Lmat}{3 + i,2 + i}\blockIndex{\Lmat}{2 + i,i}}{2}\\
            & + \frac{11}{24}\blockIndex{\Lmat}{4 + i,3 + i}\blockIndex{\Lmat}{3 + i,2 + i}\blockIndex{\Lmat}{2 + i,1 + i}\blockIndex{\Lmat}{1 + i,i}.
        \end{split}
    }
\end{equation}

Lastly for a third power
\begin{equation}
    \scalemath{0.85}{
        \begin{split}
            \blockIndex{\left(\frac{\log(\Lmat)^3}{6}\right)}{3 + i,i} &= \blockIndex{(\CompLogBlock_{1})}{3 + i,i} = \frac{\blockIndex{\Lmat}{3 + i,2 + i}\blockIndex{\Lmat}{2 + i,1 + i}\blockIndex{\Lmat}{1 + i,i}}{6},\\
            \blockIndex{\left(\frac{\log(\Lmat)^3}{6}\right)}{4 + i,i} &= \blockIndex{(\CompLogBlock_{2})}{4 + i,i}\\
            &= \frac{\blockIndex{\Lmat}{4 + i,2 + i}\blockIndex{\Lmat}{2 + i,1 + i}\blockIndex{\Lmat}{1 + i,i} + \blockIndex{\Lmat}{4 + i,3 + i}\blockIndex{\Lmat}{3 + i,1 + i}\blockIndex{\Lmat}{1 + i,i} + \blockIndex{\Lmat}{4 + i,3 + i}\blockIndex{\Lmat}{3 + i,2 + i}\blockIndex{\Lmat}{2 + i,i}}{6}\\
            & - \frac{1}{4}\blockIndex{\Lmat}{4 + i,3 + i}\blockIndex{\Lmat}{3 + i,2 + i}\blockIndex{\Lmat}{2 + i,1 + i}\blockIndex{\Lmat}{1 + i,i}.
        \end{split}
    }
\end{equation}

\section{Proofs\label{section:Appendix}}

A selection of proofs on trees will utilize a technique laid out here. The idea comes from the descriptions of the sets of trees \cref{def:talltree}, \cref{def:bushtree}, and \cref{def:nthbushtree}. As stated each path from a root of a tree must eventually end up at a tree in set that is a lower \(n\)th tall bush tree. Thus, these proofs follow the steps.
\begin{itemize}
    \item Prove the base case of single node and/or empty tree. 
    \item Proving the tall trees case \cref{def:talltree}.
    \item Proving the tall bush tree case \cref{def:bushtree} utilizing tall tree results. This exposes how symmetry impacts the results and the forest case.
    \item Proving by induction that \cref{def:nthbushtree} holds for all \(n\) symmetry trees. Thus, by induction proving the general case of arbitrary trees via forest. This typically arises from using the same techniques used to prove \cref{def:bushtree}.
\end{itemize}
Note this route to proving these statements is not unique. Any other set of terminating trees could be used. Alternatively one could use the technique of successively grafting nodes onto a tree. There are some cases where a proof can be performed by going from a general tree to forest. However, the iterative tree set proof style given provides some additional insight in symmetry. It demonstrates the recursive operation over simple cases (tall trees) and directly how symmetry gets involved (tall bush trees).

\subsection{Tree Pruning \cref{theory:pruning}}\label{proof:pruning}
\begin{proof}
    For this proof we assume that \(\y^{\star}\) is a \BSeriesText. For the base case any \(\tree \prune \emptyset\) implies the \BSeriesText itself thus, satisfied. The converse case \(\emptyset \prune \tree\) must always be zero. Trees \(\tree \prune \tree\) imply only the initial condition in the expansion.

    \item 
    \paragraph{Tall Tree Case:} \hfill \\
    For a ladder tree \(\stree \in \LadderTrees{\NTTC}\) we proceed by induction climbing up the nodes. The base case holds given the empty tree has been proven. To start splitting the base elementary differential and applying the induction hypothesis gives
    \begin{equation}
        \BDiffOp(\stree)(\y^{\star}) = \BDiffOp(\tau[\xnode])(\y^{\star})\BDiffOp(\stree_1)(\y^{\star}) = \BDiffOp(\tau[\xnode])(\y^{\star})\left(\sum_{\tree \in \NTTC} a_0(\tree \prune \stree_1) \frac{\BDiffOp(\tree)(\y_\initstate)}{\sigma(\tree)}\right) + \frac{\BDiffOp(\stree)(\y_\initstate)}{\sigma(\stree)}
    \end{equation}
    where \(a_0(\emptyset) = 0\) and otherwise \(a_0(\tree) = a(\tree)\).
    Next Taylor expanding around the initial condition and inserting the base \BSeriesText gives
    \begin{equation}
        = \sum_{p=0}^{\infty} \frac{\BDiffOp(\tau[\xnode^{p+1}])(\y_\initstate)}{p!} \left(\sum_{\tree \in \NTTC} a_0(\tree \prune \emptyset) \frac{\BDiffOp(\tree)(\y_\initstate)}{\sigma(\tree)}\right)^{p}\left(\sum_{\tree \in \NTTC} a_0(\tree \prune \stree_1) \frac{\BDiffOp(\tree)(\y_\initstate)}{\sigma(\tree)}\right) + \frac{\BDiffOp(\stree)(\y_\initstate)}{\sigma(\stree)}.
    \end{equation}
    Applying \cref{lemma:permanentTrees} to the product terms gives
    \begin{equation}
        = \sum_{p=0}^{\infty} \frac{\BDiffOp(\tau[\xnode^{p+1}])(\y_\initstate)}{p!} \sum_{\fforest = \prod_{j=1}^{n} \tree_{j}^{q_j} \in \ForestSpace{\NTTC}_{p+1}} \frac{1}{\prod_{j=1}^{n} q_j !} \permanent\left(\fforest \, \outergen{a_0(\cdot_1 \prune \cdot_2) \frac{\BDiffOp(\cdot_1)(\y_\initstate)}{\sigma(\cdot_1)}} \, \emptyset^{p} \stree_1 \right) + \frac{\BDiffOp(\stree)(\y_\initstate)}{\sigma(\stree)}.
    \end{equation}
    Pulling the elementary differential and symmetry out via \cref{property:permanentconstant} and applying the forest definition of pruning gives
    \begin{equation}
        = \sum_{p=0}^{\infty} \sum_{\fforest = \prod_{j=1}^{n} \tree_{j}^{q_j} \in \ForestSpace{\NTTC}_{p+1}} \frac{\BDiffOp(\tau[\xnode^{p+1}])(\y_\initstate) \left(\prod_{j=1}^{n} \BDiffOp(\tree_j)(\y_\initstate)^{q_j}\right)}{\prod_{j=1}^{n} q_j ! \sigma(\tree_{j})^{k_j}} a_0\left(\fforest \prune  \emptyset^{p} \stree_1 \right) + \frac{\BDiffOp(\stree)(\y_\initstate)}{\sigma(\stree)}.
    \end{equation}
    Lastly, since the pruning produces zeroed results for trees which do not have a base structure of \(\stree\) then this expression is reduced down to the tree form as
    \begin{equation}
        = \sum_{\tree \in \NTTC} a(\tree \prune \stree) \frac{\BDiffOp(\tree)(\y_\initstate)}{\sigma(\tree)}.
    \end{equation}
    Thus, by induction the ladder tree case holds. We note that during the Taylor expansion step we use the catch-all stump \(\xnode\) for the derivatives. In the case of more colors where the roots matter this is still valid. This is due to \cref{lemma:permanentmultiset} producing the leading symmetry factors. If the general \(\emptyset\) is used instead and the individual labels zero if the proper root is not provided then by the symmetry set of the forest \(S_{\gforest \emptyset^{p}}\) will divide out leading symmetries leading to the colored empty forest \(S_{\gforest \prod_{c \in \Tcolors} \emptyset^{p_c}}\). Such results will be utilized for further steps in the pruning proof.  

    \item 
    \paragraph{Tall Bush Case} \hfill\\
    Next for a tall bush tree \(\stree = \tau[\prod_{j=1}^{m} \stree_{j}^{k_j}] = \tau[\gforest] \in \TallBushTrees{\NTTC}\) we proceed as follows. First splitting the base off and applying the induction hypothesis gives
    \begin{equation}
        \frac{\BDiffOp(\stree)(\y^{\star})}{\sigma(\stree)} = \frac{\BDiffOp(\tau[\xnode^{\bar{m}}])(\y^{\star})}{\prod_{j=1}^{m} k_j!}\prod_{j=1}^{m}\left(\frac{\BDiffOp(\stree_j)(\y^{\star})}{\sigma(\stree_j)}\right)^{k_j} = \frac{\BDiffOp(\tau[\xnode^{\bar{m}}])(\y^{\star})}{\prod_{j=1}^{m} k_j!}\prod_{j=1}^{m}\left(\sum_{\tree \in \NTTC} a_0(\tree \prune \stree_j) \frac{\BDiffOp(\tree)(\y_\initstate)}{\sigma(\tree)}\right)^{k_j} + \frac{\BDiffOp(\stree)(\y_\initstate)}{\sigma(\stree)}.
    \end{equation}
    Next Taylor expanding the base gives
    \begin{equation}
        = \sum_{p=0}^{\infty}\frac{\BDiffOp(\tau[\xnode^{p + \bar{m}}])(\y_\initstate)}{\left(\prod_{j=1}^{m} k_j!\right)p!}\prod_{j=1}^{m}\left(\sum_{\tree \in \NTTC} a_0(\tree \prune \stree_j) \frac{\BDiffOp(\tree)(\y_\initstate)}{\sigma(\tree)}\right)^{k_j}\left(\sum_{\tree \in \NTTC} a_0(\tree \prune \emptyset) \frac{\BDiffOp(\tree)(\y_\initstate)}{\sigma(\tree)}\right)^{p} + \frac{\BDiffOp(\stree)(\y_\initstate)}{\sigma(\stree)}.
    \end{equation}
    Applying \cref{lemma:permanentTrees} to the product terms gives
    \begin{equation}
        = \sum_{p=0}^{\infty}\frac{\BDiffOp(\tau[\xnode^{p + \bar{m}}])(\y_\initstate)}{\left(\prod_{j=1}^{m} k_j!\right)p!} \sum_{\fforest = \prod_{j=1}^{n} \tree_{j}^{q_j} \in \ForestSpace{\NTTC}_{p + \bar{m}}} \frac{1}{\prod_{j=1}^{n} q_j !} \permanent\left(\fforest \, \outergen{a_0(\cdot_1 \prune \cdot_2) \frac{\BDiffOp(\cdot_1)(\y_\initstate)}{\sigma(\cdot_1)}} \, \emptyset^{p} \gforest \right) + \frac{\BDiffOp(\stree)(\y_\initstate)}{\sigma(\stree)}.
    \end{equation}
    Pulling the elementary differential and symmetry out via \cref{property:permanentconstant} and applying the forest definition of pruning gives
    \begin{equation}
        = \sum_{p=0}^{\infty}\sum_{\fforest = \prod_{j=1}^{n} \tree_{j}^{q_j} \in \ForestSpace{\NTTC}_{p + \bar{m}}} \frac{\BDiffOp(\tau[\xnode^{p + \bar{m}}])(\y_\initstate) \left(\prod_{j=1}^{n} \BDiffOp(\tree_j)(\y_\initstate)^{q_j}\right)}{\prod_{j=1}^{n} q_j ! \sigma(\tree_{j})^{k_j}} a_0(\fforest \prune \emptyset^{p} \gforest) + \frac{\BDiffOp(\stree)(\y_\initstate)}{\sigma(\stree)}.
    \end{equation}
    Finally, since pruning gives zeros when \(\stree \nless \tree\) then this reduces to the tree case as
    \begin{equation}
        = \sum_{\tree \in \NTTC} a(\tree \prune \stree) \frac{\BDiffOp(\tree)(\y_\initstate)}{\sigma(\tree)}.
    \end{equation}
    Thus, pruning holds for the tall bush case.

    \item 
    \paragraph{General Case:} \hfill \\
    For the general case it is proceeded over inductively by \(\stree = \tau[\forestTQs] \in \TallBushTreesN{\NTTC}{i}\) by assuming \(\TallBushTreesN{\NTTC}{i-1}\) holds true. All steps applied in the proof of the tall bush case are equally applied in this case thus, giving the same result. Therefore, the result holds.
\end{proof}

\subsection{Pruning Properties \cref{theory:pruneProperties}} \label{proof:pruneProperties}
\begin{proof}
    For \cref{theory:TriangleMatrix} all properties can be validated when proven for a power not in \(\{-1, 0, 1\}\). Given that \cref{theory:analyticalComposition} holds then it suffices for powers of \(2\) thus proving the statement.
\end{proof}

\subsection{Marked Pruning \cref{theory:pruningMarked}}\label{proof:pruningMarked}

For this proof we provide a lemma on pair right grafting over tree label brackets first.
\begin{lemma}[Pair Right Grafting Recursive Tree Brackets]\label{lemma:PRGLabel}
    Given a set of trees \(\mathrm{T}\) with a tree \(\stree \in \NTTC\) then the recursive function
    \begin{equation}
        \begin{split}
            \barwedge_{\mathrm{T}}(\stree) &= \BDiffOp(\tau[\pitchfork(\gforest)])(\y_\initstate) (\tau \ProductMarked \barwedge_{\mathrm{T}}(\gforest)),\\
            \barwedge_{\mathrm{T}}(\forestTQs) &= \prod_{i=1}^{m}(\barwedge_{\mathrm{T}}(\stree_{i}))^{k_i},\\
            \barwedge_{\mathrm{T}}(\xnode) &= \BDiffOp(\xnode)(\y_\initstate) \sum_{\ktree \in \mathrm{T}} \labelbrack{\ktree}{\xnode}\frac{\BDiffOp(\ktree)(\y_\initstate)}{\sigma(\ktree)},\\
            \barwedge_{\mathrm{T}}(\tau) &= \BDiffOp(\tau)(\y_\initstate) \labelbrack{1}{\tau},
        \end{split}
    \end{equation}
    is equivalent to 
    \begin{equation}
        \begin{split}
            \barwedge_{\mathrm{T}}(\stree) &= \stumpCountFactorial(\stree) \sum_{\hforest \in \ForestSpace{\mathrm{T}}_{\stumpCount(\stree)}}  \labelbrack{\hforest}{\stree} \BDiffOp(\stree)(\y_\initstate) \symprod \frac{\BDiffOp(\hforest)(\y_\initstate)}{\sigma(\hforest)},\\
            \barwedge_{\mathrm{T}}(\forestTQs) &= \prod_{i=1}^{m}\left( \stumpCountFactorial(\stree_{i}) \sum_{\hforest \in \ForestSpace{\mathrm{T}}_{\stumpCount(\stree_i)}} \labelbrack{\hforest}{\stree_i} \BDiffOp(\stree_{i})(\y_\initstate) \symprod \frac{\BDiffOp(\hforest)(\y_\initstate)}{\sigma(\hforest)}\right)^{k_i}.
        \end{split}
    \end{equation}
\end{lemma}
\begin{proof}
    The function \(\barwedge_{\mathrm{T}}\) can be treated as a function where the inputs are the stumps \(\xnode\). Given the inputs work over a sum of all possible trees then by linearity of \(\tau \ProductMarked \) and the application of \(\symprod\) produces all permutations of inputs then the size permutations of forest inputs are considered. Combining inputs to a forest \(\hforest\) gives the permutation counts as \(\frac{\prod_{i=1}^{m} \stumpCountFactorial(\stree_{i})}{\sigma(\hforest)}\) for the entries in the sum. Given the property of \(\tau \ProductMarked \) to collapse right entries then this takes the individual trees into forest combinations with the right label entry constructed as a tree. For the sum over all possible forest gives the final result.
\end{proof}

Now for the proof of \cref{theory:pruningMarked}.

\begin{proof}
    Since \cref{theory:pruning} states that the recursion returns the trimmings these operations replicate that exactly. The only modification is the introduction of additional output via the tree brackets. Given that the base case marks the tree with a stump then the recursive operator simply reconstructs the original tree with additional stumps. Since pair right grafting and tree brackets respect the original operation and pair right grafting only grafts back the nodes on the right entry after pruning, then the statement holds.

    For the next case of the summed result. The base case is given by the nodal definition of marked pruning. The empty set case is the sum over all trees with stump when the node is not empty giving
    \begin{equation}
        \sum_{\tree \in \NTTC} a(\tree \pruneMark \emptyset) \frac{\BDiffOp(\tree)(\y_\initstate)}{\sigma(\tree)} = \sum_{\ktree \in \NTTC} a(\labelbrack{\ktree}{\xnode_{\treeRoot(\ktree)}}) \frac{\BDiffOp(\ktree)(\y_\initstate)}{\sigma(\ktree)} + a(\labelbrack{1}{\emptyset})\BDiffOp(\emptyset)(\y_\initstate).
    \end{equation}
    \item 
    \paragraph{Ladder Case:} \hfill \\
    For the case we proceed via induction over the nodes of \(\stree \in \LadderTrees{\NTTC}\). Since the empty set case was proven then the base case is established. Starting with the established form
    \begin{equation}
        \sum_{\tree \in \NTTC} a(\tree \pruneMark \stree) \frac{\BDiffOp(\tree)(\y_\initstate)}{\sigma(\tree)}
    \end{equation}
    and expanding over the forest definition of \(\pruneMark\) for the case the roots equal gives
    \begin{equation}
        = \sum_{\tree = \tau[\prod_{i=1}^{n} \tree_{i}^{q_i}] = \tau[\fforest] \in \NTTC} a\left(\tau \ProductMarked \left(\frac{1}{p!} \permanent(\fforest \outergen{\cdot_1 \pruneMark \cdot_2} \stree_{1} \emptyset^{p})\right)\right)\frac{\BDiffOp(\tree)(\y_\initstate)}{\sigma(\tree)} + a(\labelbrack{\tree}{\xnode}) \frac{\BDiffOp(\tree)(\y_\initstate)}{\sigma(\tree)}
    \end{equation}
    where \(\stree_{1}\) is the child of \(\stree\). Moving the sum in the pair grafting and pushing the elementary differential and symmetry terms in the permanent gives
    \begin{equation}
        = a\left(\tau \ProductMarked \left(\sum_{p=0}^{\infty}\frac{\BDiffOp\left(\tau[ \xnode^{1+p}]\right)(\y_\initstate)}{p!} \sum_{\fforest \in \ForestSpace{\NTTC}_{1+p}}  \frac{1}{\prod_{j=1}^{n} q_j!} \permanent(\fforest \outergen{(\cdot_1 \pruneMark \cdot_2) \frac{\BDiffOp(\cdot_1)(\y_\initstate)}{\sigma(\cdot_1)}} \stree_{1} \emptyset^{p}) \right)\right) + a(\labelbrack{\tree}{\xnode}) \frac{\BDiffOp(\tree)(\y_\initstate)}{\sigma(\tree)}
    \end{equation}
    where the sum must go over all possible sizes of forest.
    Applying \cref{lemma:permanentTrees} gives 
    \begin{equation}
        = a\left(\tau \ProductMarked \left(\sum_{p=0}^{\infty}\frac{\BDiffOp(\tau[\xnode^{1+p}])(\y_\initstate)}{p!} \left(\sum_{\tree \in \NTTC} (\tree \pruneMark \stree_i) \frac{\BDiffOp(\tree)(\y_\initstate)}{\sigma(\tree)}\right)\left(\sum_{\tree \in \NTTC} (\tree \pruneMark \emptyset) \frac{\BDiffOp(\tree)(\y_\initstate)}{\sigma(\tree)}\right)^{p} \right)\right) + a(\labelbrack{\tree}{\xnode}) \frac{\BDiffOp(\tree)(\y_\initstate)}{\sigma(\tree)}.
    \end{equation}
    Next inserting the induction hypothesis gives
    \begin{equation}
        \begin{split}
            = & a\left(\tau \ProductMarked \left(\sum_{p=0}^{\infty}\frac{\BDiffOp\left(\tau[\xnode^{1+p}]\right)(\y_\initstate)}{p!} \left(\sum_{\stree^* \in \stree\lceil \TStumpcolors \rceil}\stumpCountFactorial(\stree^{*}) \sum_{\hforest \in \ForestSpace{\NTTC}_{\stumpCount(\stree^{*})}} \labelbrack{\hforest}{\stree^{*}} \frac{\BDiffOp(\stree^{*})(\y_\initstate)}{\sigma(\stree^{*})} \symprod \frac{\BDiffOp(\hforest)(\y_\initstate)}{\sigma(\hforest)}\right)\right. \right. \\
            \qquad & \left. \left. \left(\sum_{\ktree \in \NTTC} \labelbrack{\ktree}{\xnode} \BDiffOp(\xnode)(\y_\initstate)\frac{\BDiffOp(\ktree)(\y_\initstate)}{\sigma(\ktree)}\right)^{p}\right)\right) + a(\labelbrack{\tree}{\xnode}) \frac{\BDiffOp(\tree)(\y_\initstate)}{\sigma(\tree)}.
        \end{split}
    \end{equation}
    Pulling the stumped sum and the leading symmetries out and moving the marked prune in gives
    \begin{equation}
        \begin{split}
            = & a\left(\left(\sum_{p=0}^{\infty} \sum_{\stree^* \in \stree\lceil \TStumpcolors \rceil} \frac{\BDiffOp\left(\tau[\xnode^{1+p}]\right)(\y_\initstate)}{p!\sigma(\stree^{*})} \tau \ProductMarked \left(\stumpCountFactorial(\stree^{*}) \sum_{\hforest \in \ForestSpace{\NTTC}_{\stumpCount(\stree^{*})}} \labelbrack{\hforest}{\stree^{*}} \BDiffOp(\stree^{*})(\y_\initstate) \symprod \frac{\BDiffOp(\hforest)(\y_\initstate)}{\sigma(\hforest)}\right)\right.\right.\\
            \qquad & \left.\left.\left(\sum_{\ktree \in \NTTC} \labelbrack{\ktree}{\xnode} \BDiffOp(\xnode)(\y_\initstate)\frac{\BDiffOp(\ktree)(\y_\initstate)}{\sigma(\ktree)}\right)^{p}\right)\right) + a(\labelbrack{\tree}{\xnode}) \frac{\BDiffOp(\tree)(\y_\initstate)}{\sigma(\tree)}.
        \end{split}
    \end{equation}
    Applying the results of \cref{lemma:PRGLabel} to the product of terms and simplifying gives
    \begin{equation}
            = a\left(\sum_{p=0}^{\infty} \sum_{\stree^* \in \stree_1\lceil \TStumpcolors \rceil}\stumpCountFactorial(\tau[\stree^{*} \xnode^{p}]) \sum_{\hforest \in \ForestSpace{\NTTC}_{\stumpCount(\stree^{*}) + p}} \labelbrack{\hforest}{\tau[\stree^{*} \xnode^{p}]}  \frac{\BDiffOp(\tau[\stree^{*} \xnode^{p}])(\y_\initstate)}{\sigma(\tau[\stree^{*} \xnode^{p}])}  \symprod \frac{\BDiffOp(\hforest)(\y_\initstate)}{\sigma(\hforest)} \right) + a(\labelbrack{\tree}{\xnode}) \frac{\BDiffOp(\tree)(\y_\initstate)}{\sigma(\tree)}.
    \end{equation}
    The infinite sum formed over the extra stump inserted is simplified to stumps over the full tree giving
    \begin{equation}
        = \sum_{\ktree \in \stree\lceil \TStumpcolors \rceil}\stumpCountFactorial(\ktree) \sum_{\hforest \in \ForestSpace{\NTTC}_{\stumpCount(\ktree)}} a(\labelbrack{\hforest}{\ktree})  \frac{\BDiffOp(\ktree)(\y_\initstate)}{\sigma(\ktree)}  \symprod \frac{\BDiffOp(\hforest)(\y_\initstate)}{\sigma(\hforest)}.
    \end{equation}

    \item 
    \paragraph{Tall Bush Case:} \hfill \\
    For tall bush case we proceed by increasing the amount of children at the base node with \(\stree = \tau[\forestTQs] = \tau[\gforest] \in \TallBushTrees{\NTTC}\). Given the ladder case holds then the base case holds. Starting from
    \begin{equation}
        \sum_{\tree \in \NTTC} a(\tree \pruneMark \stree) \frac{\BDiffOp(\tree)(\y_\initstate)}{\sigma(\tree)}
    \end{equation}
    and expanding over the forest definition of \(\pruneMark\) for when the roots are the same gives
    \begin{equation}
        = \sum_{\tree = \tau[\prod_{i=1}^{n} \tree_{i}^{q_i}] = \tau[\fforest] \in \NTTC} a\left(\tau \ProductMarked \left(\frac{1}{\left(\prod_{j=1}^{m} k_j !\right) p!} \permanent\left(\fforest \outergen{\cdot_1 \pruneMark \cdot_2} \gforest \emptyset^{p}\right)\right)\right)\frac{\BDiffOp(\tree)(\y_\initstate)}{\sigma(\tree)} + a(\labelbrack{\tree}{\xnode}) \frac{\BDiffOp(\tree)(\y_\initstate)}{\sigma(\tree)}.
    \end{equation}
    Moving the sum in the pair grafting and pushing the elementary differential and symmetry terms in the permanent gives
    \begin{equation}
        = a\left(\tau \ProductMarked \left(\sum_{p=0}^{\infty} \frac{\BDiffOp(\tau[\stumpForestConvert(\gforest)\xnode^{p}])(\y_\initstate)}{\left(\prod_{j=1}^{m} k_j !\right) p!} \sum_{\fforest \in \ForestSpace{\NTTC}_{\bar{k} + p}} \frac{1}{\prod_{j=1}^{n} q_j!} \permanent\left(\fforest \outergen{(\cdot_1 \pruneMark \cdot_2) \frac{\BDiffOp(\cdot_1)(\y_\initstate)}{\sigma(\cdot_1)}} \gforest\, \emptyset^{p}\right) \right)\right) + a(\labelbrack{\tree}{\xnode}) \frac{\BDiffOp(\tree)(\y_\initstate)}{\sigma(\tree)}.
    \end{equation}
    Applying \cref{lemma:permanentTrees} gives 
    \begin{equation}
        \begin{split}
            &= a\left(\tau \ProductMarked \left(\sum_{p=0}^{\infty}\frac{\BDiffOp(\tau[\stumpForestConvert(\gforest)\xnode^{p}])(\y_\initstate)}{\left(\prod_{j=1}^{m} k_j !\right) p!} \prod_{j=1}^{m}\left(\sum_{\tree \in \NTTC} (\tree \pruneMark \stree_i) \frac{\BDiffOp(\tree)(\y_\initstate)}{\sigma(\tree)}\right)^{k_i}\left(\sum_{\tree \in \NTTC} (\tree \pruneMark \emptyset) \frac{\BDiffOp(\tree)(\y_\initstate)}{\sigma(\tree)}\right)^{p} \right)\right)\\
            &+ a(\labelbrack{\tree}{\xnode}) \frac{\BDiffOp(\tree)(\y_\initstate)}{\sigma(\tree)}.
        \end{split}
    \end{equation}
    Next applying the previous result since all subtrees are ladders gives
    \begin{equation}
        \begin{split}
            &= a\left(\tau \ProductMarked \left(\sum_{p=0}^{\infty}\frac{\BDiffOp(\tau[\stumpForestConvert(\gforest)\xnode^{p}])(\y_\initstate)}{\left(\prod_{j=1}^{m} k_j !\right) p!}\prod_{j=1}^{m}\left(\sum_{\ktree \in \stree_{j}\lceil \TStumpcolors \rceil}\stumpCountFactorial(\ktree) \sum_{\hforest \in \ForestSpace{\NTTC}_{\stumpCount(\ktree)}} a(\labelbrack{\hforest}{\ktree})  \frac{\BDiffOp(\ktree)(\y_\initstate)}{\sigma(\ktree)}  \symprod \frac{\BDiffOp(\hforest)(\y_\initstate)}{\sigma(\hforest)}\right)^{k_i}\right.\right.\\
            \qquad & \left.\left.\left(\sum_{\ktree \in \NTTC} \labelbrack{\ktree}{\xnode} \BDiffOp(\xnode)(\y_\initstate)\frac{\BDiffOp(\ktree)(\y_\initstate)}{\sigma(\ktree)}\right)^{p} \right)\right) + a(\labelbrack{\tree}{\xnode}) \frac{\BDiffOp(\tree)(\y_\initstate)}{\sigma(\tree)}.
        \end{split}
    \end{equation}
    Next forming a multinomial expansion over the possible stumped characters \(\stree_j\lceil \TStumpcolors \rceil\) gives the forest sum as 
    \begin{equation}
        \begin{split}
            = & a\left(\tau \ProductMarked \left(\sum_{p=0}^{\infty}\frac{\BDiffOp(\tau[\stumpForestConvert(\gforest)\xnode^{p}])(\y_\initstate)}{\left(\prod_{j=1}^{m} k_j !\right) p!}\prod_{i=1}^{m} \sum_{\prod_{j=1}^{\bbar_{i}}\streebar_{i,j}^{\bbar_{i,j}} \in  \ForestSpace{\stree_{i}\lceil \TStumpcolors \rceil}_{k_i}}\frac{k_i!}{\prod_{j=1}^{\bbar_{i}} \bbar_{i,j}!}\prod_{j=1}^{\bbar_{i}}\right.\right.\\
            \qquad & \left.\left. \left(\stumpCountFactorial(\streebar_{i,j})\sum_{\hforest \in \ForestSpace{\NTTC}_{\stumpCount(\streebar_{i,j})}} a(\labelbrack{\hforest}{\streebar_{i,j}})  \frac{\BDiffOp(\streebar_{i,j})(\y_\initstate)}{\sigma(\streebar_{i,j})}  \symprod \frac{\BDiffOp(\hforest)(\y_\initstate)}{\sigma(\hforest)}\right)^{\bbar_{i,j}}\left(\sum_{\ktree \in \NTTC} \labelbrack{\ktree}{\xnode} \BDiffOp(\xnode)(\y_\initstate)\frac{\BDiffOp(\ktree)(\y_\initstate)}{\sigma(\ktree)}\right)^{p} \right)\right)\\
            & + a(\labelbrack{\tree}{\xnode}) \frac{\BDiffOp(\tree)(\y_\initstate)}{\sigma(\tree)}.
        \end{split}
    \end{equation}
    Given all the character sums are unique then expanding the sums out and collecting the factorial terms gives
    \begin{equation}
        \begin{split}
            = & a\left(\sum_{p=0}^{\infty}\sum_{\tilde{\gforest} = \prod_{i=1}^{m}\prod_{j=1}^{\bbar_{i}}\streebar_{i,j}^{\bbar_{i,j}} \in \ForestSpace{\gforest\lceil \TStumpcolors \rceil}_{\bar{k}}} \frac{\BDiffOp(\tau[\stumpForestConvert(\gforest)\xnode^{p}])(\y_\initstate)}{\left(\prod_{i=1}^{m}\prod_{j=1}^{\bbar_{i}} \bbar_{i,j}!\sigma(\streebar_{i,j})^{\bbar_{i,j}}\right) p!} \right.\\
            \qquad &  \tau \ProductMarked \left( \prod_{i=1}^{m}\prod_{j=1}^{\bbar_{i}}\left(\stumpCountFactorial(\streebar_{i,j})\sum_{\hforest \in \ForestSpace{\NTTC}_{\stumpCount(\streebar_{i,j})}} a(\labelbrack{\hforest}{\streebar_{i,j}})  \BDiffOp(\streebar_{i,j})(\y_\initstate) \symprod \frac{\BDiffOp(\hforest)(\y_\initstate)}{\sigma(\hforest)}\right)^{\bbar_{i,j}}\right.\\
            \qquad & \left.\left.\left(\sum_{\ktree \in \NTTC} \labelbrack{\ktree}{\xnode} \BDiffOp(\xnode)(\y_\initstate)\frac{\BDiffOp(\ktree)(\y_\initstate)}{\sigma(\ktree)}\right)^{p} \right)\right) + a(\labelbrack{\tree}{\xnode}) \frac{\BDiffOp(\tree)(\y_\initstate)}{\sigma(\tree)}.
        \end{split}
    \end{equation}

    Applying \cref{lemma:PRGLabel} to the product terms and simplifying gives
    \begin{equation}
        \begin{split}
            =& a\left(\sum_{p=0}^{\infty} \sum_{\prod_{i=1}^{m}\prod_{j=1}^{\bbar_{i}}\streebar_{i,j}^{\bbar_{i,j}} \in \ForestSpace{\gforest\lceil \TStumpcolors \rceil}_{\bar{k}}} \stumpCountFactorial\left(\tau\left[\prod_{i=1}^{m}\prod_{j=1}^{\bbar_{i}}\streebar_{i,j}^{\bbar_{i,j}}\xnode^{p}\right]\right) \sum_{\hforest \in \ForestSpace{\NTTC}_{\stumpCount(\tau\left[\prod_{i=1}^{m}\prod_{j=1}^{\bbar_{i}}\streebar_{i,j}^{\bbar_{i,j}}\xnode^{p}\right])}} \right.\\
            \qquad &\left.\labelbrack{\hforest}{\tau\left[\prod_{i=1}^{m}\prod_{j=1}^{\bbar_{i}}\streebar_{i,j}^{\bbar_{i,j}}\xnode^{p}\right]} \frac{\BDiffOp(\tau[\prod_{i=1}^{m}\prod_{j=1}^{\bbar_{i}}\streebar_{i,j}^{\bbar_{i,j}}\xnode^{p}])(\y_\initstate)}{\sigma(\tau[\prod_{i=1}^{m}\prod_{j=1}^{\bbar_{i}}\streebar_{i,j}^{\bbar_{i,j}}\xnode^{p}])} \symprod \frac{\BDiffOp(\hforest)(\y_\initstate)}{\sigma(\hforest)} \right) + a(\labelbrack{\tree}{\xnode}) \frac{\BDiffOp(\tree)(\y_\initstate)}{\sigma(\tree)}.
        \end{split}
    \end{equation}
    The infinite sum over the stumps is enveloped into the forest giving the set of all possible graft of stumps on the tree giving
    \begin{equation}
        = \sum_{\ktree \in \lceil \TStumpcolors \rceil} \stumpCountFactorial(\ktree) \sum_{\hforest \in \ForestSpace{\NTTC}_{\stumpCount(\ktree)}} \frac{\BDiffOp(\ktree)(\y_\initstate)}{\sigma(\ktree)} \symprod \frac{\BDiffOp(\hforest)(\y_\initstate)}{\sigma(\hforest)}.
    \end{equation}
    Thus, the results hold for tall bush trees.
    \item 
    \paragraph{General Case:} \hfill \\
    For the general case induction is applied to reach the tall bush case. Applying the exact same proof as was done over the tall bush trees for \(\stree = \tau[\forestTQs] \in \TallBushTreesN{\NTTC}{i}\) by assuming \(\TallBushTreesN{\NTTC}{i-1}\) holds gives the same results. Thus, the general case holds.
\end{proof}

\subsection{Replacement Operation \cref{theory:treereplacement}}\label{proof:treereplacement}
\begin{proof}
    To start if any input is given as \(\emptyset\) then only \(\emptyset\) is possible thus, the base conditions on replacement.
    In general if any tree for left input contains a node in \(\mathcal{U}\) then this must always return \(0\). Since there is no replacement to end up with the final elementary differential.
    Next for the case of nodes. If the node \(\psi \notin \mathcal{U}\) then the only possible result is itself since there is no replacement. In the case of \(\psi \in \mathcal{U}\) then an expansion gives the result
    \begin{equation}
        \BDiffOp(\psi)(\y_\initstate) = \sum_{\utree \in \NTTi{\mathcal{A}} \setminus \{\emptyset\}}  a(\labelbrack{\utree}{\psi}) \frac{\BDiffOp(\utree)(\y_\initstate)}{\sigma(\utree)}
    \end{equation}
    by definition. Thus, for replacement then \(\utree \replaceNode{\mathcal{U}} \psi = \labelbrack{\utree}{\psi}\) must hold for all trees without nodes in the replacement set.
    \item 
    \paragraph{Ladder Case:} \hfill \\
    Given the ladder tree \(\psi[\laddertree]\) where \(\laddertree \ne \emptyset\) then using the induction hypothesis gives
    \begin{equation}
        \BDiffOp(\psi[\laddertree])(\y_\initstate) = \BDiffOp(\psi[\xnode])(\y_\initstate)\BDiffOp(\laddertree)(\y_\initstate) = \BDiffOp(\psi[\xnode])(\y_\initstate)\left(\sum_{\ktree \in \NTTi{\mathcal{A}}} a(\ktree \replaceNode{\mathcal{U}} \laddertree) \frac{\BDiffOp(\ktree)(\y_\initstate)}{\sigma(\ktree)}\right).
    \end{equation}
    If \(\psi \notin \mathcal{U}\) then this resolves to
    \begin{equation}
         = \sum_{\ktree \in \NTTi{\mathcal{A}}} a(\ktree \replaceNode{\mathcal{U}} \psi[\laddertree]) \frac{\BDiffOp(\ktree)(\y_\initstate)}{\sigma(\ktree)}.
    \end{equation}
    In the case that \(\psi \in \mathcal{U}\) then expanding the base node into its series form of a single stump gives
    \begin{equation}
        = \left(\sum_{\stree \in \NTTi{\mathcal{A}}\lceil\xnode\rceil} a(\labelbrack{\stree}{\psi}) \frac{\BDiffOp(\stree)(\y_\initstate)}{\sigma(\stree)}\right)\left(\sum_{\ktree \in \NTTi{\mathcal{A}}} a(\ktree \replaceNode{\mathcal{U}} \laddertree) \frac{\BDiffOp(\ktree)(\y_\initstate)}{\sigma(\ktree)}\right).
    \end{equation}
    Expanding the product into double sum and applying the tree bracket replacement property gives
    \begin{equation}
        = \sum_{\stree \in \NTTi{\mathcal{A}}\lceil\xnode\rceil} \sum_{\ktree \in \NTTi{\mathcal{A}}} a(\labelbrack{\ktree}{\stree} \replaceNode{\mathcal{U}} \labelbrack{\laddertree}{\psi}) \frac{\BDiffOp(\stree)(\y_\initstate)}{\sigma(\stree)}\frac{\BDiffOp(\ktree)(\y_\initstate)}{\sigma(\ktree)}.
    \end{equation}
    Treating the label as zeroing out more stump appended to \(\stree\) this gives the results of \cref{theory:pruningMarked} thus, reducing to
    \begin{equation}
        = \sum_{\utree \in \NTTi{\mathcal{A}}}  a(\utree \replaceNode{\mathcal{U}} \psi[\laddertree]) \frac{\BDiffOp(\utree)(\y_\initstate)}{\sigma(\utree)}.
    \end{equation}
    Thus, for linear trees the result holds.

    \item 
    \paragraph{Tall Bush Case:} \hfill \\
    For the tall bush trees \(\bushtree = \psi[\forestLadder]\) where \(\laddertree_{i} \ne \emptyset\,\forall\, i\in\{1..n\}\) then
    \begin{equation}
        \frac{\BDiffOp(\bushtree)(\y_\initstate)}{\sigma(\bushtree)} = \frac{\BDiffOp(\psi[\dforest])(\y_\initstate)}{\left(\prod_{j=1}^{n} q_j !\right)} \left( \prod_{j=1}^{n} \left(\BDiffOp(\laddertree_{j})(\y_\initstate)\right)^{q_j}\right)
    \end{equation}
    where \(\dforest = \stumpForestConvert(\forestLadder)\).
    Applying the ladder tree results to the children gives
    \begin{equation}
         = \frac{\BDiffOp(\psi[\dforest])(\y_\initstate)}{\left(\prod_{j=1}^{n} q_j !\right)} \left( \prod_{j=1}^{n} \left(\sum_{\ktree \in \NTTi{\mathcal{A}}} a(\ktree \replaceNode{\mathcal{U}} \laddertree_{j}) \frac{\BDiffOp(\ktree)(\y_\initstate)}{\sigma(\ktree)}\right)^{q_j}\right).
    \end{equation}
    Reducing this by applying \cref{lemma:permanentTrees} gives
    \begin{equation}
         = \frac{\BDiffOp(\psi[\dforest])(\y_\initstate)}{\left(\prod_{j=1}^{n} q_j !\right)} \left( \sum_{\gforest \in \ForestSpace{\NTTi{\mathcal{A}}}_{\bar{n}}} \frac{1}{\prod_{j=1}^{m} k_j!} \permanent\left(\gforest \, \outergen{a(\cdot_{1}  \replaceNode{\mathcal{U}}  \cdot_{2}) \frac{\BDiffOp(\cdot_1)(\y_\initstate)}{\sigma(\cdot_1)}} \, \fforest \right)\right)
    \end{equation}
    where \(\gforest = \forestTQs\). Next applying \cref{lemma:permanentscale} and pulling the label \(a\) out of the outer product and permanent gives 
    \begin{equation}
         = \frac{\BDiffOp(\psi[\dforest])(\y_\initstate)}{\left(\prod_{j=1}^{n} q_j !\right)} \left( \sum_{\gforest \in \ForestSpace{\NTTi{\mathcal{A}}}_{\bar{n}}} a(\permanent\left(\gforest \, \outergen{\replaceNode{\mathcal{U}}} \, \fforest \right)) \left(\prod_{j=1}^{m} \frac{\BDiffOp(\stree_{j})(\y_\initstate)^{k_j}}{\sigma(\stree)^{k_j} k_j!}\right)\right).
    \end{equation}
    Pushing the base factorial in then utilizing the forest definition of replacement gives 
    \begin{equation}
         = \BDiffOp(\psi[\dforest])(\y_\initstate)\left( \sum_{\gforest \in \ForestSpace{\NTTi{\mathcal{A}}}_{\bar{n}}} a\left(\gforest \, \replaceNode{\mathcal{U}} \, \fforest \right) \left(\prod_{j=1}^{m} \frac{\BDiffOp(\stree_{j})(\y_\initstate)^{k_j}}{\sigma(\stree)^{k_j} k_j!}\right)\right).
    \end{equation}
    Now for the case of when \(\psi \notin \mathcal{U}\). Since the forest replacement is zero when the number of children don't match then the set is replaced with the full bush tree set. Given the base node recursion condition then pushing the children elementary differential back on the stumped base gives 
    \begin{equation}
         = \sum_{\utree \in \NTTi{\mathcal{A}}} a(\utree \replaceNode{\mathcal{U}} \bushtree) \frac{\BDiffOp(\utree)(\y_\initstate)}{\sigma(\utree)}.
    \end{equation}
    For the case of when \(\psi \in \mathcal{U}\) then expanding the base node to its series expansion gives 
    \begin{equation}
        \stumpCountFactorial(\psi[\dforest])\left(\sum_{\rtree \in \NTTi{\mathcal{A}}\lceil\dforest\rceil} a(\labelbrack{\rtree}{\psi}) \frac{\BDiffOp(\rtree)(\y_\initstate)}{\sigma(\rtree)}\right)\left( \sum_{\gforest \in \ForestSpace{\NTTi{\mathcal{A}}}_{\bar{n}}} a\left(\gforest \, \replaceNode{\mathcal{U}} \, \fforest \right) \left(\prod_{j=1}^{m} \frac{\BDiffOp(\stree_{j})(\y_\initstate)^{k_j}}{\sigma(\stree)^{k_j} k_j!}\right)\right).
    \end{equation}
    Expanding the product and applying the tree pair relationship gives
    \begin{equation}
        \sum_{\rtree \in \NTTi{\mathcal{A}}\lceil\dforest\rceil} \sum_{\gforest \in \ForestSpace{\NTTi{\mathcal{A}}}_{\bar{n}}} \stumpCountFactorial(\psi[\dforest]) a(\labelbrack{\gforest}{\rtree} \replaceNode{\mathcal{U}} \labelbrack{\fforest}{\psi}) \frac{\BDiffOp(\rtree)(\y_\initstate)}{\sigma(\rtree)}\left(\prod_{j=1}^{m} \frac{\BDiffOp(\stree_{j})(\y_\initstate)^{k_j}}{\sigma(\stree)^{k_j} k_j!}\right).
    \end{equation}
    Treating the label as zeroing when more and less than \(\bar{n}\) stumps are appended to \(\rtree\) then this gives the results of \cref{theory:pruningMarked} thus, reducing to
    \begin{equation}
        \sum_{\utree \in \NTTi{\mathcal{A}}} a(\utree \replaceNode{\mathcal{U}} \bushtree) \frac{\BDiffOp(\utree)(\y_\initstate)}{\sigma(\utree)}.
    \end{equation}
    Thus, the results hold for arbitrary bush trees.
    \item 
    \paragraph{General Case:} \hfill \\
    For the general case taking the tree from \(\TallBushTreesN{\NTTC}{i}\). Assuming the result holds for \(\TallBushTreesN{\NTTC}{i-1}\) then applying all the same steps as with the tall bush tree gives the main result for \(\TallBushTreesN{\NTTC}{i}\). Thus, by induction the statement holds for arbitrary trees.
\end{proof}

\subsection{Runge-Kutta Time Expansion \cref{theorem:RKTime}} \label{proof:RKTime}
\begin{proof}
    Adding additional steps into the proof from  \cite[Theorem 3.6B]{Butcher_2016_book}. The differential equation \cref{eq:diffTimeSum} is Taylor expanded around \(\z_0\) gives
    \begin{equation}
        \begin{split}
            \fun\left(\y, \z_0 + \sum_{i=1}^{m} g_i(t)\,\basisK_i\right) &= \fun\left(\y, \z_0\right) + \sum_{i=1}^{m} g_i(t) \,\fun_{\z}\left(\y, \z_0\right)\basisK_i\\
            & + \frac{1}{2!} \sum_{i=1}^{m} g_i(t)^2\, \fun_{\z,\z}\left(\y, \z_0\right)(\basisK_i , \basisK_i)\\
            & + \sum_{i=2}^{m}\sum_{j=1}^{i-1} g_i(t)\,g_i(t)\, \fun_{\z,\z}\left(\y, \z_0\right)(\basisK_i , \basisK_j) + \dots.
        \end{split}
    \end{equation}
    For all evaluations of the RHS the derivative with respect to \(\z\) will all have \(k\) leaves of the time functions. Thus, encapsulated by the new node in the \BSeriesText differential operator. All the leading factorials in the expansion are pushed down to the elementary differentials into symmetry. Thus, giving the set of recursive labels 
    \begin{equation}
        \begin{split}
            \timeLabel{\theta}{t_\initstate}{t_\finalstate}^{\reductionFun}(\tree) &= 1 + A \left(\sum_{\hforest = \prod_{i}^{m} \lambda_i^{p_i}, p_i \ge 0} (\BSeriesDer{\tau \ButcherProduct \hforest} \timeLabel{\theta}{t_\initstate}{t_\finalstate}^{\reductionFun})(\tree) \right),\\
            \timeLabel{\Theta}{t_\initstate}{t_\finalstate}^{\reductionFun}(\tree) &= 1 + b^T \left(\sum_{\hforest = \prod_{i}^{m} \lambda_i^{p_i}, p_i \ge 0} (\BSeriesDer{\tau \ButcherProduct \hforest} \timeLabel{\theta}{t_\initstate}{t_\finalstate}^{\reductionFun})(\tree) \right).\\
        \end{split}
    \end{equation}
    Since the time dependent functions are only leaves then only one term of the sum is non-zero. This gives the updated recursion as
    \begin{equation}
        \scalemath{0.95}{
        \begin{split}
            \timeLabel{\theta}{t_\initstate}{t_\finalstate}^{\reductionFun}\left(\tau\left[\forestTQt \left(\prod_{i}^{m} \lambda_i^{p_i}\right)\right]\right) &= 1 + A\left(\left(\prod_{i=1}^{m} \timeLabel{\theta}{t_\initstate}{t_\finalstate}^{\reductionFun}(\lambda_i)\right)\left(\BSeriesDer{\tau} \timeLabel{\theta}{t_\initstate}{t_\finalstate}^{\reductionFun}\right)(\tau[\forestTQt])  \right),\\
            \timeLabel{\Theta}{t_\initstate}{t_\finalstate}^{\reductionFun}\left(\tau\left[\forestTQt \left(\prod_{i}^{m} \lambda_i^{p_i}\right)\right]\right) &= 1 + b^T\left(\left(\prod_{i=1}^{m} \timeLabel{\theta}{t_\initstate}{t_\finalstate}^{\reductionFun}(\lambda_i)\right)\left(\BSeriesDer{\tau} \timeLabel{\theta}{t_\initstate}{t_\finalstate}^{\reductionFun} \right)(\tau[\forestTQt]) \right).\\
        \end{split}
        }
    \end{equation}
    Applying the definition of \(\BSeriesDer{\tau}\) and of the time functions for \(\lambda_i\) nodes then \(\timeLabel{\theta}{t_\initstate}{t_\finalstate}^{\reductionFun} = \timeLabel{\phi}{t_\initstate}{t_\finalstate}^{\reductionFun}\) and \(\timeLabel{\Theta}{t_\initstate}{t_\finalstate}^{\reductionFun} = \timeLabel{\Phi}{t_\initstate}{t_\finalstate}^{\reductionFun}\).
\end{proof}

\subsection{Partial Exact Solution \cref{theorem:partialExact}} \label{proof:partialExact}
\begin{proof}
    For the case of trees of set \(\tree \in \NTTi{\Tcolors_1}\) then these are the individual methods labels. Thus, \(\varkappa\comp{k,a}(\tree) = \tree\) must hold for trees in the space. Since the recursion grafts back on to the root after proceeding back down, then the case holds.

    Next utilizing the same notation from \cref{corollary:BSeriesCompSuperSet}. Given the case of trees from \(\tree \in \NTTi{\Tcolors_1 < \Tcolors_2}\) then only trees of colors \(\Tcolors_2\) are traversed after leaving a subtree of colors \(\Tcolors_1\). Since the coupling terms are replaced with a time dependent exact solution, then only time derivatives are added by the node \(\btree{[t]}\). Thus, giving traversals over \(\beta\comp{k,a}\) by forming the product of all time terms to graft on after returning. To finalize the transformation a ratio of symmetries is taken since 
    \begin{equation}
         \frac{b}{\sigma}(\tree) = \frac{\timeLabel{a}{0}{1}\comp{0}}{\sigma}(\varkappa\comp{1,a}(\tree))
    \end{equation}
    must hold.

    Lastly for the case of traversing \(\tree \in (\NTTC \setminus \NTTi{\Tcolors_1 < \Tcolors_2})\). This implies the interaction with the interpolation on states. Thus, the recursive calling of \(\frac{\timeLabel{a}{0}{u}\comp{k}}{\sigma}(\varkappa\comp{k+1,a}(\tree))\) since the tree still needs to be transformed to represent the underlying system. Given the system solutions utilized the reduced nodes as time dependent terms then higher time derivatives of \(\timeLabel{a}{0}{u}\comp{k}\) still influence the same tree. Thus, the Taylor expansion w.r.t. the time dependent label must be taken giving
    \begin{equation}
        \sigma(\tree) \sum_{i=1}^{\infty} \btree{[t]}^{i} \frac{\ddiff^{i}}{\ddiff u^{i}} \left( \frac{\timeLabel{a}{0}{u}\comp{k}}{\sigma}(\varkappa\comp{k+1,a}(\tree))\right)\Big|_{u=0}
    \end{equation}
    where the \(\btree{[t]}^{i}\) represent the higher time actions and the one over factorial terms are introduced at the base of the recursion by \(\sigma\). Given that this recursively calls the higher interpolations, then by induction this holds for all trees.
\end{proof}

\subsection{TII Expansion \cref{theory:ittLabelExact}}\label{proof:ittLabelExact}
\begin{proof}
    By inserting the coefficients \(\frac{1}{\gammaFun{\Tcolors \setminus \mathcal{V}}}\) in the underlying labels by \cref{theorem:reductionForceRepresentation} then the analytical underlying equation is obtained. Next replacing all terms with the integral representation from \cref{theory:iterlabelgamma} gives
    \begin{equation}
        \sigma(\stree) \sum_{\rtree \in \NTTC}  \frac{\Dt^{\treeorder{\rtree} - \treeorder{\stree}}}{\sigma(\rtree)} \timeElemOp{t_\initstate}^{\reductionFun}(\rtree \replaceNode{\mathcal{U}} \stree) \treeIterInt{t_\initstate}{t_\finalstate}^{\delta_{\cdot, \Tcolors \setminus \mathcal{V}}}(\rtree).
    \end{equation}
    To proceed, introducing the nodes \(\vnode_{i,j} \in \mathcal{V}^*\) by replacing nodes \(\unode_{i} \in \mathcal{U}\) where it corresponds to the \(i\)th node in \(\vnode_{i} \in \mathcal{V}\). The second index of \(\vnode_{i,j}\) corresponds to its amount of time derivatives extending from it. Given the product definition of TII \cref{def:treeIterativeIntegral} and removing \(j\) time nodes extending from \(\vnode_{i,j}\) gives a new definition \(\varrho(\vnode_{i,j})(\tau) = \frac{(\Dt \, \tau)^{j}}{j!}\) thus giving
    \begin{equation}
        \sigma(\stree) \sum_{\rtree \in \NTTC}  \frac{1}{\sigma(\rtree)} \timeElemOp{t_\initstate}^{\reductionFun}(\rtree \replaceNode{\mathcal{U}} \stree) \treeIterInt{t_\initstate}{t_\finalstate}^{\varrho}(\rtree).
    \end{equation}
    Next by taking the leading elementary time derivatives terms \(\timeElemOp{t_\initstate}^{\reductionFun}(\rtree \replaceNode{\mathcal{U}} \stree)\) and inserting them into their respective \(\vnode_{i,j}\) node locations for \(j\) derivative by \cref{property:scalaritt} gives \(\lambda(\vnode_{i,j})(\tau) = \reductionFun(\vnode_{i})^{(j)}(t_\initstate)\frac{(\Dt \, \tau)^{j}}{j!}\) as 
    \begin{equation}
        \sigma(\stree) \sum_{\rtree \in \NTTC}  \frac{1}{\sigma(\rtree)} \treeIterInt{t_\initstate}{t_\finalstate}^{\lambda}(\rtree).
    \end{equation}
    Lastly, utilizing \cref{property:linitt} to sum all terms of \(j\) to \(\vnode_i\) gives \(\reductionFun(\vnode_{i}) = \sum_{j=0}^{\infty} \reductionFun(\vnode_{i})^{(j)}(t_\initstate)\frac{(\Dt \, \tau)^{j}}{j!}\) the Taylor series of time dependent coefficient thus, matching the TII expansion.
\end{proof}

\subsection{Integral Analytical Composition \cref{theory:analyticalComposition}}\label{proof:analyticalComposition}
\begin{proof}
    We first prove recursively for the case of when \(\utree = \emptyset\) then we extend to the general case.
    \item 
    \paragraph{Tall Tree Case:} \hfill \\
    Let \(\laddertree_{n} \in \NTTC\) be some ladder tree such that \(\treeorder{\laddertree} = n \ge 1\). Applying \cref{property:domainSplitTII} to the base integral gives
    \begin{equation}
        \treeIterInt{a}{b}^{\reductionFun}(\laddertree_{n}) + \int_{b}^{c} \reductionFun(\psi_{1})(\tau) \treeIterInt{a}{\tau}(\laddertree_{n-1}) d\tau,
    \end{equation}
    and by repeating the application gives
    \begin{equation}
        \treeIterInt{a}{b}^{\reductionFun}(\laddertree_{n}) + \int_{b}^{c}\reductionFun(\psi_{1})(\tau_1) \left(\treeIterInt{a}{b}^{\reductionFun}(\laddertree_{n-1}) + \int_{b}^{\tau_1} \reductionFun(\psi_{2})(\tau_2) \treeIterInt{a}{\tau_2}^{\reductionFun}(\laddertree_{n-2}) \dint \tau_2 \right)\dint \tau_1,
    \end{equation}
    where \(\treeIterInt{a}{b}^{\reductionFun}(\laddertree_i)\) is constant thus can be flushed through the integral. Furthermore, the term \(\int_{b}^c \reductionFun(\psi)(\tau_1) \dint \tau_1 = \treeIterInt{b}{c}^{\reductionFun}(\laddertree_{1})\) is simplified thus 
    \begin{equation}
        \treeIterInt{a}{b}^{\reductionFun}(\laddertree_{n}) + \treeIterInt{a}{b}^{\reductionFun}(\laddertree_{n-1})\treeIterInt{b}{c}^{\reductionFun}(\laddertree_{1}) + \int_{b}^{c} \reductionFun(\psi_{1})(\tau_1) \int_{b}^{\tau_1} \reductionFun(\psi_{2})(\tau_2) \treeIterInt{a}{\tau_2}(l_{n-2}) \dint \tau_2 \dint \tau_1.
    \end{equation}
    Inductively as the integrals proceed up then it will always prune continuing up the ladder tree. Once the cuts reach the top of the ladder \(n - 1\) then the case of a single node is reached
    \begin{equation}
        + \int_{b}^{c} \reductionFun(\psi_{1})(\tau_1) \int_{b}^{\tau_1} \reductionFun(\psi_{2})(\tau_2) \dots \int_{b}^{\tau_{n - 1}}  \reductionFun(\psi_{n})(\tau_{n-1})\left(\treeIterInt{a}{b}^{\reductionFun}(\laddertree_{1}) + \treeIterInt{b}{\tau_{n-1}}^{\reductionFun}(\laddertree_{1})   \right) \dint \tau_{n - 1} \dots \dint \tau_2 \dint \tau_1.
    \end{equation}
    Thus, giving the result of the hypothesis for the case of ladder trees as
    \begin{equation}
        \treeIterInt{a}{c}^{\reductionFun}(\laddertree_{n}) = \sum_{\stree \in \LadderTrees{\NTTC}} \treeIterInt{a}{b}^{\reductionFun}(\laddertree_{n} \prune \stree) \treeIterInt{b}{c}^{\reductionFun}(\stree).
    \end{equation}
    Extending to the case of when \(\utree \ne \emptyset\) and is a ladder tree then given \(\utree = \laddertree_{m}\) gives the result
    \begin{equation}
        \begin{split}
            \treeIterInt{a}{c}^{\reductionFun}(\laddertree_{n} \prune \laddertree_{m}) &= \treeIterInt{a}{c}^{\reductionFun}(\laddertree_{n - m}) = \sum_{\stree \in \LadderTrees{\NTTC}} \treeIterInt{a}{b}^{\reductionFun}(\laddertree_{n - m} \prune \stree)\treeIterInt{b}{c}^{\reductionFun}(\stree)\\
            & = \sum_{\stree \in \LadderTrees{\NTTC}} \treeIterInt{a}{b}^{\reductionFun}(\laddertree_{n} \prune \stree)\treeIterInt{b}{c}^{\reductionFun}(\stree \prune \laddertree_{m}).\\
        \end{split}
    \end{equation}

    \item 
    \paragraph{Tall Bush Case:} \hfill \\
    To start focusing on the case of \(\utree = \emptyset\) given the tall bush tree as \(\bushtree_{0} = \psi[\fforest] = \psi[\laddertree_1^{q_1} \laddertree_2^{q_2} \dots \laddertree_n^{q_n}] \in \TallBushTrees{\NTTC}\). Starting with \cref{property:domainSplitTII}.
    \begin{equation}
        \treeIterInt{a}{c}^{\reductionFun}(\bushtree_{0}) = \treeIterInt{a}{b}^{\reductionFun}(\bushtree_{0}) + \int_{b}^{c} \reductionFun(\psi)(\tau_1) \prod_{i=1}^{n} \treeIterInt{a}{\tau_1}^{\reductionFun}(\laddertree_i)^{q_i} \dint \tau_1
    \end{equation}
    Since these are ladder trees, then reinserting in the result gives
    \begin{equation}
        \treeIterInt{a}{b}^{\reductionFun}(\bushtree_{0}) + \int_{b}^{c} \reductionFun(\psi)(\tau_1) \prod_{i=1}^{\bar{n}} \left(\sum_{\stree_i \in \LadderTrees{\NTTC}} \treeIterInt{a}{b}^{\reductionFun}(\laddertree_i \prune \stree_i)\treeIterInt{b}{\tau_1}^{\reductionFun}(\stree_i) \right) \dint \tau_1
    \end{equation}
    Given \cref{lemma:permanentTrees} then this can be converted to the permanent form as
    \begin{equation}
        \treeIterInt{a}{b}^{\reductionFun}(\bushtree_{0}) + \int_{b}^{c} \reductionFun(\psi)(\tau_1) \sum_{\gforest \in \ForestSpace{\LadderTrees{\NTTC}}} \frac{1}{(\prod_{j=1}^{m} k_j!)p!} \permanent(\fforest \outergen{\treeIterInt{a}{b}^{\reductionFun}(\cdot_1 \prune \cdot_2)\treeIterInt{b}{\tau_1}^{\reductionFun}(\cdot_2)} \gforest) \dint \tau_1
    \end{equation}
    where \(p = \sum_{i=1}^{n}q_i - \sum_{j=1}^{m}k_j \ge 0\).
    Since the columns all contain the same scalar \(\treeIterInt{b}{\tau_1}^{\reductionFun}(\cdot_2)\), then by applying \cref{lemma:permanentscale} pulls the scalar out of the expression. Furthermore, since the label \(\treeIterInt{a}{b}^{\reductionFun}(\cdot_1 \prune \cdot_2)\) is elementwise, then the function can surround the permanent result over pruning. Thus, since the pruning terms are constant, then they are pulled out of the base integral giving 
    \begin{equation}
        \treeIterInt{a}{b}^{\reductionFun}(\bushtree_{0})\treeIterInt{b}{c}^{\reductionFun}(\emptyset) + \sum_{\gforest \in \ForestSpace{\LadderTrees{\NTTC}}}  \treeIterInt{a}{b}^{\reductionFun}\left(\frac{1}{(\prod_{j=1}^{m} k_j!)p!} \permanent(\fforest \outergen{\prune} \gforest)\right)  \left(\int_{b}^{c} \reductionFun(\psi)(\tau_1) \prod_{j=1}^{m}\treeIterInt{b}{\tau_1}^{\reductionFun}(s_j)^{k_j} \dint \tau_1\right).
    \end{equation}
    Thus, by the definition of pruning (\cref{theory:pruning}) and TII (\cref{def:treeIterativeIntegral}) then this reduces to the result of 
    \begin{equation}
        \treeIterInt{a}{c}^{\reductionFun}(\bushtree_{0}) = \sum_{\stree \in \TallBushTrees{\NTTC}} \treeIterInt{a}{b}^{\reductionFun}(\bushtree_{0} \prune \stree) \treeIterInt{b}{c}^{\reductionFun}(\stree).
    \end{equation}
    Next for the case of when \(\utree = \tilde{\bushtree}_{0} = \psi[\gforest] = \psi[\tilde{\laddertree}_1^{k_1} \tilde{\laddertree}_2^{k_2} \dots \tilde{\laddertree}_m^{k_m}]\) then from the permanent co-factor expansion
    \begin{equation}
        \treeIterInt{a}{c}^{\reductionFun}(\bushtree_{0} \prune \tilde{\bushtree}_{0}) = \sum_{\sigma \in S_{\gforest}} \prod_{j=1}^{\bar{n}} \treeIterInt{a}{c}^{\reductionFun}(\fforest_{\pi(j)} \prune \gforest_{\sigma(j)})
    \end{equation}
    for some arbitrary integer mapping \(\pi\).
    Since all subtrees have had the result proven, then
    \begin{equation}
        = \sum_{\sigma \in S_{\gforest}} \prod_{j=1}^{\bar{n}} \left[\sum_{\stree \in \LadderTrees{\NTTC}} \treeIterInt{a}{b}^{\reductionFun}(\fforest_{\pi(j)} \prune \stree) \treeIterInt{b}{c}^{\reductionFun}(\stree \prune \gforest_{\sigma(j)})\right]
    \end{equation}
    Next applying \cref{lemma:permanentTrees} to push the product up and convert the terms to permanents in to sums over \(\hforest\) gives
    \begin{equation}
        \begin{split}
            &= \sum_{\sigma \in S_{\gforest}} \sum_{\hforest \in \ForestSpace{\LadderTrees{\NTTC}}} \sum_{\theta \in S_{\hforest}} \prod_{j=1}^{\bar{n}} \treeIterInt{a}{b}^{\reductionFun}(\fforest_{\pi(j)} \prune \hforest_{\theta(j)}) \treeIterInt{b}{c}^{\reductionFun}(\hforest_{\theta(j)} \prune \gforest_{\sigma(j)})\\
            &= \sum_{\hforest \in \ForestSpace{\LadderTrees{\NTTC}}} \frac{1}{(\prod_{i=1}^{m} x_i!)z!} \permanent(\fforest \outergen{\treeIterInt{a}{b}^{\reductionFun}(\cdot \prune \cdot)} \hforest) \left[ \sum_{\theta \in S_{\hforest}} \sum_{\sigma \in S_{\gforest}} \prod_{j=1}^{\bar{n}} \treeIterInt{b}{c}^{\reductionFun}(\hforest_{\theta(j)} \prune \gforest_{\sigma(j)})\right]\\
            &= \sum_{\hforest \in \ForestSpace{\LadderTrees{\NTTC}}} \frac{1}{(\prod_{i=1}^{m} x_i!)z!} \permanent(\fforest \outergen{\treeIterInt{a}{b}^{\reductionFun}(\cdot \prune \cdot)} \hforest) \frac{1}{(\prod_{i=1}^{m} q_i!)p!} \permanent(\hforest \outergen{\treeIterInt{b}{c}^{\reductionFun}(\cdot \prune \cdot)} \gforest).\\
        \end{split}
    \end{equation}
    Just as before pulling the labels and utilizing the definition of pruning gives the result 
    \begin{equation}
        = \sum_{\stree \in  \TallBushTrees{\NTTC}} \treeIterInt{a}{b}^{\reductionFun}(\bushtree_{0} \prune \stree) \treeIterInt{b}{c}^{\reductionFun}(\stree \prune \utree).
    \end{equation}

    \item 
    \paragraph{General Case:} \hfill \\
    With the base case of the tall bush tree given and inductively assuming the result is true for the children of \(\tree \in \TallBushTreesN{\NTTC}{n}\) then all the same steps as before are applied in the case of \(\utree \in \TallBushTreesN{\NTTC}{n}\) and \(\utree = \emptyset\). Thus, by induction the overall result holds.
\end{proof}

\end{document}